\documentclass{amsart}

\usepackage[all]{xy}

\usepackage{latexsym}
\usepackage{amssymb}
\usepackage{amsthm}

\usepackage{mathrsfs}
\let\mathcal\mathscr

\theoremstyle{plain}
\newtheorem{prop}{Proposition}
\newtheorem{theo}[prop]{Theorem}
\newtheorem{lemm}[prop]{Lemma}
\newtheorem{coro}[prop]{Corollary}
\theoremstyle{definition}
\newtheorem*{defi}{Definition}
\newtheorem*{rema}{\sc Remark}

\newtheorem*{exams}{Examples}
\newtheorem*{remas}{\sc Remarks}

\def\gr{\textrm{gr}}
\newcommand{\CC}{\mathbb{C}}
\def\disc{{disc}}

\def\hdb#1#2#3{\left\la #3\right\ra_{#2}}
\def\hdbd#1#2#3{\left\la #3\right\ra^{#1}_{#2}}

\def\TTT{{\mathcal{T}}}
\newcommand{\Tw}{\mathrm{Tw}}
\newcommand{\la}{\langle}
\newcommand{\ra}{\rangle}
\newcommand{\fD}{fD}
\newcommand{\pa}{\partial}
\newcommand{\qi}{\xrightarrow{\sim}}
\newcommand{\g}{\mathfrak{g}}
\newcommand{\BV}{\mathcal{B}\mathcal{V}}
\newcommand{\GV}{\mathcal{G}\mathcal{V}}

\newcommand{\qBV}{\mathrm{q}\BV}
\newcommand{\qR}{\mathrm{q}R}
\newcommand{\qPo}{\mathrm{q}\Po}

\newcommand{\KK}{\mathbb{K}}
\newcommand{\ZZ}{\mathbb{Z}}

\newcommand{\End}{\mathrm{End}}

\newcommand{\G}{\mathcal{G}}

\newcommand{\NN}{\mathbb{N}}

\newcommand{\Sy}{\mathbb{S}}

\newcommand{\Com}{\mathcal{C}om}

\newcommand{\Po}{\mathcal{P}}
\newcommand{\F}{\mathcal{F}}

\newcommand{\ac}{\scriptstyle \text{\rm !`}}

\newcommand{\Id}{\mathrm{Id}}
\newcommand{\Qo}{\mathcal{Q}}

\newcommand{\Co}{\mathcal{C}}

\newcommand{\ot}{\otimes}

\newcommand{\Hom}{\mathrm{Hom}}

\newcommand{\Coder}{\mathrm{Coder}}
\newcommand{\Codiff}{\mathrm{Codiff}}
\newcommand{\sgn}{{sgn}}

\newcommand{\epi}{\twoheadrightarrow}
\newcommand{\mono}{\rightarrowtail}

\def\Lie{\mathit{{\mathcal L}\!ie}}
\def\Com{\mathit{{\mathcal C}\!om}}
\def\Ass{\mathit{{\mathcal A}\!ss}}

\newcommand{\Y}{\vcenter{\xymatrix@M=0pt@R=6pt@C=6pt{
\ar@{-}[dr] &  &\ar@{-}[dl]  \\
 &\ar@{-}[d] &  \\  & &}}}
\newcommand{\YY}{\vcenter{\xymatrix@M=0pt@R=6pt@C=6pt{
\ar@{-}[dr] &  &\ar@2{-}[dl]  \\
 &\ar@2{-}[d] &  \\  & &}}}
\newcommand{\YYY}{\vcenter{\xymatrix@M=0pt@R=6pt@C=6pt{
\ar@{-}[dr] &  &\ar@3{-}[dl]  \\
 &\ar@3{-}[d] &  \\  & &}}}
\newcommand{\cop}{\vcenter{\xymatrix@M=0pt@R=6pt@C=6pt{
 & \ar@{-}[d] & \\
 &\ar@{-}[dr] \ar@{-}[dl] &  \\  & &}}}

\newcommand{\copL}{\xymatrix@M=0pt@R=6pt@C=6pt{
 & \ar@{-}[d] & \\
 &\ar@{-}[dr] \ar@{-}[dl] &  \\  & &\\  & &\\  & &}}
\newcommand{\YL}{\vcenter{\xymatrix@M=0pt@R=6pt@C=6pt{
\ar@{-}[dr] &  &\ar@{-}[dl]  \\
 &\ar@{-}[d] &   \\  & &\\  & &\\  & &}}}
\newcommand{\YYL}{\vcenter{\xymatrix@M=0pt@R=6pt@C=6pt{
\ar@{-}[dr] &  &\ar@2{-}[dl]  \\
 &\ar@2{-}[d] &  \\  & &\\  & &\\  & &}}}
\newcommand{\LYY}{\vcenter{\xymatrix@M=0pt@R=6pt@C=6pt{
\ar@{-}[dr] &  &\ar@2{-}[dl]  \\
 &\ar@2{-}[d] &  \\  & &\\  & &\\  & &}}}

\newcommand{\XX}{\vcenter{\xymatrix@M=0pt@R=6pt@C=6pt{\ar@{-}[ddrr]&&\ar@{-}[ddll] \\ && \\ &&   }}}

\newcommand{\Ta}{\vcenter{\xymatrix@M=0pt@R=6pt@C=6pt{ \ar@{-}[dddrrr] && \ar@{-}[dl] &&  \\
&&& \ar@{-}[dl]  &  \\ &&&&  \ar@{-}[dl]  \\&&&  \ar@{-}[d] &
\\&&&& }}}
\newcommand{\Tb}{\vcenter{\xymatrix@M=0pt@R=6pt@C=6pt{  & \ar@{-}[dr]&&\ar@{-}[dl] \\
\ar@{-}[dr]&&\ar@{-}[dl]& \\&\ar@{-}[dr]&&\ar@{-}[dl]
\\&&\ar@{-}[d]& \\&&& }}}
\newcommand{\Tc}{\vcenter{\xymatrix@M=0pt@R=6pt@C=6pt{   \ar@{-}[dr]&&\ar@{-}[dl]& \\
&\ar@{-}[dr]&& \ar@{-}[dl] \\\ar@{-}[dr]&&\ar@{-}[dl]&
\\&\ar@{-}[d]&& \\&&& }}}
\newcommand{\Td}{\vcenter{\xymatrix@M=0pt@R=6pt@C=6pt{ && \ar@{-}[dr]&&\ar@{-}[dddlll] \\
 &\ar@{-}[dr]&&& \\ \ar@{-}[dr]&&&& \\& \ar@{-}[d]&&& \\&&&&  }}}
\newcommand{\Te}{\vcenter{\xymatrix@R=3pt@C=3pt{\ar@{-}[drdr] &&\ar@{-}[dl]  *=0{}
\ar@{-}[dr]&& \ar@{-}[ddll] \\ &&& *=0{}& \\&& *=0{} \ar@{-}[d]&&
\\&&&& }}}
\newcommand{\TaC}{\vcenter{\xymatrix@M=0pt@R=6pt@C=6pt{ \ar@{-}[ddddddrrrrrr] && \ar@{-}[dl] && && \\
&&& \ar@{-}[dl]  &&&  \\ &&&&  \ar@{-}[dl]&&  \\&&& &&&
\\&&&&\ar@{-}[dl]&& \\&&&&&\ar@{-}[dl]&\\&&&&&& }}}
\newcommand{\TreeL}{\vcenter{\xymatrix@M=0pt@R=5pt@C=5pt{ \ar@{-}[dr] &
&\ar@{-}[dl] & &  \\
& \ar@{-}[dr] & &\ar@{-}[dl]  & \\
& &\ar@{-}[d] & & \\
& & \\ & & }}}
\newcommand{\TreeR}{\vcenter{\xymatrix@M=0pt@R=5pt@C=5pt{
 & &\ar@{-}[dr] & & \ar@{-}[dl]  \\
& \ar@{-}[dr] & &\ar@{-}[dl]  & \\
& &\ar@{-}[d] & & \\
& & \\ & & }}}



\thanks{The first author was partially supported by grants 
MTM2007-63277, MTM2010-15831, MTM2010-20692, SGR1092-2009,  
and the second by
MTM2007-63277, MTM2010-15831, SGR119-2009.
The third author is supported by ANR grant JCJC06 OBTH}

\title{Homotopy Batalin--Vilkovisky algebras}

\author{Imma G\'alvez-Carrillo}
\address{Departament de Matem\`atiques, Edifici C
       \\Universitat Aut\`onoma de Barcelona
       \\08193 Bellaterra (Barcelona), Spain
\newline{\emph{From February 1, 2011:}} 
        Departament de Matem\`atica Aplicada III
      \\Universitat Polit\`ecnica de Catalunya
      \\Escola d'Enginyeria de Terrassa 
      \\Carrer Colom 1\\08222 Terrassa (Barcelona)\\Spain
}
\email{igalvez@mat.uab.cat}

\author{Andy Tonks}
\address{STORM, London Metropolitan University\\
166--220 Holloway Road, London N7 8DB, UK}
\email{a.tonks@londonmet.ac.uk}

\author{Bruno Vallette}
\address{Laboratoire J.A. Dieudonn\'e\\Universit\'e de Nice Sophia-Antipolis\\
Parc Valrose\\06108 Nice\\ Cedex 02\\France}
\email{brunov@unice.fr}

\begin{document}
\maketitle
\begin{abstract}
This paper provides an explicit cofibrant resolution of the operad encoding Batalin--Vilkovisky algebras. Thus it defines the notion of \emph{homotopy Batalin--Vilkovisky algebras} with the required homotopy properties.

To define this resolution, we extend the theory of Koszul duality to operads and properads that are defined by quadratic
and linear relations. The operad encoding Batalin--Vilkovisky
algebras is shown to be Koszul in this sense. This allows us to prove a Poincar\'e--Birkhoff--Witt Theorem for such an operad and to give an explicit small quasi-free resolution for it.

This particular resolution enables us to describe the deformation theory and homotopy theory of BV-algebras and of homotopy BV-algebras. We show that any topological conformal field theory  carries a homotopy BV-algebra structure which lifts the BV-algebra structure on homology. The same result is proved for the singular chain complex of the double loop space of a topological space endowed with an action of the circle. We also prove the cyclic Deligne conjecture with this cofibrant resolution of the operad $\BV$. We develop the general obstruction theory for algebras over the Koszul resolution of a properad and apply it to extend a conjecture of Lian--Zuckerman, showing that certain vertex algebras have an explicit homotopy BV-algebra structure.
\end{abstract}


  \tableofcontents


\section*{Introduction}

The main goal of this paper is to develop the homotopy theory of Batalin--Vilkovisky algebras, BV-algebras for short, and to apply it to algebra, topology, geometry and mathematical physics. 

The ``disadvantage with monoids is that they do not live in homotopy theory'' said Saunders MacLane in 1967.
Given two homotopy equivalent topological spaces or chain complexes, such that one of
them is endowed with a monoid or an associative algebra structure, we may transfer the binary product to the other in an obvious way. But this transferred product will not be associative in general; it carries higher homotopies and the resulting structure is a \emph{homotopy associative}, or simply \emph{$A_\infty$}, space or algebra~\cite{Stasheff63}. The category of $A_\infty$-algebras includes the category of associative algebras and is stable under homotopy constructions. To this end, homotopy and algebra meet.

To control the combinatorics of these higher homotopies, MacLane introduced the notion of \emph{prop} \cite{MacLane65}
of which \emph{operads} are a particular example \cite{May72}. Operads and props are used to encode algebraically the operations acting in certain categories. Props model multilinear operations with many inputs and many outputs, whereas operads model multilinear operations with many inputs but only one output. One can do homological algebra and homotopy theory on the level of operads themselves, and in this context there is a notion of \emph{cofibrant operad} that is to operads what projective modules are to modules. It is proved that categories of algebras over cofibrant operads enjoy good homotopy properties \cite{BoardmanVogt73, BergerMoerdijk03}. For example, the operad which encodes $A_\infty$-algebras is cofibrant.

The notion of Batalin--Vilkovisky algebra plays an important role in geometry, topology and mathematical
physics. Unfortunately, as with associative algebras, the operad which governs them is not
cofibrant. Therefore the first question answered by this paper is that of providing an explicit
cofibrant resolution of this operad, algebras over which are naturally
called \emph{homotopy Batalin--Vilkovisky algebras}. This conceptual definition ensures that these algebras share nice homotopy properties.

To define this resolution we extend to operads the inhomogeneous Koszul duality theory for algebras, after \cite{Priddy70,GinzburgKapranov94, GetzlerJones94}, and apply it to the operad $\BV$ which encodes BV-algebras. The resulting explicit cofibrant resolution allows us to give four equivalent definitions of a homotopy BV-algebra, each of which we make explicit for the different applications in the text. Our approach via Koszul duality theory gives as a corollary an algebraic theorem: the Poincar\'e--Birkhoff--Witt theorem for the operad $\BV$ itself.

Another natural question is that of lifting to the chain complex level an algebraic structure (such as a BV-algebra) that is given on the homology. The above arguments show that one could hope for a homotopy BV-algebra structure on the chain complex in general. Since we define the notion of homotopy BV-algebra by a cofibrant operad, we can prove the following lifting  results in mathematical physics and in algebraic topology: any topological conformal field theory and any double loop space of a topological space endowed with an action of the circle carry a homotopy BV-algebra structure which lifts the BV-algebra structure of \cite{Getzler94} on homology. We also show that the cyclic Deligne Conjecture is true with the homotopy $\BV$ operad provided here. These proofs rely on the formality of the framed little discs operad \cite{Severa09, GiansiracusaSalvatore09} and the fact that the homotopy $\BV$ operad is cofibrant. 

Since the resolution we provide here for the operad $\BV$ comes from an extended Koszul
duality theory, we can make the deformation theory of homotopy BV-algebras explicit.
We describe a dg Lie algebra whose Maurer--Cartan elements are in one-to-one correspondence
with homotopy BV-algebra structures. Given such an element, the associated twisted dg Lie
algebra defines the cohomology of the homotopy BV-algebra. This also allows us to study the
obstruction theory for homotopy BV-algebras. In the relative case, we make explicit the
obstructions to lift a homotopy Gerstenhaber algebra structure to a homotopy BV-algebra structure.
We apply the general obstruction theory to prove an extended version of Lian--Zuckerman conjecture:
any  topological vertex operator algebra, with $\NN$-graded conformal weight, admits a homotopy BV-algebra
structure which extends Lian--Zuckerman operations  and  which lifts  the BV-algebra structure on the BRST
homology. This homotopy BV-algebra structure is \emph{explicit}, unlike in the
previous applications, since the method used is different. 

Using the Koszul dual cooperad of the operad encoding BV-algebras, we can develop several
important constructions in homotopy theory for BV-algebras and homotopy BV-algebras:
bar and cobar constructions, $\infty$-morphisms, transfer of homotopy BV-algebras structures
under chain homotopy equivalences and Massey products. This allows us to formulate the formality
conjecture for the Hochschild cochains in the BV case, extending Kontsevich's result
\cite{Kontsevich03}. Other applications of homotopy BV-algebras are expected in
String Topology \cite{ChasSullivan99}, Frobenius manifolds \cite{Manin99},
Chiral algebras \cite{BeilinsonDrinfeld04} and in the Geometric Langlands Programme \cite{FrenkelBenZvi05}.

In \cite{Vallette07}, we introduced the notion of \emph{properad} which lies between the notions of operad and prop. Properads faithfully encode categories of bialgebras with products and coproducts, as do props, but they are simpler. Hence it was possible to develop a Koszul duality theory on that level, which cannot be done for props since they lack, for the moment, the required homological constructions, such as the bar construction.

The second part of this paper is composed of three appendices which introduce
the general theories for properads that  we have applied to
the operad $\BV$  in the first part. The first appendix describes the Koszul duality theory for
inhomogeneous quadratic properads, and thus operads; as a direct corollary we
prove the Poincar\'e--Birkhoff--Witt Theorem for Koszul properads.
In the second, we develop the homotopy theory for algebras over such a Koszul
operad (bar and cobar constructions, homotopy algebras, $\infty$-morphisms,
transfer of structure and Massey products). The third appendix describes the
deformation and obstruction theory of algebras over a Koszul properad.

\section*{Operadic homological algebra}

In this section, we recall the basic notions of homological algebra for operads such as twisting morphisms and bar and cobar constructions. We give the four equivalent definitions of homotopy algebras over a Koszul operad that we will use throughout the paper, one for each section.

\smallskip

We refer the reader  to the papers
\cite{GinzburgKapranov94,GetzlerJones94, Fresse04} and to the books \cite{MarklShniderStasheff02, LodayVallette09}
for further details.

\subsection{Conventions}

Throughout the text, we work over the underlying category of differential graded modules, or chain complexes, over a field $\KK$ of characteristic $0$. Because of the relationship with algebraic topology, we use the homological degree convention:  the differential maps lower the degree by $1$. 
Hence, we work with dg operads, dg cooperads, dg algebras and dg coalgebras, which will often be called operads, cooperads, algebras and coalgebras. Recall that a dg operad is a monoid in the monoidal category of dg $\Sy$-modules equipped with the composite product $\circ$. A dg cooperad  is a comonoid in this monoidal category. 
Notice that most of the operads and cooperads present in this paper are reduced, that is, they satisfy $\Po(0)=0$.
The degree of an element or map is denoted by $|x|$. We adopt the usual opposite grading for the linear dual of a chain complex, $(V^*)_n:=(V_{-n})^*$. Let $s\KK$ denote a copy of the field concentrated in degree 1 and let
$s\,:\, V\to sV$ denote the suspension operator, where $(sV)_n=s\KK\otimes
V_{n-1}$ and $s(v)=1\otimes v$.

\subsection{Twisting morphisms}\label{Twisting morphisms}

Let $\Co$ be a dg cooperad and let $\Qo$ be a dg operad.
Recall from \cite{BergerMoerdijk03} that the collection $\Hom(\Co, \Qo):=\{ \Hom(\Co(n), \Qo(n))\}_{n\in \NN}$ forms an operad called the \emph{convolution operad}. This structure induces a dg preLie algebra structure and hence a dg Lie algebra structure on the direct product of the $\Sy$-equivariant maps
$$\Hom_\Sy(\Co, \Qo):=\big(\prod_{n\in \NN} \Hom_{\Sy_n}({\Co}(n), \Qo(n)), \partial, [\;,\,] \big), $$
    see \cite{KapranovManin01}.

In this \emph{convolution dg Lie algebra}, we consider the Maurer--Cartan equation $\partial (\alpha) + \frac{1}{2}[\alpha, \alpha]=0$, whose degree $-1$ solutions are called \emph{twisting morphisms} by analogy with algebraic topology. The set of twisting morphisms is denoted $\Tw(\Co, \Qo)$.\\

A cooperad $\Co$ is called \emph{coaugmented} when the counit map splits by a morphism of cooperads $I \to \Co$ called the \emph{coaugmentation map}, with cokernel denoted by $\overline\Co$. When $\Co$ is coaugmented, we require that twisting morphisms factor through $\overline \Co\to \Qo$.

\subsection{Bar and cobar constructions}

The twisting morphisms functor,
$\Tw \, :\, \textrm{coaugm. dg coop.} \times \textrm{dg op.} \to \textrm{sets}$, is represented by the following functors:
$$\Omega \ : \ \{\textrm{coaugmented dg cooperads}\} \rightleftharpoons
\{\textrm{dg operads}\}\ : \ \mathrm{B}$$
called the \emph{bar construction} $\mathrm{B}$ and the \emph{cobar construction} $\Omega$. \\

The cobar construction $\Omega \Co$ of a coaugmented dg cooperad $(\Co, d_\Co)$ is defined by the free operad $\F(s^{-1}\overline \Co)$ on the homological desuspension of $\overline \Co$. It is endowed with the derivation $d=d_2+d_1$ given by the unique  derivation $d_2$ which extends the partial decomposition map of the cooperad $\Co$ plus the unique derivation $d_1$ which extends the differential $d_\Co$ of $\Co$. The bar construction $\mathrm{B} \Qo$ is defined dually by the cofree cooperad $\F^c(s{\Qo})$, see \cite[Section~$2$]{GetzlerJones94} and \cite[Chapter~$5$]{LodayVallette09}.\\

A \emph{conilpotent cooperad} $(\Co, \Delta)$ is a coaugmented cooperad such that the iterated powers of the decomposition map $\Delta$ minus its primitive part applied to any element eventually vanish. The bar and cobar constructions form an adjunction between dg operads and conilpotent dg cooperads, which is explicitly given by the set of twisting morphisms \cite[Theorem~$2.17$]{GetzlerJones94}
$$\Hom_{\textrm{dg op.}}\left(\Omega \Co, \Qo\right) \cong
\Tw(\Co, \Qo) \cong  \Hom_{\textrm{dg coaugm.
coop.}}\left(\Co, \mathrm{B} \Qo\right).$$

\subsection{Homotopy algebras}

When an operad $\Po$ is Koszul, the cobar construction of its Koszul dual cooperad $\Po^{\ac}$ provides a quasi-free resolution of $\Po$, see Appendix~\ref{Appendix Koszul}.
$$\Po_\infty:=\Omega \Po^{\ac} \xrightarrow{\sim} \Po$$
We define \emph{homotopy $\Po$-algebras} to be the algebras over the cofibrant operad $\Po_\infty$; thus they share nice homotopy properties \cite{BoardmanVogt73,BergerMoerdijk03}. A homotopy $\Po$-algebra structure on a dg module $A$ is a morphism of dg operads $\Omega \Po^{\ac} \to \End_A$.

By the preceding subsection applied to $\Co=\Po^{\ac}$ and $\Qo=\End_A$, a  homotopy $\Po$-algebra structure on $A$ is equivalently given by a twisting morphism in $\Tw(\Po^{\ac}, \End_A)$. Recall that the free $\Po^{\ac}$-coalgebra on a dg module $A$ is  $\Po^{\ac}(A):=\bigoplus_{n\in \NN} \Po^{\ac}(n)\otimes_{\Sy_n}A^{\otimes n}$. Using the adjunction $\Hom_{\Sy_n}(\Po^{\ac}(n), \Hom(A^{\otimes n}, A))\cong \Hom(\Po^{\ac}(n)\otimes_{\Sy_n}A^{\otimes n}, A)$, one proves that
$$\Hom_\Sy(\Po^{\ac}, \End_A)\cong \Hom( \Po^{\ac}(A), A)  \cong \Coder(\Po^{\ac}(A)),$$
the space of coderivations on $\Po^{\ac}(A)$. Under this isomorphism, the set of twisting morphisms $\Tw(\Po^{\ac}, \End_A)$ is in one-to-one correspondence with square zero coderivations of degree $-1$, or \emph{codifferentials}, $\Codiff(\Po^{\ac}(A))$ on $\Po^{\ac}(A)$.

\begin{theo}[Rosetta Stone \cite{GinzburgKapranov94, GetzlerJones94, VanderLaan02, MerkulovVallette08I}]\label{4 def theo}
For any Koszul operad $\Po$ in the sense of the Appendix~\ref{Appendix Koszul}, the set of homotopy $\Po$-algebra structures on a dg module $A$ is equal to
\begin{align*}\Hom_{\emph{dg op.}}\left(\Omega \Po^{\ac}, \End_A\right) \cong
\Tw(\Po^{\ac}, \End_A)& \cong  \Hom_{\emph{dg
coaugm. coop.}}\left(\Po^{\ac}, \mathrm{B} \End_A \right)\\&\cong\Codiff(\Po^{\ac}(A)).
\end{align*}
\end{theo}

The operad $\BV$ encoding Batalin--Vilkovisky algebras is Koszul in the sense of Appendix~\ref{Appendix Koszul}, though it is not quadratic or Koszul in the classical sense. In each section of the present paper, we make the notion of \emph{homotopy BV-algebra} explicit using one of the above equivalent definitions each time. In Section~\ref{ResolutionBV}, we make the quasi-free resolution $\Omega \BV^{\ac}$ explicit. In Section~\ref{Homotopy BV-algebra Explicit}, we use the definition in terms of codifferentials on $\BV^{\ac}(A)$.
In Section~\ref{Def theo for BV}, we develop the deformation theory of homotopy BV-algebras using the Lie theoretic description with twisting morphisms. And in Section~\ref{BV Homotopy theory}, we use the last two definitions to study the homotopy theory for homotopy BV-algebras. For instance, the third definition will be shown to provide an explicit formula for Massey products on the homology groups of any BV-algebra or homotopy BV-algebra.

\section{Resolution of the  operad $\BV$}\label{ResolutionBV}

We make  a quasi-free resolution $\BV_\infty$ of the operad $\BV$ encoding Batalin--Vilkovisky
algebras explicit. An algebra over this resolution is a \emph{homotopy BV-algebra}.
The results of this section are an important example of (but can be read independently from) the \emph{inhomogeneous quadratic} Koszul theory  developed  in Appendix~\ref{Appendix Koszul}. This allows us to prove a Poincar\'e--Birkhoff--Witt theorem for the operad $\BV$ itself. At the end of the section, it is shown that the genus $0$ operadic part of a topological conformal field theory, TCFT for short, is encoded in this algebraic notion of homotopy BV-algebra.

\subsection{Batalin--Vilkovisky algebras}

The Koszul--Quillen sign rule says that the image under $f\otimes g$ of $x\otimes y$ is equal to
$(f\otimes g)(x\otimes y)=(-1)^{|x||g|}f(x)\otimes g(y)$. Throughout the paper, we adopt the economical convention of denoting multilinear operations on $A$ as elements in the endomorphism operad $\End_A$, that is, without making them act on the elements of $A$ explicitly. For instance, a binary product $\bullet$ is (graded) symmetric if  the permutation $(12)$ acts trivially on it: $(\textrm{-} \bullet \textrm{-}).(12)=\textrm{-} \bullet \textrm{-}$. This relation, applied to homogeneous elements, gives $x\bullet y =(-1)^{|x||y|}y\bullet x$.

 A skew-symmetric bracket $[\;,\,]$ of degree $0$ on the suspension $sA$ of a space $A$ is equivalent to a symmetric bracket $\langle\; ,  \, \rangle$  of degree $+1$ on $A$ under the formula
$\langle\, \textrm{-} ,  \textrm{-}\,  \rangle= s^{-1}[s(\textrm{-}), s(\textrm{-})]$ because
$$\langle\; ,  \, \rangle.(12)= (s^{-1}[\;, \,]\circ s\otimes s).(12)=
-s^{-1}([\;, \,] .(12))\circ s\otimes s= s^{-1}[\;, \,]\circ s\otimes s=\la \;, \, \ra,$$
which, applied to homogeneous elements, says
\begin{eqnarray*}
\langle y, x \rangle = (-1)^{|y|}s^{-1}[sy , sx]=(-1)^{|y|+1+(|x|+1)(|y|+1)}s^{-1}[sx , sy]=(-1)^{|x||y|}\la x, y\ra.
\end{eqnarray*}

\begin{defi}[Batalin--Vilkovisky algebras]
A \emph{Batalin--Vilkovisky algebra}, or \emph{BV-algebra} for short,  is a differential graded vector space $(A, d_A)$ endowed with
\begin{itemize}
\item[$\triangleright$]  a symmetric binary product $\bullet$ of degree
$0$,

\item[$\triangleright$]  a symmetric  bracket $\langle\; ,  \, \rangle$  of degree $+1$,

\item[$\triangleright$]  a unary operator $\Delta$ of degree
$+1$,
\end{itemize}
such that $d_A$ is a derivation with respect to each of them and such that
\begin{itemize}
\item[$\rhd$] the product $\bullet$ is associative,

\item[$\rhd$] the bracket satisfies the Jacobi identity

$$\langle \langle\; ,  \, \rangle,  \, \rangle \ +\ \langle \langle\; ,  \, \rangle,  \, \rangle.(123) \ +\ \langle \langle\; ,  \, \rangle,  \, \rangle.(321)\ = \ 0,$$

\item[$\rhd$] the product $\bullet$ and the bracket $\langle\;
, \,\rangle$ satisfy the Leibniz relation
$$\langle\, \textrm{-} , \textrm{-}\bullet \textrm{-}\, \rangle\ =\
(\langle\,\textrm{-}, \textrm{-}\,\rangle\bullet \textrm{-})\
+ \ (\textrm{-}\bullet \langle\,\textrm{-},\textrm{-}\,\rangle).(12),   $$

\item[$\rhd$] the operator satisfies $\Delta^2=0$,

 \item[$\rhd$] the bracket is the obstruction to $\Delta$ being a derivation with respect to  the product $\bullet$
$$\langle\,\textrm{-} , \textrm{-}\,\rangle\ =\ \Delta \circ (\textrm{-} \bullet \textrm{-})\ -\
(\Delta(\textrm{-}) \bullet \textrm{-})   \ - \ (\textrm{-} \bullet
\Delta(\textrm{-})),$$

\item[$\rhd$] the operator $\Delta$ is a graded derivation with
respect to the bracket
$$\Delta (\langle\,  \textrm{-}, \textrm{-}\, \rangle) \ + \ \langle\Delta(\textrm{-}),
\textrm{-}\,\rangle \ +\ \langle\, \textrm{-}, \Delta(\textrm{-})\rangle \ = \ 0.$$

\end{itemize}
\end{defi}

\begin{remas}$\  $

\begin{itemize}

\item[$\diamond$] In the literature (\cite{Getzler94,TamarkinTsygan00} for instance), one usually defines a BV-algebra with a degree $1$ bracket $[\;, \, ]$ such that $[x, \, y]=-(-1)^{(|x|+1)(|y|+1)}[y,x]$ and satisfying a Jacobi relation with different signs. This definition is equivalent to the one above under $[x, \, y]=(-1)^{|x|}\la x,y\ra$.

\item[$\diamond$]
The relations are homogeneous with respect to the degree.

\item[$\diamond$]
The Leibniz relation is equivalent to the fact that the operators
$\textrm{ad}_x(\textrm{-}):=\la x, \textrm{-}\,\ra$ are derivations, of
degree $|x|+1$, with respect to the product $\bullet$.

\item[$\diamond$]
A vector space endowed just with a symmetric product
$\bullet$ and a symmetric bracket $\la\; , \, \ra$ which satisfy
the first three relations is called a \emph{Gerstenhaber algebra} after  \cite{Gerstenhaber63}.

\item[$\diamond$]
The last relation is a direct consequence of the other
axioms but we keep it in the definition.
\end{itemize}
\end{remas}
The operad encoding BV-algebras is the operad defined by generators and relations
$$\BV:=\F(V)/(R),$$
where $\F(V)$ denotes the free operad on the $\Sy$-module
$$V:=\bullet\, \KK_2 \;\oplus\; \la\; ,\, \ra\, \KK_2
\; \oplus\;
\Delta\, \KK,$$ with $\KK_2$  the trivial representation of the symmetric group $\Sy_2$.
The space of relations $R$ is the sub-$\Sy$-module
of $\F(V)$ generated by the  relations `$\rhd$' given above. The
basis elements $\bullet$, $\,\la\:,\,\ra$, $\,\Delta$ are of degree 0, 1, 1. Since
 the relations are homogeneous, the operad $\BV$
is graded by this degree, termed the {\em homological degree}.

We denote by $\Com$ the operad generated by the symmetric product $\bullet$ and the associativity relation. We denote by $\Lie_1$ the operad generated by the symmetric bracket $\la \; , \, \ra$ and the Jacobi relation; it is the operad encoding Lie algebra structures on the suspension of a dg module. The operad $\G$ governing Gerstenhaber algebras is defined similarly.

\subsection{The quadratic analogue of $\BV$}

The free operad $\F(V)$ may also be given a \emph{weight grading},
where each generator $\bullet$, $\la \; ,\,\ra$, $\,\Delta$ has weight
one. Equivalently, the weight is the number of internal vertices in
the tree representing an element. The homogeneous component of weight $n$ is denoted
$\F(V)^{(n)}$. The ideal of relations $(R)$ is generated by the {\em inhomogeneous
  quadratic} subspace
$$R\;\;\subset\;\; \F(V)^{(1)}\oplus \F(V)^{(2)}=V \oplus \F(V)^{(2)},$$
so the operad $\BV$ is not weight graded. There are two ways to define a graded operad from $\BV$.

Let $\mathrm{q}\, :\,  \, \F(V) \epi \F(V)^{(2)}$ be the projection
of the free operad onto its quadratic part
and let $\qR$ denote the image of $R$ under $\mathrm{q}$. We consider the
following quadratic operad
$$\mathrm{q}\BV:=\F(V)/(\mathrm{q}R).$$
The quadratic relations of $\BV$ are not modified under $\mathrm{q}$.
The only relation which is not quadratic is the penultimate, which measures
the obstruction for the operator $\Delta$ to be a derivation with respect to the product
$\bullet$. It becomes
$$\mathrm{q}\left(
\vcenter{\xymatrix@M=2pt@R=8pt@C=8pt{
 \ar@{-}[dr] & &\ar@{-}[dl]   \\
 &  {\scriptstyle \bullet }     \ar@{-}[d] &  \\
& {\scriptstyle \Delta} & }} -
\vcenter{\xymatrix@M=2pt@R=8pt@C=8pt{
  \ar@{-}[dr]{ \scriptstyle \Delta} & &\ar@{-}[dl]   \\
 &  {\scriptstyle \bullet }     \ar@{-}[d] &  \\
&  & }} -
\vcenter{\xymatrix@M=2pt@R=8pt@C=8pt{
  \ar@{-}[dr]& & {\scriptstyle \Delta} \ar@{-}[dl]   \\
 &  {\scriptstyle \bullet }     \ar@{-}[d] &  \\
&  & }}  -
\vcenter{\xymatrix@M=2pt@R=8pt@C=8pt{
 \ar@{-}[dr] & &\ar@{-}[dl]   \\
 &  {\scriptstyle \la\, , \, \ra}     \ar@{-}[d] &  \\
& & }}
\right)
 =
 \vcenter{\xymatrix@M=2pt@R=8pt@C=8pt{
 \ar@{-}[dr] & &\ar@{-}[dl]   \\
 &  {\scriptstyle \bullet }     \ar@{-}[d] &  \\
& {\scriptstyle \Delta} & }} -
\vcenter{\xymatrix@M=2pt@R=8pt@C=8pt{
  \ar@{-}[dr]{ \scriptstyle \Delta} & &\ar@{-}[dl]   \\
 &  {\scriptstyle \bullet }     \ar@{-}[d] &  \\
&  & }} -
\vcenter{\xymatrix@M=2pt@R=8pt@C=8pt{
  \ar@{-}[dr]& & {\scriptstyle \Delta} \ar@{-}[dl]   \\
 &  {\scriptstyle \bullet }     \ar@{-}[d] &  \\
&  & }}$$
The operad $\mathrm{q}\BV$ is bigraded, by the
 homological degree and the weight grading.

The weight grading of $\F(V)$
induces a filtration
of $\BV$,
$$F_n\BV:= \pi\big(\bigoplus_{k\leq n}
\F(V)^{(k)}\big),$$
 where $\pi$ denotes the natural
projection $\F(V) \epi \BV$.
As usual the associated graded operad, denoted
by $\gr \BV$, is isomorphic to $\BV$ as an $\Sy$-module.
Since the space of quadratic relations $\mathrm{q}R$ is also zero in $\gr\BV$, there
is a natural epimorphism of bigraded operads
$$ \qBV\;\epi\; \gr \BV. $$

\begin{theo}\label{qBV Koszul}
The operad $\qBV$ is Koszul.
\end{theo}

\begin{proof}
This follows from  \cite[Proposition~$8.4$]{Vallette07} by the distributive law method. Observe that the operad $\qBV$ is obtained from the Koszul operad $\G$ of Gerstenhaber algebras and $D:=\KK[\Delta]/(\Delta^2)$ the algebra of dual numbers, an operad concentrated in arity $1$, by means of the distributive laws
\begin{eqnarray*}
\vcenter{\xymatrix@M=2pt@R=8pt@C=8pt{
 \ar@{-}[dr] & &\ar@{-}[dl]   \\
 &  {\scriptstyle \bullet }     \ar@{-}[d] &  \\
& {\scriptstyle \Delta} & }} &\longmapsto&
\vcenter{\xymatrix@M=2pt@R=8pt@C=8pt{
  \ar@{-}[dr]{ \scriptstyle \Delta} & &\ar@{-}[dl]   \\
 &  {\scriptstyle \bullet }     \ar@{-}[d] &  \\
&  & }}+
\vcenter{\xymatrix@M=2pt@R=8pt@C=8pt{
  \ar@{-}[dr]& & {\scriptstyle \Delta} \ar@{-}[dl]   \\
 &  {\scriptstyle \bullet }     \ar@{-}[d] &  \\
&  & }}\\
\vcenter{\xymatrix@M=2pt@R=8pt@C=8pt{
 \ar@{-}[dr] & &\ar@{-}[dl]   \\
 &   {\scriptstyle \la \, , \, \ra}    \ar@{-}[d] &  \\
& {\scriptstyle \Delta} & }} &\longmapsto&
\vcenter{\xymatrix@M=2pt@R=8pt@C=8pt{
  \ar@{-}[dr]{ \scriptstyle \Delta} & &\ar@{-}[dl]   \\
 &   {\scriptstyle \la\, , \, \ra}     \ar@{-}[d] &  \\
&  & }} +
\vcenter{\xymatrix@M=2pt@R=8pt@C=8pt{
  \ar@{-}[dr]& & {\scriptstyle \Delta} \ar@{-}[dl]   \\
 &   {\scriptstyle \la\, , \, \ra}    \ar@{-}[d] &  \\
&  & }}
\end{eqnarray*}
These define an injective map $D \circ \G \mono \G \circ D \cong \qBV$ of graded $\Sy$-modules and hence \cite[Proposition~$8.4$]{Vallette07} proves that the operad $\qBV$ is Koszul. (The fact that the algebra of dual numbers here is not concentrated in degree zero does not affect the spectral sequence used in the proof there: we have $E^1_{p q} =({D^{\ac}}^{(p)} \circ D)_{p+q}$ and $E^2_{pq}=0$ unless $p=q=0$ since $d^1$ is the differential of the Koszul complex.)
\end{proof}

\begin{rema}
The distributive law theorem for operads due to Markl in \cite{Markl96D}
cannot be used here because it does not apply to operads with unary operations. The general methods of
\cite[Section 8]{Vallette07} are the same, but the results include unary operations since the proofs are based on different filtrations.
\end{rema}

As a direct corollary of the distributive law methods, we get the following explicit  description
of the $\Sy$-modules underlying $\qBV$ and its Koszul dual $\qBV^{\ac}$, with their
homological grading.

\begin{prop}\label{form BV dual}
There are isomorphisms of degree graded $\Sy$-modules
\begin{eqnarray*}
&&\qBV\cong\G\circ D\cong\Com \circ \Lie_1 \circ \KK[\Delta]/(\Delta^2) \\
 &&  \qBV^{\ac}\cong D^{\ac}\circ \G^{\ac}\cong \KK[\delta] \circ
 \Lie_1^{\ac} \circ \Com^{\ac},
 \end{eqnarray*}
where $\KK[\delta]\cong T^c(\delta)$ is the cofree coalgebra on a degree $2$
generator $\delta$, that is, a cooperad concentrated in arity $1$.
\end{prop}

\begin{proof}
It follows from \cite[Lemma~$8.1$ and
Proposition~$8.2$]{Vallette07}.
\end{proof}

When the number of generators of a binary  operad $\Po$ is finite, it is convenient to work with the Koszul dual operad
$$\Po^!=(\F(V)/(R))^!:=\F(V^*\otimes_H \sgn_{\Sy_2})/(R^\perp).$$
Here the Hadamard product $\otimes_H$
is the arity-wise tensor product of $\Sy$-modules.
In the finitely generated case, the Koszul dual cooperad $\Po^{\ac}$ is isomorphic to
$$\Po^{\ac}\cong S^c(\Po^!)^*:=S^c
\otimes_H (\Po^!)^*,$$
where $(\Po^!)^*$ denotes the component-wise linear dual of the Koszul dual
operad  and where $S^c:=\End^c_{s^{-1}\KK}$ is the \emph{suspension} cooperad of endomorphisms
of the one dimensional vector space $s^{-1}\KK$ concentrated in degree $-1$ (see Section~$2.4$ of \cite{GetzlerJones94} and
Section~$2$ of \cite{Vallette08}). 
We can also consider the desuspension of  operads by taking the Hadamard product $S^{-1}\Po:=S^{-1}\otimes_H \Po$ with the  operad  $S^{-1}:=\End_{s\KK}$.
Observe the Hadamard product of two (co)operads is again a (co)operad. 
\begin{coro}
The underlying $\Sy$-module of the cooperad $\qBV^{\ac}$ is isomorphic to
$$\qBV^{\ac}\cong \KK[\delta]\circ S^c\Com_1^c \circ S^c\Lie^c,$$
where $\Lie^c\cong\Lie^*$ is the cooperad encoding Lie coalgebras and where $\Com_1^c\cong \Com_{-1}^*$ is the cooperad encoding cocommutative
 coalgebra structures on the suspension of a dg module.
 The degree of the elements in $$\KK.\delta^d \otimes
S^c\Com_1^c(t) \otimes S^c\Lie^c(p_1)  \otimes \ldots \otimes
S^c\Lie^c(p_t)\subset \qBV^{\ac}$$ is $n+t+2d-2$, for
$n=p_1+\cdots+p_t$.
\end{coro}

\begin{proof}
From \cite{GinzburgKapranov94,GetzlerJones94}, we have $\Com^!= \Lie$ and ${\Lie_1}^!=\Com_{-1}$.
\end{proof}

Since $\qBV$ is a Koszul operad, its  quadratic model is given by the cobar
construction of its Koszul dual cooperad
$$\Omega\,  \qBV^{\ac}:=(\F(s^{-1}\overline{\qBV}^{\ac}),d_2)\xrightarrow{\sim} \qBV\epi \gr \BV\cong \BV.$$

The $\Sy$-module $s^{-1}\overline{\qBV}^{\ac}$ seems to be a good
candidate for the space of  generators a quasi-free
resolution of the operad $\BV$ itself. In the next section, we
perturb the quadratic differential $d_2$ by a
linear term $d_1$ in order to obtain a quasi-isomorphism with the operad $\BV$.

\subsection{The Koszul dual of $\BV$}\label{koszul-dual}

We consider the map $\varphi\, :\, \qR \to V$ defined by
$$ \vcenter{\xymatrix@M=2pt@R=8pt@C=8pt{
 \ar@{-}[dr] & &\ar@{-}[dl]   \\
 &  {\scriptstyle \bullet }     \ar@{-}[d] &  \\
& {\scriptstyle \Delta} & }} -
\vcenter{\xymatrix@M=2pt@R=8pt@C=8pt{
  \ar@{-}[dr]{ \scriptstyle \Delta} & &\ar@{-}[dl]   \\
 &  {\scriptstyle \bullet }     \ar@{-}[d] &  \\
&  & }} -
\vcenter{\xymatrix@M=2pt@R=8pt@C=8pt{
  \ar@{-}[dr]& & {\scriptstyle \Delta} \ar@{-}[dl]   \\
 &  {\scriptstyle \bullet }     \ar@{-}[d] &  \\
&  & }}  \quad\mapsto\quad
\vcenter{\xymatrix@M=2pt@R=8pt@C=8pt{
 \ar@{-}[dr] & &\ar@{-}[dl]   \\
 &  {\scriptstyle \la\; , \, \ra}     \ar@{-}[d] &  \\
& & }}$$ and $0$ otherwise, so that the graph of $\varphi$  is equal
to the space of relations $R$ (see Appendix~\ref{filtered properad} for the general theory).

We use the notation $\odot$ for the `commutative' tensor product, that
is, the quotient of the tensor product under the permutation of
terms. In particular, we denote by $\delta^d\otimes L_1 \odot \cdots \odot L_t$
a generic element  of $\KK[\delta]\circ S^c\Com_1^c \circ S^c\Lie^c$
 with $L_i\in S^c\Lie^c$, for $i=1,\ldots, t$; the elements of $S^c\Com_1^c$ being implicit.

\begin{lemm}\label{dvarphi}
There is a unique square-zero coderivation $d_\varphi$
on the cooperad $\qBV^{\ac}$ which extends $\varphi$.
It is explicitly given by
$$d_\varphi(\delta^d\otimes L_1 \odot \cdots \odot L_t)=\sum_{i=1}^t (-1)^{\varepsilon_i}  \delta^{d-1}
\otimes L_1\odot \cdots\odot L_i'\odot L_i''\odot \cdots
\odot L_t ,$$ where $L_i'\odot L_i''$ is Sweedler-type
notation for  the image of $L_i$ under the binary part
 $$S^c\Lie^c \to
S^c\Lie^c(2)\otimes (S^c\Lie^c\otimes S^c\Lie^c)$$
of the decomposition
map of the cooperad $S^c\Lie^c$.
The sign is
given by $\varepsilon_i=(|L_1|+\cdots+|L_{i-1}|)$. The image of
$d_\varphi$  is equal to $0$ when $d=0$ or when $L_i\in
S^c\Lie^c(1)=I\KK$ for all $i$.
\end{lemm}

\begin{proof}
The first claim follows from Lemma~\ref{coderivation varphi}.  It is straightforward to verify that
$(R\otimes V + V \otimes R)\cap \F(V)^{(2)}\subset R\cap \F(V)^{(2)}$
here.

To make $d_\varphi$ explicit, we dualize everything and we consider the
unique derivation ${}^td_\varphi$  on $(\qBV^{\ac})^*$  extending
${}^t\varphi$. The transpose map ${}^t\varphi$ is equal to
${}^t\varphi(c)=\delta^*\otimes I \otimes l\in \KK[\delta^*]\circ S^{-1}\Com_{-1}
\circ S^{-1}\Lie$, where $c$ denotes the generator of $S^{-1}\Com_{-1}(2)$ and
where $l$ denotes the generator of $S^{-1}\Lie(2)$. Therefore,
${}^td_\varphi$ is equal to
$${}^td_\varphi((\delta^*)^d\otimes L_1^* \odot \cdots \odot L_t^*)=\!\!\!\!
\sum_{1\leq i<j\leq t}\!\!\!\!
(-1)^{\varepsilon_{ij}} (\delta^*)^{d+1}\otimes [L_i^*, L_j^* ]\odot L_1^*
\odot \cdots \odot \hat{L_i^*} \odot \cdots \odot \hat{L_j^*}
\odot \cdots \odot L_t^*,$$ where $[L_i^*, L_j^*]$ is the
image $\gamma (l ; L_i^*, L_j^*)$ under the composition map $S^{-1}\Lie(2)\otimes (S^{-1}\Lie
\otimes S^{-1}\Lie) \to S^{-1}\Lie$ of the operad $S^{-1}\Lie$ with
the generating element $l$ in $S^{-1}\Lie(2)$. The sign $\varepsilon_{ij}$ is
automatically given by Koszul--Quillen sign rule, that is here
$$\varepsilon_{ij}=(|L_i^*|+|L_j^*|).(|L_1^*|+\cdots+|L_{i-1}^*|)+|L_j^*|.(|L_{i+1}^*|+\cdots+|L_{j-1}^*|) .$$
Finally, we dualize once again  to get $d_\varphi$.
\end{proof}
Since $S^{-1}\Lie(A)\cong s\Lie(s^{-1}A)$, this last formula shows that the differential
${}^td_\varphi$ on $(\qBV^{\ac})^*$ is equal, up to the powers of $\delta^*$,
to the Chevalley--Eilenberg boundary map defining the homology of the free Lie algebra.

We can now define the Koszul dual of the inhomogeneous quadratic
Batalin--Vilkovisky operad.
\begin{defi}
The \emph{Koszul dual cooperad} of the operad $\BV$ is the differential graded cooperad
$$\BV^{\ac}:=(\qBV^{\ac}, d_\varphi).$$
\end{defi}

\subsection{The Koszul resolution of the operad $\BV$}

Since the definition of the cobar construction extends to dg cooperads, we can consider the cobar construction $\Omega \BV^{\ac}$ of the Koszul dual dg cooperad of $\BV$.

\begin{defi}[Operad $\BV_\infty$]
We denote by $\BV_\infty$ the quasi-free operad given by the cobar construction on $\BV^{\ac}$:
$$\BV_\infty := \Omega \BV^{\ac}\cong(\F(s^{-1}\overline{\qBV}^{\ac}),\, d=d_1+d_2). $$
\end{defi}
The space of generators of this quasi-free operad is isomorphic to
$\KK[\delta]\circ S^c\Com_1^c \circ S^c\Lie^c$, up to coaugmentation. The differential $d_1$ is
the unique derivation which extends the internal differential
$d_\varphi$ and the differential $d_2$ is the unique derivation
which extends the partial coproduct of the cooperad $\BV^{\ac}$. (Recall that the partial coproduct of a cooperad is the component of the coproduct in $\F(\Co)^{(2)}\subset \Co \circ \Co$, given by trees with $2$ vertices.) 
Therefore, the total derivation $d=d_1+d_2$  squares to zero
 and it faithfully encodes the full data of dg cooperad on
$\BV^{\ac}$.

 \begin{theo}\label{MainBV}
The operad $\BV_\infty$ is a resolution  of the operad $\BV$
$$\BV_\infty=\Omega \BV^{\ac} =\left(\F(s^{-1}\overline{\qBV}^{\ac}),\,d=d_1+d_2\right)\;\qi \;\BV.$$
\end{theo}

\begin{proof}
It is an application of  Theorem~\ref{MainTHM}.
\end{proof}

As a direct corollary of the proof of this theorem, we get the
following seemingly unrelated result.

\begin{theo}[Poincar\'e--Birkhoff--Witt Theorem for the operad $\BV$] \label{PBW for BV}
The natural epimorphism $\qBV \epi \gr \BV$ is an isomorphism of
bigraded operads. Thus there are isomorphisms of
$\Sy$-modules
$$\BV \cong \gr\BV \cong \qBV$$
which preserve the homological degree.
\end{theo}

\begin{proof}
It is an application of Theorem~\ref{PBW for P}.
\end{proof}

This algebraic result allows us to make the free Batalin--Vilkovisky algebra explicit.

\begin{prop}[\cite{Getzler94}]
The underlying module of the free Batalin--Vilkovisky algebra on  a module $X$ is isomorphic to
$$\BV(X)\cong \qBV(X) \cong\G\circ D(X)\cong\Com \circ \Lie_1 (X\;\oplus \;\Delta(X)). $$
\end{prop}

The result is also a corollary of Proposition~$4.8$ of \cite{Getzler94}, but the method is different.

\begin{rema}
Our construction of $\BV^{\ac}$ and $\BV_\infty$ is based on a particular presentation of
$\BV$  in terms of generators and relations. This method cannot be applied to the other quadratic operad $\qBV'$,
associated to the presentation of $\BV$ in which the
final `redundant' derivation relation between $\Delta$ and the bracket $\la\:, \, \ra$ is omitted.
The two presentations of $\BV$ are equivalent, but the induced quadratic operads $\qBV$ and
$\qBV'$ are not isomorphic. In fact the derivation
relation between $\Delta$ and $\la\;, \, \ra$ holds in $\gr \BV'$ but not
in $\qBV'$, so in weight two $\qBV'$ is strictly larger
than $\gr \BV'$ and there is no PBW theorem. To apply
this method one must always consider the {\em  maximal} space of relations
$R$, see Proposition~\ref{Max Relations} in  Appendix~\ref{Appendix Koszul} for further details.
\end{rema}

\subsection{Homotopy BV-algebra}

A \emph{homotopy BV-algebra}, or \emph{$\textrm{BV}_\infty$-algebra}, or \emph{BV-algebra up to homotopy}, is an algebra over the operad $\BV_\infty$.
Such a structure on a dg module $A$ is given by a morphism of dg
operads $$\BV_\infty \to \End_A.$$
Any BV-algebra is a special case of a homotopy BV-algebra; this is
the case when
the structure map $\BV_\infty \to \End_A$ factors as the
 quasi-isomorphism followed by the $\BV$ structure map,
$$\BV_\infty \qi \BV \to \End_A.$$
We will return to the explicit definition of homotopy BV-algebras
in Section \ref{Homotopy BV-algebra Explicit}.

\subsection{Relations with the framed little {\disc} operad and double loop spaces}\label{fD and 2-Omega}

We recall from Section~$4$ of \cite{Getzler94} (see also \cite{SalvatoreWahl03}) the definition of the \emph{framed little discs operad $\fD$}. Let $D$ denote the unit disc in the plane. The space $\fD(n)$ is the space of maps from the disjoint union of $n$ copies of $D$  to $D$,  whose restriction to each  disc
is the composite of a translation and multiplication by an element of $\mathbb{C}^\times$. The interior of the images of the  discs are supposed to be disjoint. It can be thought of as the space of configurations of $n$ discs, with one marked point on each boundary, inside the unit disc $D$. Since we work over a field, the singular chains $C_\bullet(\fD)$ is a dg operad and its homology groups $H_\bullet(\fD)$ form a graded operad.

\begin{prop}[\cite{Getzler94}]\label{Homology Fd}
The homology of the framed little {\disc} operad is isomorphic to the $\BV$ operad,
$$H_\bullet(\fD)\cong \BV.$$
\end{prop}

Therefore, the framed little {\disc} operad gives a topological model for the notion of  homotopy BV-algebras. We make this relation precise with the following key result.

\begin{prop}[\cite{Severa09, GiansiracusaSalvatore09}]\label{Formality Fd}
The framed little {\disc} operad is formal, i.e.\ the dg operads
$C_\bullet(\fD)$ and $H_\bullet(\fD)$
are linked by a zig-zag of quasi-isomorphisms.
$$
C_\bullet(\fD)=
X_1
\stackrel{\sim}\longleftarrow
X_2
\stackrel{\sim}\longrightarrow
\cdots \stackrel{\sim}\longleftarrow
X_{2r  }
\stackrel{\sim}\longrightarrow
X_{2r+1}=H_\bullet(\fD)\cong\BV
$$
\end{prop}

\begin{theo}
There is a quasi-isomorphism of dg operads between the operads $\BV_\infty$ and $C_\bullet(\fD)$, the chains of the framed little {\disc} operad, which lifts the resolution $\BV_\infty \qi \BV$.
\end{theo}

\begin{proof}
Let $X\xrightarrow{\sim} C_\bullet(\fD)$  be a cofibrant replacement of $C_\bullet(\fD)$ in the model category of dg operads (see \cite{Hinich97, BergerMoerdijk03}). We consider the  following zig-zag coming from Proposition \ref{Formality Fd}
$$X\xrightarrow{\sim}
C_\bullet(\fD)=
X_1
\stackrel{\sim}\longleftarrow
X_2
\stackrel{\sim}\longrightarrow
\cdots \stackrel{\sim}\longleftarrow
X_{2r}
\stackrel{\sim}\longrightarrow
X_{2r+1} =\BV
$$
in the homotopy category of dg operads $\mathsf{Ho}(X, \BV)$. In this model category every dg operad is fibrant, so is  $\BV$. Since $X$ is cofibrant and $\BV$ fibrant, there is a morphism of dg operads $\xi\, :\, X \to \BV$ which factors
$$\xymatrix{ C_\bullet(\fD)=
X_1 & \ar[l]_(0.3){\sim} X_2\ar[r]^{\sim}  & \cdots & \ar[l]_{\sim} X_{2r} \ar[r]^(0.3){\sim} & X_{2r+1} =\BV\\
X \ar[u]^{\sim} \ar@{->>}[rrrr]^{\sim}& & & & \BV \ar@{=}[u]
}$$
(see Proposition~$5.11$ of \cite{DwyerSpalinski95}). Since all the maps are quasi-isomorphisms, then $\xi$ is a quasi-isomorphism. And since the dg operad $\BV$ has trivial differential, $\xi$ is an epimorphism. Finally, the morphism of dg operads $\xi$ is a fibration. The operad $\BV_\infty$ is quasi-free and non-negatively graded, so it is a cofibrant replacement of $\BV$. In conclusion, there is a quasi-isomorphism of dg operads which lifts
$$\xymatrix@R=20pt@C=30pt{ \BV_\infty \ar[d]_{\sim} \ar@{..>}[dr]^{\sim}& & \\
\BV \ar[r]& X \ar@{->>}[l]_\sim \ar[r]^(0.4)\sim & C_\bullet(\fD)
}$$
 \end{proof}

This quasi-isomorphism exists by an abstract model category
argument and the formality theorem. So it would be desirable to have a better understanding of it, for
instance in terms of a morphism to the cellular chains on a cellular decomposition of the framed little {\disc} operad (see \cite{Kaufmann04}). This would lead to a better understanding of the Grothendieck--Teichm\"uller group, see \cite{Tamarkin02, Tamarkin07}. We leave this interesting question for future study or to the reader. \\

The main motivation for the introduction of the framed little discs operad by Getzler \cite{Getzler94} was to extend the recognition principle of Boardman--Vogt \cite{BoardmanVogt73} and May \cite{May72}: the framed little discs operad acts on the double loop space of any topological space endowed with an action of the circle and the operad action characterises this type of spaces \cite{SalvatoreWahl03}. Hence the homology of the double loop space of any topological space endowed with an action of the circle carries a BV-algebra structure by Proposition~\ref{Homology Fd}, where the product is the Pontryagin product and the bracket is the Browder bracket.

\begin{coro}
Let $X$ be a topological space endowed with an action of the circle. The singular chain complex $C_\bullet(\Omega^2 X)$ of the double loop space of $X$ carries a homotopy BV-algebra structure
which lifts the BV-algebra structure on homology.
\end{coro}

\subsection{Relations with the Riemann sphere operad and topological conformal field theory}\label{BVinfty and TCFT}

Recall from Segal \cite{Segal04} that the \emph{properad of Riemann
  surfaces} $\mathcal{R}$ is defined by the moduli space of
isomorphism classes of connected Riemann surfaces of arbitrary genus
with biholomorphic maps from the disjoint union of $n+m$
{\disc}s. Again, the images of the {\disc}s are supposed to have
disjoint interiors. The images of the first $n$ {\disc}s form the
inputs and the images of the last $m$ {\disc}s form the outputs. The
properad structure is given by sewing along the
boundaries of the {\disc}s. (In the literature, one usually
considers the {\em prop} freely associated to this properad but this
yields no more information.) A space which is an algebra over $\mathcal{R}$ is
called a \emph{Conformal Field Theory}, or \emph{CFT} for short. The
singular chains $C_\bullet(\mathcal{R})$ of this topological properad
is a dg properad, and a dg module which is an algebra over $C_\bullet(\mathcal{R})$ is  called a \emph{Topological Conformal Field Theory}, or \emph{TCFT} for short.

If we consider only Riemann spheres (surfaces of genus $0$) with one
output {\disc}, this collection forms a topological operad denoted $R$
of which the framed little {\disc} operad $\fD$ is a deformation
retract. The retraction is obtained by shrinking the complement of
the output {\disc} to a unit {\disc}, leaving a configuration of $n$ {\disc}s in the
unit {\disc} with one marked point on each boundary. This result, with the help of Proposition~\ref{Homology Fd}, allowed Getzler \cite{Getzler94} to prove that the homology groups of any TCFT carry a BV-algebra structure.

\begin{theo}\label{TCFTBV}
There is a quasi-isomorphism of dg operads $\BV_\infty \qi
C_\bullet(R)$ such that the following diagram of quasi-isomorphisms commutes
$$\xymatrix@R=20pt@C=30pt{ C_\bullet(R) \ar@{->>}@<-0.5ex>[r] & \ar@{>->}@<-0.5ex>[l]\ \  C_\bullet(\fD) \ar@{->>}[d] \\
\BV_\infty \ar[r] \ar[ur] \ar[u] & \BV \cong H_\bullet(\fD).}$$
\end{theo}

\begin{coro}
Any TCFT carries a homotopy BV-algebra structure which lifts the BV-algebra structure on homology.
\end{coro}

Thus an important part of a TCFT structure is encoded in the algebraic notion of homotopy BV-algebra structure. For the same kind of result at the level of vertex algebras, we refer the reader to Section~\ref{Vertex Algebras}.

\section{Homotopy BV-algebra}\label{Homotopy BV-algebra Explicit}

We translate the results on resolutions of
operads from the first section to give the explicit definitions of
homotopy BV-algebras in terms of operations and relations. The operations carried by
these homotopy algebras correspond to the generators of the quasi-free
resolutions we have described. It remains to give explicit descriptions of
the axioms the operations must satisfy. They are given by the relation $(d_1+d_2)^2=0$.

Throughout this section we will consider, as examples of these structures, the Hochschild
cochain complex of associative algebras.

\subsection{Homotopy Gerstenhaber algebra}

It was proved in Getzler--Jones \cite{GetzlerJones94} that the operad
$\G$ encoding Gerstenhaber algebras is Koszul, thereby defining the notion of homotopy Gerstenhaber algebras. To give a homotopy Gerstenhaber algebra structure on a graded vector space
$A$ is to give a morphism of dg operads
$$\Omega \G^{\ac} =\bigl(\F(s^{-1}\overline{\G}^{\ac}),\,d_2\bigr)\;\longrightarrow\; (\End_A, d_{\End_A})$$
where $\G^{\ac}$ is the Koszul dual cooperad.
Recall from \cite[Proposition 4.2.14]{GinzburgKapranov94} that this
is equivalent to giving a map
$$
m:\G^{\ac}(A)\to A
$$
which extends to a {\em square zero} coderivation (see Tamarkin--Tsygan \cite[Section~$1$]{TamarkinTsygan00} for the form of these operations).
Here $\G^{\ac}(A)$ is the cofree Gerstenhaber coalgebra, up to (de)suspension,
$$
\G^{\ac}(A)\;\;=\;\;s^{-1}\G^c(sA)\;\;=\;\;s^{-2}\Com^c(s\Lie^c(sA)),
$$
where $\G^c$ is the component-wise linear dual to $\G$.
The square zero condition gives a family of quadratic relations
on the components of the map $m$. We will make these relations completely explicit, using
the laws for extending $m$ to a coderivation with respect to
the Gerstenhaber coalgebra structure of $\G^{\ac}(A)$ (compare with \cite[Lemma~$3.4$ \& Theorem~$3.6$]{Ginot04}, where the formula for $iii)$ has to be slightly modified).
\medskip

Given elements $a_k\in sA$, a permutation $\sigma\in\Sy_n $ and an interval $[i,j]=\{i,i+1,\dots ,j\}$, $i\leq j\leq n$, 
we will use the notation:
\begin{align*}
a_{[i,j]}&:=
a_{i}\otimes a_{i+1}\otimes\dots\otimes
a_{j}
\\
a^\sigma_{[i,j]}&:=
a_{\sigma^{-1}(i)}\otimes a_{\sigma^{-1}(i+1)}\otimes\dots\otimes
a_{\sigma^{-1}(j)}
\end{align*}
in $(sA)^{\otimes(j-i+1)}$.

\medskip

The cofree Lie coalgebra $\Lie^c(sA)$ on the dg module $sA$ is given by
$$
\Lie^c(sA)=\bigoplus_{p\geq1}\overline{(sA)^{\otimes p}}.
$$
Here $\overline{(sA)^{\otimes p}}$ is the quotient of the $(sA)^{\otimes p}$ by the
images
of the {\em shuffle maps}
\begin{align*}\textit{sh}_{i,p-i}\colon (sA)^{\otimes i}\otimes
(sA)^{\otimes (p-i)}&\longrightarrow (sA)^{\otimes p}
,\qquad\qquad (1\leq i\leq p-1),
\\
\textit{sh}_{i,p-i}(a_{[1,i]}\otimes a_{[i+1,p]})&=\!\!\!\!
\sum_{\sigma \in Sh_{i,p-i}}\!\!\!\! (-1)^\varepsilon
a^\sigma_{[1,p]}.
\end{align*}
The sum is over the set $Sh_{i,p-i}$ of all $(i,p-i)$-shuffles, and $(-1)^\varepsilon$ is
the Koszul sign associated to the reordering of the tensor factors $a_k\in sA$. This result is due to
Ree \cite{Ree58}. We refer to  \cite[Chapter~$1$]{LodayVallette09} for more details and
the following operadic interpretation: the morphism of operads $\Ass \epi \Com$, equivalent to the
forgetful functor from commutative and associative algebras to associative algebras, induces a morphism
on the level of the Koszul dual cooperads $\Ass^{\ac} \epi \Com^{\ac}$. After desuspension, it
gives a morphism of cooperads $\Ass^c \epi \Lie^c$, whose kernel is equal to the image of the shuffle maps.

\medskip

We now write $W=\Lie^c(sA)$ and turn to the construction of
$\Com^c(sW)$. It is well known that
the free cocommutative coalgebra $\Com^c(sW)$
on a graded space $sW$ is
given by
$$
\Com^c(sW)=\bigoplus_{t\ge 1} (sW)^{\odot t} \cong \bigoplus_{t\ge 1} s^t\bigwedge{\!\!}^tW
.$$
Here $(sW)^{\odot t}$ is the space of symmetric tensors in $(sW)^{\otimes
  t}$ which we identify with the
exterior power $\bigwedge^tW$, up to sign and a change of grading
as discussed in Section $1.1$.
Furthermore
any family of maps $m'_t:\bigwedge^tW\to W$ of degrees
$t-2$
extends to a coderivation
$\widetilde m$ on $\Com^c(sW)$ by
\begin{align*}
\widetilde m(w_1\wedge \dots\wedge w_t)\,=\!\!\!\!\!
\sum_{I\sqcup J=\{1,\dots,t\}
}\!\!\!\!\!\!\!\!(-1)^{\varepsilon+\varepsilon'+\varepsilon''}
\;s^{\otimes|J|}\;m'_{|I|}(w_{i_1}\wedge\dots\wedge
w_{i_{|I|}})\wedge w_{j_1}
\wedge\dots\wedge
w_{j_{|J|}}
.
\end{align*}
Here $(-1)^{\varepsilon+\varepsilon'}$
is the shifted Koszul sign, and $\varepsilon''=|I|\cdot|J|$.
The shifted Koszul sign associated to a permutation $\sigma$ of tensor
factors is the usual Koszul sign multiplied by the sign of $\sigma$
itself. Partitions $I\sqcup J$ can equivalently be regarded as
(un)shuffles of $\{1,\dots,t\}$.

\smallskip

Throughout the text, we will simply denote any tuple $(q_1, \ldots, q_t)$ by $\mathsf{q}$.

\def\ELL{L}
\begin{defi}[Straight shuffles]
Consider integers $t\geq1$ and $\ell_k,r_k\geq0$, $q_k\geq1$,
$p_k=\ell_k+q_k+r_k$ for each $k=1,\dots,t$, and let
$$P_k=\sum_{1\leq i\leq k-1}p_i,\qquad\qquad
\ELL_k=\sum_{1\leq i\leq t}\ell_i+\sum_{1\leq i\leq k-1}{q}_i.
$$
Then  a {\em straight
$(\mathsf{l}, \mathsf{q}, \mathsf{r}, \mathsf{p})$-shuffle}
is a $\mathsf{p}$-shuffle (i.e. a $(p_1,\dots,p_t)$-shuffle) $\sigma$ satisfying the following extra property
for each $k=1,\dots,t$:
$$\sigma[P_k+\ell_k+1,P_k+\ell_k+{q}_k]=[\ELL_k+1,\ELL_{k+1}].
$$
By a {\em straight
$(\mathsf{q}, \mathsf{p})$-shuffle}
we mean a straight $(\mathsf{l}, \mathsf{q}, \mathsf{r}, \mathsf{p})$-shuffle for some
values of $\ell_k,r_k\geq0$ with $\ell_k+r_k=p_k-q_k$.
\end{defi}

Less formally, a straight shuffle may be thought of as a permutation
$\sigma$ of
a concatenation $X$ of $t$ strings of lengths $p_k$, each of which contains a
non-empty distinguished interval of length $q_k$, with $l_k$ elements on the left and $r_k$ elements on the right. For example,
$$
X\;\;=\;\;1\;\underline{2}\;3\;4\;|\;\underline{5\;
6}\;7\;8\;|\;9\;\underline{10\;11}\;
|\;\underline{12\;13}.
$$
For the permutation to be a straight shuffle it must satisfy following conditions:
\begin{itemize}
\item[$\diamond$]
the orders of the elements within the $t$ strings are preserved by $\sigma$ (\emph{shuffle}),
\item[$\diamond$]
the distinguished intervals appear unchanged and contiguously in the image (\emph{straight}).
\end{itemize}
For example, one possible straight shuffle of $X$ is
$$
\quad\qquad 9\;1\;\;\underline{2\;5\;6\;10\;11\;12\;13}\;\;7\;3\;8\;4.
$$
More precisely, it is a straight $\big((1,0,1,0), (1,2,2,2), (2,2,0,0),(4,4,3,2)\big)$-shuffle.

\begin{defi}[Straight shuffle extension]
For any map 
$$
m_\mathsf{q}=m_{{q}_1,\dots,{q}_t}\colon
\,\overline{(sA)^{\otimes {q}_1}}\;\wedge\dots\wedge\; \overline{(sA)^{\otimes {q}_t}}
\longrightarrow sA,\qquad\quad\;\;
\qquad\qquad(q_k\ge 1),\qquad\qquad\qquad\qquad
$$
and integers $p_k\geq q_k$ we define the {\em straight shuffle extension}
\begin{align*}
m^\mathsf{p}_\mathsf{q}=m^{p_1,\dots,p_t}_{{q}_1,\dots,{q}_t}\colon
&
\overline{(sA)^{\otimes p_1}}\;\wedge\dots\wedge \;\overline{(sA)^{\otimes p_t}}
\longrightarrow \overline{(sA)^{\otimes p'}},
\qquad(p'=1+\sum_{k=1}^t (p_k-{q}_k)),
\\[2mm]
&\overline{a_{[1,\,p_1]}}\wedge
\dots\wedge\overline{
a_{[P_t+1,\,P_{t+1}]}
}
\longmapsto  \\[2mm]
&
\!\!\!\!\!\!\!\!\!\!\!\!\!\!\!\!\!\!\!\!\!\!\!\!\!\!\!\!\!\!\!\!
\sum_\sigma
(-1)^{\varepsilon+\varepsilon'}
\,\overline{\!
a^\sigma
_{[1,\ELL_1]}\otimes m_{{q}_1,\dots,{q}_t}(
\overline{a^\sigma_{[\ELL_1+1,\ELL_2]}}
\wedge\dots\wedge\overline{a^\sigma_{[\ELL_t+1,\ELL_{t+1}]}}
)\otimes a^\sigma_{[\ELL_{t+1}+1,P_{t+1}]\!\!}
}\;\;
.
\end{align*}
Here $\sigma$ runs over all straight $(\mathsf{q},\mathsf{p})$-shuffles and the integers $P_k$, $\ELL_k$ are as above. The sign $(-1)^\varepsilon$ is the Koszul sign associated to the reordering of the $a_i\in sA$, and $\varepsilon'\,=\,t\cdot|a^\sigma_{[1,\ELL_1]}|$.
\end{defi}

{From} the straight shuffle of $X$ above, for example, we see that one of the terms in the expression
$$
m_{1,2,2,2}^{4,4,3,2}(\overline{a_{[1,4]}}\wedge\overline{a_{[5,8]}}\wedge\overline{a_{[9,11]}}\wedge\overline{a_{[12,13]}})\;\;\in\;\;\overline{(sA)^{\otimes 7}}
$$
is
$$
\pm\;\overline{a_9\otimes a_1\otimes m_{1,2,2,2}({a_{2}}\wedge\overline{a_{[5,6]}}\wedge\overline{a_{[10,11]}}\wedge\overline{a_{[12,13]}})\otimes a_7\otimes a_3\otimes a_8\otimes a_4}.
$$

\medskip

\begin{remas}{~}\par
\begin{enumerate}
\item[$\diamond$]
The straight shuffle extension
$m_\mathsf{q}^\mathsf{p}$ is a multilinear map of the same homological
degree as $m_\mathsf{q}$, and inherits
the same (skew-)symmetry properties.

\item[$\diamond$]
If we allowed $p_j>0$ but $q_j=0$ for some $1\leq j\leq t$ then the straight shuffle extension would vanish, since elements in the image lie in that of the shuffle map
$\textit{sh}_{p_j,p'-p_j}$.

\item[$\diamond$]
The straight shuffle extension is well defined on the quotients $\overline{(sA)^{\otimes p_k}}$.
That is, if  $1\leq i\leq p_j-1$ then
the image of the composite $m_{{q}_1,\dots,{q}_t}^{p_1,\dots,p_t}(1\wedge\textit{sh}_{i,p_j-i}\wedge1)$
in $\overline{(sA)^{\otimes p'}}$
 is zero.
An element in the image is a sum of terms, indexed by all
straight $(\ell_k,{q}_k,r_k,p_k)_{k=1}^t$
shuffles $\sigma$ and all $(i,p_j-i)$-shuffles $\tau$.
For each $(\sigma,\tau)$, let $x$ and $y$ be the lengths of the subintervals
of $[1,i]$ and $[i+1,p_j]$ given by the
intersections with $\tau^{-1}[\ell_j+1,\ell_j+q_j]$.
Now the terms for given $x,y\geq1$ may be grouped together into terms containing a factor
$m_{{q}_1,\dots,{q}_t}
(1\wedge\textit{sh}_{x,y}\wedge1)$, and hence these vanish.
Collecting the terms for which $x$ or $y$ is zero gives elements in the image of the straight shuffle extensions $m_{q_1,\dots,x,y,\dots,q_t}^{p_1,\dots,i,p_j-i,\dots,p_t}$, which by the previous remark also vanish.
\end{enumerate}\end{remas}

\begin{exams}{~}\par

\begin{itemize}
\item[$\diamond$]
For $t=1$ we observe that $m'=\sum_{q\leq p} m_q^p$ reduces to the usual formula for the extension of
a degree $-1$ map
$\sum m_q:\Lie^c(sA)\to sA$ as a coderivation on the cofree Lie coalgebra $\Lie^c(sA)$.

\item[$\diamond$]
For $t=2$ a result of
Fresse~\cite[Lemma 1.3.5]{Fresse06} says that a Lie bracket $m_{1,1}:sA\wedge sA\to sA$
extends to a {\em shuffle bracket} $\Lie^c (sA)\wedge \Lie^c (sA)\to \Lie^c (sA)$. The straight shuffle extension
$\sum m^{p_1,p_2}_{1,1}$ coincides with this shuffle
bracket.
\end{itemize}
\end{exams}

Now the cofree $\G^{\ac}$-coalgebra on a graded space $A$ is given by
\begin{align*}
\G^{\ac}(A)
&=s^{-2}\Com^c(s\Lie^c(sA))
\cong\bigoplus_{p_1\leq\dots\leq p_t}
s^{t-2}\bigl(\overline{(sA)^{\otimes p_1}}
\wedge\dots\wedge
\overline{(sA)^{\otimes p_t}}\bigr),
\end{align*}
so that a degree $-1$ map $m:\G^{\ac}(A)\to A$ is equivalent to  a family of maps
\begin{align}\nonumber 
m_\mathsf{q}=m_{q_1,\dots,q_t}\colon
\overline{(sA)^{\otimes q_1}}\wedge\dots\wedge \overline{(sA)^{\otimes q_t}}
\longrightarrow sA,\qquad 
\quad t,\,q_k\geq1,
\end{align}
of degrees $t-2$.
\begin{lemm}
Let $m\colon \G^{\ac}(A)\to A$ be  a degree $-1$ map  and let  $m_\mathsf{q}\colon
\overline{(sA)^{\otimes q_1}}\wedge\dots\wedge \overline{(sA)^{\otimes q_t}}
\longrightarrow sA$  be the associated equivalent family of maps.

The map
$$
m'=\sum m^\mathsf{p}_\mathsf{q}\ \colon\  \G^{\ac}(A)\to \Lie^c(sA)
$$
defined by the straight shuffle extensions is the unique coderivation with respect to the Lie coalgebra structures which
extends $m$. Moreover the map
  $\widetilde m:\G^{\ac}(A)\to \G^{\ac}(A)$
defined by $$
\widetilde m(w_1\wedge\dots\wedge w_t)=\!\!\!\!
\sum_{\substack{I\sqcup J=\{1,\dots,t\}
}}
\!\!\!\!
(-1)^{\varepsilon+\varepsilon'+\varepsilon''}\,
s^{|J|}
m'
(w_{{i_1}}\wedge\dots\wedge w_{{i_{|I|}}})\wedge w_{{j_1}}\wedge\dots\wedge w_{{j_{|J|}}},
$$
is the unique extension of $m'$ as a coderivation of cocommutative coalgebras.

Finally $\widetilde m:\G^{\ac}(A)\to \G^{\ac}(A)$ is the unique extension of $m$ to a degree $-1$ coderivation of $\G^{\ac}$-coalgebras.
\end{lemm}
\begin{proof}
Let us introduce the following non-commutative version of the operad $\G$: the operad $\mathbb{G}$ is the operad
encoding algebras defined by a degree $0$  associative product $\star$ and by a degree $1$
Lie bracket  $\langle\; ,  \, \rangle$   satisfying the  Leibniz relation
$$\langle\, \textrm{-} , \textrm{-}\star \textrm{-}\, \rangle\ =\
(\langle\,\textrm{-}, \textrm{-}\,\rangle\star \textrm{-})\
+ \ (\textrm{-}\star \langle\,\textrm{-},\textrm{-}\,\rangle).(12),   $$

This operad was already introduced in \cite{Merkulov04} where it was shown that it does not satisfy the distributive law assumption. Hence it is actually isomorphic to a non-trivial quotient of $\Ass \circ \Lie_1$.

Recall that the morphism of operads $\Ass \epi \Com$ induces the morphism between the Koszul dual cooperads
$\pi\, :\, \Ass^{\ac} \epi \Com^{\ac}$, equivalent to the morphism of cooperads $ \Ass^c \epi \Lie^c$ after desuspension.
 These morphisms extend to a morphism of operads
$ \mathbb{G} \epi \G$ and to a morphism of cooperads
$\Pi \, : \, \mathbb{G}^{\ac} \epi
\G^{\ac}$ as follows. The Koszul dual operad of $\G$ is $\G$ itself \cite{GetzlerJones94} and the Koszul dual operad
$\mathbb{G}^!$ of $\mathbb{G}$ is the operad generated by an associative product $\ast$ of degree $1$ and by a
commutative and associative product $\bullet$ of degree $0$, which satisfy the left and right Leibniz relations:
\begin{eqnarray*}
\big(\textrm{-} \ast (\textrm{-} \bullet \textrm{-})\big)=
\big((\textrm{-} \ast \textrm{-}) \bullet \textrm{-}\big)+
\big(\textrm{-} \bullet (\textrm{-} \ast \textrm{-})\big).(12),\\
\big((\textrm{-} \bullet \textrm{-}) \ast \textrm{-}\big) = \big((\textrm{-} \ast \textrm{-}) \bullet \textrm{-}\big).(23)
+\big(\textrm{-} \bullet (\textrm{-} \ast \textrm{-})\big).
\end{eqnarray*}
A $\mathbb{G}^!$-algebra induces a $\G^!=\G$-algebra by anti-symmetrizing the associative product to define a Lie bracket. This functor is in one-to-one correspondence with a morphism of operads $\G^! \to \mathbb{G}^!$, which induces the morphism of cooperads $\Pi \, : \, \mathbb{G}^{\ac} \epi
\G^{\ac}$ by linear duality and suspension. Once again, the operad $\mathbb{G}^!$ does not satisfy the distributive law assumption, but $\mathbb{G}^{\ac}$ is a (non-trivial) sub-$\Sy$-module of $S^c\Com_1^c \circ S^c\Ass^c$ and the morphism of cooperads $\Pi \, : \, \mathbb{G}^{\ac} \epi \G^{\ac}$ is equal to the composite
$$\mathbb{G}^{\ac} \mono S^c\Com_1^c \circ S^c\Ass^c  \epi  S^c\Com_1^c \circ S^c\Lie^c \cong  \G^{\ac}.$$

The unique Lie coderivation which extends $m\colon \G^{\ac}(A)\to A$ is equal to the composite:
$$ \G^{\ac}(A)\xrightarrow{\widetilde{\Delta}_{\G^{\ac}}(A)} S^c\Lie^c\circ (A;\G^{\ac}(A))
\xrightarrow{S^c\Lie^c\circ (A;m)}  S^c\Lie^c(A), $$
where $\widetilde{\Delta}_{\G^{\ac}} \, : \,  \G^{\ac}\to S^c\Lie^c\circ (I;\G^{\ac})$ is the part
of the decomposition map
$\G^{\ac}\to \G^{\ac}\circ\G^{\ac}$ of the cooperad $\G^{\ac}$ composed
only by elements of $S^c\Lie^c$ on the left hand
side of the composite product $\circ$ and by identities $I$ and only one element from $\G^{\ac}$ on the right hand side.
The morphisms of cooperads $\pi \, :\, \Ass^{\ac} \epi \Lie^{\ac}$ and
$\Pi \, :\, \mathbb{G}^{\ac} \epi \G^{\ac}$ induce the following commutative diagram

$$\xymatrix@C=50pt{\mathbb{G}^{\ac}\ar[r]^(0.35){\widetilde{\Delta}_{\mathbb{G}^{\ac}}}
\ar@{->>}[d]^{\Pi}  & S^c\Ass^c\circ (I;\mathbb{G}^{\ac}) \ar@{->>}[d]^{\pi\circ(I,\Pi)}\\
\G^{\ac} \ar[r]^(0.35){\widetilde{\Delta}_{\G^{\ac}}}& S^c\Lie^c\circ (I;\G^{\ac}).}$$

Therefore, it is enough to understand the formula on the level of $\mathbb{G}^{\ac}$ and then pass to the quotient
by the image of the shuffle maps. To understand how
$m (\Pi(A)) \, :\, \mathbb{G}^{\ac}(A) \to \G^{\ac}(A) \to A$ extends to a coderivation of coassociative coalgebras
$\mathbb{G}^{\ac}(A)\to S^c \Ass^c (A)$, we make the linear dual of the map $\widetilde{\Delta}_{\mathbb{G}^{\ac}}(A)$
explicit, up to signs and (de)suspensions, as follows.

Let $\mathbb{P}$ be the non-commutative analogue of the operad governing Poisson algebras; it is defined in the same
way as the operad $\mathbb{G}$ but with a degree $0$ Lie bracket. Its Koszul dual operad $\mathbb{P}^!$ of $\mathbb{P}$ is the operad modelling algebras endowed with an associative product $\ast$ and a commutative and associative product $\bullet$, both of degree $0$, which satisfy the left and right Leibniz relations aforementioned.

The underlying $\Sy$-module of the operad  $\mathbb{P}^!$ is a (non-trivial) quotient of the composite product
$\Com \circ \Ass$. Therefore the free $\mathbb{P}^!$-algebra on $A$ is a quotient of the free commutative
algebra on the free associative algebra on $A$ and any of its elements can be written $A^1\odot \cdots \odot A^t$, with
$A^1, \ldots, A^t \in \bar{T}(A):=\sum_{n\ge 1} A^{\otimes n}$, the tensor module on $A$. For any element $a\in A$, the Leibniz relations imply
\begin{eqnarray*}
&& a \ast  (A^1\odot \cdots \odot A^t)= \sum_{k=1}^t A^1\odot \cdots \odot (  a \otimes A^k ) \odot \cdots
\odot A^t, \ \textrm{and}\\
 && (A^1\odot \cdots \odot A^t)\ast a= \sum_{k=1}^t A^1\odot \cdots \odot (  A^k \otimes a) \odot \cdots
\odot A^t.
\end{eqnarray*}

Hence, the $\ast$-product $\Ass\circ (A,\mathbb{P}^!(A) ) \to \mathbb{P}^!(A)$ of elements of
$A$ with one of $\mathbb{P}^!(A)$ is given by the following formula, by the associativity of the product $\ast$.
$$ \ast\big(a'_1\otimes \cdots \otimes a'_r \otimes (A^1\odot \cdots \odot A^t) \otimes a''_1\otimes
\cdots \otimes a''_s\big)=
\sum (a'_{I_1}\otimes A^1 \otimes a''_{J_1}) \odot \cdots \odot (a'_{I_t}\otimes A^t \otimes a''_{J_t}),
$$
where the sum runs over the partitions $I_1\sqcup \ldots \sqcup I_t =[1, r]$ and
$J_1\sqcup \ldots \sqcup J_t =[1, s]$, with possibly empty sets, and where $a'_{I_k}$ is equal to
$a'_{i_1}\otimes \cdots \otimes a'_{i_l}$, with ${I_k}=\lbrace i_1 < \ldots <i_l \rbrace$ and similarly
for the $a''$. Finally, this map is the exact linear dual to the one given by the straight shuffles.

From this it follows
that the extension $\widetilde m$ is also a coderivation with respect to the cobracket,
and with respect to the coproduct by definition. Uniqueness holds since $\G^{\ac}(A)$ is cofree.
\end{proof}

The formula for $\widetilde m$ given by the straight shuffle extensions
$m_\mathsf{q}^\mathsf{p}$ enables us to expand the condition $m\circ \widetilde m=0$ for $m$ to define a codifferential on $\G^{\ac}(A)$. Hence we can give an explicit presentation of the axioms for
homotopy Gerstenhaber algebras:

\begin{prop}\label{explicit Ginfty}
A homotopy Gerstenhaber algebra is a dg module $(A, d_A)$ with a
family of maps
$$
m_{p_1,\dots,p_t}\colon
\overline{(sA)^{\otimes p_1}}\wedge\dots\wedge \overline{(sA)^{\otimes p_t}}
\longrightarrow sA,\qquad p_1\leq\dots\leq p_t,\quad t,\,p_k\geq1,
$$
of degree $t-2$, with $m_1=d_A$ and  which satisfy the following conditions for
$w_k\in \overline{(sA)^{\otimes p_k}}$,
$k=1,\dots,t$:
$$
\!\!\!\sum_{\substack{I\sqcup J=\{1,\dots,t\}\\[0.6mm]
I=\{i_1<\dots<i_r\}\neq\varnothing
\\[0.6mm]
J=\{j_1<\dots<j_{t-r}\}
}}
\!\!\!\!\!\!\!\!\!\!\!\!\!\!
(-1)^{\varepsilon+\varepsilon'+\varepsilon''}\,m_{p',p_{j_1},\dots,p_{j_{t-r}}}
\!\bigl(
\!\!\!\!\!\!
\sum_{q_1,\dots,{q}_r\geq1 \atop q_k\leq p_{i_k}}\!\!\!\!\!
m^{p_{i_1},\dots,p_{i_r}}_{{q}_{1},\dots,{q}_{r}}
(w_{{i_1}}\wedge\dots\wedge w_{{i_r}})\wedge w_{{j_1}}\wedge\dots\wedge w_{{j_{t-r}}}
\bigr)=0.
$$
Here $p'=1+\sum (p_{i_k}-{q}_{k})$, $\varepsilon''=r(t-r)$ and $(-1)^{\varepsilon+\varepsilon'}$  is the shifted Koszul sign
associated to the reordering of the $w_k$. 
\end{prop}

\subsection*{Homotopy Gerstenhaber algebra structure on the Hochschild cochain complex (Deligne conjecture)}

Let $(A, \mu)$ be an associative algebra. Its deformation complex is
given by the Hochschild cochain complex $CH^\bullet(A):=\Hom(A^{ \ot \bullet},
A)=\End_A$. To fit in with our homological degree convention the degree of an element of $\Hom(A^{\ot n}, A)$ is equal to $-n$. This space forms a non-symmetric operad, therefore the direct sum of its components is endowed with a Lie bracket $[f ,  g ]$ and non-symmetric brace operations $\{f\}\{g_1, \ldots , g_k\}$. More precisely, if we denote by $f \circ_i g$ the composite of the map $f$ with the map $g$ at the ${i}^{\textrm{th}}$-input, the first brace operation is defined by
 $$\{f\}\{g\}\;\;:=\;\;\sum_{i=1}^n (-1)^{(i-1)(m-1)} f \circ_i g\;\;=\;\; \sum_{i=1}^n (-1)^{(i-1)(m-1)}\!\!\!\!
\vcenter{\xymatrix@R=9pt@C=9pt@M=3pt{&&&& \\ &&
*+[F-,]{g}\ar@{-}[ul] \ar@{-}[ur]\ar@{-}[u] \ar@{-}[d]^{\!i} &&
\\&&*+[F-,]{f}\ar@{-}[urr]\ar@{-}[ull] \ar@{-}[ul]\ar@{-}[ur] \ar@{-}[d]&&
\\&&&&}}$$
for $f\in \Hom(A^{\ot n}, A)$ and $g\in \Hom(A^{\ot m}, A)$. Since it
is a pre-Lie operation, it induces a Lie bracket by
anti-symmetrization. With our degree convention, this Lie bracket is
an operation of degree $+1$. We refer the reader to
\cite{Gerstenhaber63} for details on this classical topic.

The boundary map of the Hochschild cochain complex is given by
$\partial_\mu(f):=[\mu, f]$, using $\mu\in \Hom(A^{\otimes 2}, A)$. One can also define the cup product using the second brace operation:

$$f\cup g = \lbrace \mu  \rbrace \lbrace f,g  \rbrace:=\mu \circ (f\otimes g)=  \vcenter{\xymatrix@R=10pt@C=10pt@M=3pt{&&&& \\ &*+[F-,]{f}\ar@{-}[ul]\ar@{-}[u]\ar@{-}[ur]&&*+[F-,]{g} \ar@{-}[ur]\ar@{-}[ul]  & \\
 &&*+[F-,]{\mu}\ar@{-}[ul]\ar@{-}[ur] \ar@{-}[d]&&
\\ &&&&}}$$

The cup product is associative and has degree $0$. In \cite{GerstenhaberVoronov95}, Gerstenhaber and Voronov introduced an operad generated by a binary (cup) product and general $(1+k)$-ary (brace) operations satisfying the same relations as the cup product and the brace operations on $CH^\bullet(A)$, see also \cite{GetzlerJones94}. They call any algebra over this operad, a ``homotopy Gerstenhaber'' algebra. However this operad is neither equal nor equivalent to the quasi-free Koszul resolution $\G_\infty$ of the operad $\G$ of Gerstenhaber algebras, which gives the notion of Gerstenhaber algebra up to homotopy.
Therefore, we choose to denote it $\GV$ after Gerstenhaber--Voronov and call the associated algebras, \emph{Gerstenhaber--Voronov algebras}. Hence, the previous constructions yield a morphism of operads $\GV \to \End_{CH^\bullet(A)}$.

Using the brace operations, Gerstenhaber proved in
\cite{Gerstenhaber63} that the commutativity of the cup product and
the Leibniz relation between the cup product and the Lie bracket hold,
but only up to homotopy. Therefore, the cohomology $HH^\bullet(A)$ carries a Gerstenhaber algebra structure. This led Deligne to ask if it was possible to lift this structure to a homotopy Gerstenhaber algebra structure on $CH^\bullet(A)$. This problem was solved by Tamarkin, among others, who proved the following theorem.

\begin{theo}[\cite{Tamarkin98}]\label{Tamarkin Ginfty}
There is a homotopy Gerstenhaber algebra structure on Hochschild complex $CH^\bullet(A)$ which lifts the Gerstenhaber algebra structure on its cohomology. Moreover, this structure factors through the operad $\GV$:
$$\G_\infty \to \GV \to \End_{CH^\bullet(A)}.$$
\end{theo}

Notice that this $\G_\infty$-structure is not completely explicit
since it relies on Etingof--Kazhdan quantisation of Lie bialgebra or
equivalently on Drinfeld's associators (see \cite{TamarkinTsygan00}).

\subsection{Homotopy quadratic BV-algebra}

Since $\qBV$ is a Koszul operad in the classical sense, the notion of homotopy quadratic BV-algebra
 may be defined by a slight modification of Proposition
\ref{explicit Ginfty}. For a dg module $A$,
\begin{align*}
\qBV^{\ac}(A)&=\KK[\delta]\otimes \G^{\ac}(A)
\\&\cong\bigoplus_{d=0}^\infty \bigoplus_{p_1 \leq \dots\leq p_t}
s^{2d+t-2}\bigl(\overline{(sA)^{\otimes p_1}}
\wedge\dots\wedge
\overline{(sA)^{\otimes p_t}}\bigr),
\end{align*}
where the $2d$-fold suspension arises from the power $\delta^d$.\\

\begin{prop}
A {\em homotopy quadratic BV-algebra} is a dg module $(A, d_A)$ together
with a family of maps
$$
m_{p_1,\dots,p_t}^d\colon
\overline{(sA)^{\otimes p_1}}\wedge\dots\wedge \overline{(sA)^{\otimes p_t}}
\longrightarrow sA,\qquad d\geq0,\,t\geq 1,\,p_1\leq\dots\leq p_t,
$$
of degree $2d+t-2$, with $m_1^0=d_A$.
The quadratic relations that these operations must satisfy are expressed by saying
that for each $d\geq 0$ and $t,\,p_1,\dots,p_t\geq1$,
the following expression is zero:
$$
\sum_{r=1}^t
\!\!\!\sum_{\substack{
d'+d''=d
\\[0.6mm]
I\sqcup J=\{1,\dots,t\}
}}
\!\!\!\!\!\!\!\!\!
(-1)^{\varepsilon+\varepsilon'+\varepsilon''}\,m^{d''}_{p',p_{j_1},\dots,p_{j_{t-r}}}
\!\bigl(
\!\!\!\!\!\!
\sum_{q_1,\dots,{q}_r\geq1 \atop q_k\leq p_{i_k}}\!\!\!\!\!
m^{p_{i_1},\dots,p_{i_r};d'}_{{q}_{1},\dots,{q}_{r}}
(w_{{i_1}}\wedge\dots\wedge w_{{i_r}})\wedge w_{{j_1}}\wedge\dots\wedge w_{{j_{t-r}}}
\bigr).
$$
Here $I=\{i_1<\dots<i_r\}$,
$\;J=\{j_1<\dots<j_{t-r}\}$, $\;p'=1+\sum(p_{i_k}-{q}_{k})$, $\;w_k\in
\overline{(sA)^{\otimes p_k}}$, and
$m^{p_1,\dots,p_r;d'}_{{q}_{1},\dots,{q}_{r}}$ is the straight shuffle extension
of $m^{d'}_{{q}_{1},\dots,{q}_{r}}$ to
$\overline{(sA)^{\otimes p_1}}\wedge\dots\wedge \overline{(sA)^{\otimes p_r}}$.
The sign $(-1)^{\varepsilon+\varepsilon'}$  is the shifted Koszul sign
associated to the reordering of the $w_k$,
and $\varepsilon''=r(t-r)$.
\end{prop}

\subsection*{Homotopy qBV-algebra structure on Hochschild cochain complex}

Let us pursue the example of the Hochschild complex of an associative algebra $A$. We suppose now that $A$ has a unit $1$. This extra datum allows us to define the following square-zero unary operator
$$\Delta(f)\;\;
:=\;\;\{f\}\{1\}\;\;=\;\;\sum_{i=1}^n (-1)^{(i-1)} f \circ_i 1\;\;=\;\; \sum_{i=1}^n (-1)^{(i-1)}\!\!\!\!\!\!\!\!\!
\vcenter{\xymatrix@R=10pt@C=10pt@M=3pt{ &&
1 \ar@{-}[d]^{i} &&
\\&&*+[F-,]{f}\ar@{-}[urr]\ar@{-}[ull] \ar@{-}[ul]\ar@{-}[ur] \ar@{-}[d]&&
\\&&&&}}. $$
The operator  $\Delta$ defines a map from $\Hom(A^{\otimes n}, A)$ to $\Hom(A^{\otimes n-1}, A)$, so $\Delta$ has degree $+1$ with our degree convention. It easy to check that the cup product and the brace operations
commute with $\Delta$. 

\begin{prop}\label{qBV Deligne}
There is a  homotopy qBV-algebra structure on the Hochschild complex $CH^\bullet(A)$ of any unital associative 
algebra $A$, which is given by  the homotopy Gerstenhaber algebra of Theorem~\ref{Tamarkin Ginfty} and the operator $\Delta$. 
\end{prop}

\begin{proof}
Consider the homotopy Gerstenhaber algebra structure provided by
Theorem~\ref{Tamarkin Ginfty}.
 The only other non-trivial higher operation is
$m^1_1=\Delta$; the other $m_{p_1,\ldots,p_t}^k$ are null for $k\ge
1$. To prove that these operations satisfy the relations of a homotopy
qBV-algebra, one has just to check that the $m_{p_1,\ldots,p_t}$
commute with $\Delta$. Since the homotopy Gerstenhaber algebra factors through
$\GV$ and since $\Delta$ commute with the operations coming from $\GV$, this concludes the proof.
\end{proof}

Notice that many of higher operations vanish in this $\qBV_\infty$-algebra structure (and indeed in the homotopy Gerstenhaber algebra structure of Theorem~\ref{Tamarkin Ginfty}). Notice also that  the operator $\Delta$ is well defined but vanishes on the cohomology level. (The reader may check that the degree $+2$ operator sending a cochain $f$ to $\lbrace f \rbrace \lbrace 1, 1\rbrace$ defines a homotopy between $\Delta$ and the zero map).


\subsection{Homotopy BV-algebras}

A {\em homotopy BV-algebra} is given by the same operations as a
homotopy quadratic BV-algebra, but with a linear perturbation of
the relations they have to satisfy.

\begin{theo}\label{Explicit Def BVinfty}
A \emph{homotopy BV-algebra} is dg module $(A, d_A)$ together
with a family of maps
$$
m_{p_1,\dots,p_t}^d\colon
\overline{(sA)^{\otimes p_1}}\wedge\dots\wedge \overline{(sA)^{\otimes p_t}}
\longrightarrow sA,\qquad t\geq1,\,d\geq0,\,p_1\leq\dots\leq p_t,
$$
of degree $2d+t-2$, with $m_1^0=d_A$.
The relations $\mathsf{R}^d_{p_1,\dots,p_t}$
that these operations must satisfy are expressed by saying
that for each $d\geq 0$, $\;t,\,p_1,\dots,p_t\geq1$ and $w_r=\overline{a^r_{[1,p_r]}}\in
\overline{(sA)^{\otimes p_r}}$,
the following expression is zero:
\begin{align*}&
\sum_{r=1}^t
\!\!\!\sum_{\substack{
d'+d''= d
\\[0.6mm]
I\sqcup J=\{1,\dots,t\}
}}
\!\!\!\!\!\!\!\!\!
(-1)^{\varepsilon+\varepsilon'+\varepsilon''}\,m^{d'}_{p',p_{j_1},\dots,p_{j_{t-r}}}
\!\bigl(
\!\!\!\!\!\!
\sum_{q_1,\dots,{q}_r\geq1 \atop q_k\leq p_{i_k}}\!\!\!\!\!
m^{p_{i_1},\dots,p_{i_r};d''}_{{q}_{1},\dots,{q}_{r}}
(w_{{i_1}}\wedge\dots\wedge w_{{i_r}})\wedge w_{{j_1}}\wedge\dots\wedge w_{{j_{t-r}}}
\bigr)\\
&+\sum_{r=1}^t\sum_{p'_r+p''_r=p_r}(-1)^{\varepsilon'''}
\,m^{d-1}_{p_1,\dots,p_{r-1},p'_r\!,p''_r,p_{r+1},\dots,p_t}
(w_1\wedge\dots\wedge\overline{a^r_{[1,p'_r]}}\wedge\overline{a^r_{[p'_r+1,p_r]}}\wedge\dots\wedge w_t)
\end{align*}
Here $I=\{i_1<\dots<i_r\}$,
$\;J=\{j_1<\dots<j_{t-r}\}$, $\;p'=1+\sum(p_{i_k}-{q}_{k})$, and
$\,m^{p_1,\dots,p_r;d'}_{{q}_{1},\dots,{q}_{r}}$ is the straight shuffle
extension of $m^{d'}_{{q}_{1},\dots,{q}_{r}}$ to
$\overline{(sA)^{\otimes p_1}}\wedge\dots\wedge \overline{(sA)^{\otimes p_r}}$.
The sign $(-1)^{\varepsilon+\varepsilon'}$  is the shifted Koszul sign
associated to the reordering of the $a_k$,
$\;\varepsilon''=r(t-r)\,$ and
$\,\varepsilon'''=|w_1|+\dots+|w_{r-1}|$.
\end{theo}

If $d=0$, or if each $p_k=1$, the linear terms are not present in the
relations. In particular:
\begin{itemize}
\item[$\diamond$] the operations
$m^0_{p_1,\dots,p_t}$ give $A$ the structure of a homotopy Gerstenhaber algebra,
\item[$\diamond$] the operations
$m^0_{1,\dots,1}$ give $sA$ the structure of a (strong) homotopy Lie algebra,
\item[$\diamond$] the operations
$m^0_{p}\;$ give $A$ the structure of a homotopy commutative algebra.
\end{itemize}
That is, there are inclusions of the operads
$\G_\infty$, $(\Lie_1)_\infty$ and ${\mathcal C}_\infty$ into $\BV_\infty$,
which are split by corresponding projections. (Recall that a homotopy commutative algebra or ${\mathcal C}_\infty$-algebra is an ${\mathcal A}_\infty$-algebra whose operations vanish on the image of the sum of the non-trivial shuffles, see \cite[Section~$13.1.13$]{LodayVallette09}).

\begin{rema}
One more naturally regards the structure of a homotopy BV-algebra $A$ in terms of maps on $A$ itself rather than on the suspension $sA$, and we may consider a homotopy BV-algebra as given by multilinear maps
$$
m_{p_1,\dots,p_t}^d\colon
A^{\otimes p_1}\otimes\dots\otimes A^{\otimes p_t}
\longrightarrow A,\qquad t\geq1, \, d\geq0, \,p_k\geq 1
$$
of degrees $2d+n+t-3$, where $n=\sum_{k=1}^t p_k$. These must be symmetric under permutation of the blocks $A^{\otimes p_k}$, vanish on the images of the shuffle maps on each block, and satisfy the above conditions $\mathsf{R}^d_{p_1,\dots,p_t}$ with the appropriate changes of sign.
\end{rema}

\subsection*{Homotopy BV-algebra structure on Hochschild cochain complex (cyclic Deligne conjecture)}

If we consider a Frobenius
algebra, i.e.\ a unital associative algebra endowed with a symmetric invariant non-degenerate bilinear form, then one can transfer Connes' boundary map from
Hochschild homology to define an operator $\Delta$ on cohomology \cite{Tradler02, Ginzburg06, Menichi07}. In this case, $HH^\bullet(A)$
carries a BV-algebra structure. The so-called Cyclic Deligne Conjecture amounts to proving that this structure can be lifted to a homotopy BV-algebra structure on the cochain level $CH^\bullet(A)$.

\begin{theo}[Cyclic Deligne conjecture]
Let $A$ be  a Frobenius algebra. There is a homotopy BV-algebra structure on its Hochschild cochain complex which lifts the BV-algebra structure on Hochschild cohomology and such that
$$\BV_\infty \qi C_\bullet(\fD) \to {\End}_{CH^\bullet(A)}$$
\end{theo}

\begin{proof}
Costello proved in
\cite{Costello07} that the Riemann sphere operad acts on
$CH^\bullet(A)$ and Kaufmann proved in
\cite{Kaufmann04} that $C_\bullet(\fD)$ acts on $CH^\bullet(A)$. Finally, we conclude the proof with Theorem~\ref{TCFTBV}.
\end{proof}

This conjecture was proved with various topological models in
\cite{Kaufmann04, TradlerZeinalian06, Costello07,
  KontsevichSoibelman06}. To the best of our knowledge, the operads involved in these proofs are not cofibrant, whereas the dg operad $\BV_\infty$ is. So this provides a canonical model for the cyclic Deligne conjecture.\\

Now we would like to conjecture, as in \cite{TamarkinTsygan00}, that
this homotopy BV-algebra structure on the Hochschild cochain complex  is formal under some assumptions on the algebra $A$, like $A=C^\infty(M)$ the algebra of smooth functions on a manifold. To phrase this conjecture, we need the notion of $\infty$-morphism for homotopy BV-algebras, which comes in Section~\ref{BV Homotopy theory}.

\subsection{Related definitions of homotopy BV-algebras in the literature}



Getzler showed in~\cite[Proposition~$1.2$]{Getzler94} that a BV-algebra
may be equivalently defined as an associative and commutative
dg algebra $(A,d_A,\bullet)$ equipped with a degree 1 square-zero
chain map $\Delta$
of {\em order $\leq$ 2}.
That is, the bracket defined by
\begin{align*}
\la a,b\ra 
&\;\;:=\;\;\Delta(a\bullet b)-(\Delta
a)\bullet b-(-1)^{|a|}
a\bullet (\Delta b)
\\
\intertext{satisfies the Leibniz relation}
0\;\;=\;\;\hdb{}3{a,b,c}&\;\;:=\;\;
\hdb{}{}{a,b\bullet c} -
\hdb{}{}{a,b}\bullet c -
(-1)^{|b|(|a|+1)} b\bullet \hdb{}{}{a,c} .
\end{align*}
More generally a notion of order may be
defined as follows.
\begin{defi}
Let $\Delta$ be a linear map on a commutative associative algebra $A$.
Define a sequence of  {\em higher brackets}
$\hdb\Delta n{\dots}:A^{\odot n}\to A$ of degree $|\Delta|\;$
by $\hdb\Delta1a:=\Delta(a)$ and, recursively,
\begin{align}
\nonumber
\hdb\Delta{n+1}{a_1,\dots,a_{n-1},a_n,a_n'}
&\;\;:=
\\&\makebox[-18ex]{}\nonumber
 \hdb\Delta n{a_1,\dots,a_{n-1},a_n\bullet a_n'}
 -\hdb\Delta n{a_1,\dots,a_n}\bullet a_n'
 -(-1)^\varepsilon a_n\bullet\hdb\Delta n{a_1,\dots,a_{n-1},a_n'}
\end{align}
where $\varepsilon =|a_n|(|\Delta|+|a_1|+\cdots+|a_{n-1}|)$.
Then $\Delta$  has {\em order $\leq n$} if $\hdb\Delta{n+1}{\dots}$ is identically zero.
\end{defi}
Several definitions of the order of a differential operator appear in
the literature~\cite{Akman97jpaa,Grothendieck67ega4,Koszul85asterisque},
and the sequence of brackets is  termed the {\em Koszul hierarchy}.
When needed, we will denote the higher brackets induced by an operator $\Delta$ by
$\hdbd\Delta{n}{\dots}$. In the context of commutative associative algebras the various definitions
coincide (see for example~\cite{AkmanIonescu08jhrs}) and may be restated as
\begin{align*}
\hdb\Delta n{a_1,\dots,a_n}\;\;\;=
\sum_{\substack{
I\sqcup J=\{1,\dots,n\}
\\[0.6mm]
|I|=r\geq1}}
\!\!\!
(-1)^{n-r+\varepsilon}\,
\Delta(a_{{i_1}}\bullet\dots\bullet a_{{i_r}})
\bullet a_{{j_1}}\bullet\dots\bullet a_{{j_{n-r}}},
\end{align*}
Here $\varepsilon$ is the Koszul sign associated to the
reordering of the $a_k$. It is now clear the brackets are graded
symmetric, and the brackets and order are well-defined even for
non-homogeneous $\Delta$.

%
\begin{prop}[\cite{Bering96,AkmanIonescu08jhrs,Cattaneo07,Kravchenko00,TVoronov03}]
If $\Delta$ is a square-zero operator of homogeneous degree $k$ on a commutative associative algebra
$A$, then the hierarchy of higher brackets $\{\hdbd{s^{-k-1}\Delta} n\dots\}$ associated to $s^{-k-1}\Delta$ defines a
$L_\infty$-algebra structure on $(s^{-1}A,s^{-k-1}\Delta)$, the desuspension of $A$ with $s^{-k-1}\Delta$ for differential.

Under this correspondence, a homogeneous square-zero operator $\Delta$ has order $\leq n$ if and only if the induced $L_\infty$-algebra structure on $s^{-1}A$ is an $L_n$-algebra, that is, the structure operations $\lbrace l_k\rbrace$ vanish for $k> n$.
\end{prop}

In case the operator $\Delta$ is not homogeneous we may consider the
following definition, see \cite{Kravchenko00}:
\begin{defi}
A \emph{commutative $\BV_\infty$-algebra} consists of a dg commutative algebra
$(A, d_A, \bullet)$ and a square-zero operator $\Delta$ whose homogeneous components $\Delta_k:A_*\to
A_{*+ k}$
satisfy:
\begin{itemize}
\item[$\diamond$] $\Delta_{-1}=d_A$,
\item[$\diamond$] $\Delta_k=0$ if $k$ is even or $k<-2$,
\item[$\diamond$] $\Delta_k$ has order $\le (k+3)/2$.
\end{itemize}
\end{defi}

This generalises the notion of BV-algebra, for which the
operator is just $\Delta=d_A+\Delta_1$, with $\Delta_1$ of order $\le2$,
though it conserves the strict associativity of the product.
For any $n\ge 0$, the operator $s^{-2n}\Delta_{2n-1}$ of order $\leq n$ induces an $L_n$ algebra on $s^{-1}A$ with differential $s^{-2n}\Delta_{2n-1}$.
The exact relation between the two definitions is the following.

\begin{prop}
Specifying a commutative $\BV_\infty$-algebra is equivalent to
specifying a homotopy BV-algebra in which all structure maps are
zero except
possibly $m_2^0$ and the operations $m^d_{1,\dots,1}$.
\end{prop}

The operad for commutative $\BV_\infty$-algebras is a quotient of
$\BV_\infty$, given by identifying generators
$s^{-1}\delta^d\otimes L_1\odot\dots\odot L_t$ to zero unless
each $L_i$ has arity 1, or unless $d=0$, $t=1$ and $L_1$ has arity 2.

\begin{proof}
Let $A$ be a homotopy BV-algebra in which
all structure maps are identically zero except $m_2^0$ and the maps
$m_{1,\dots,1}^d$. In this case  the only relations
$\mathsf{R}_{p_1,\dots,p_t}^d$ in which not all terms vanish are
$\mathsf{R}^0_2$, $\mathsf{R}^0_3$, the $\mathsf{R}^d_{1,\dots,1}$ and the $\mathsf{R}^d_{1,\dots,1,2}$.
First the relations $\mathsf{R}^0_2$ and $\mathsf{R}^0_3$ are equivalent to saying that $\bullet:=m^0_2$ is a chain map, which defines an associative and commutative product of degree $0$. Let us denote $\Delta_{2d-1}:=m^d_1$ the operator of degree $2d-1$ for $d\ge 0$, with $\Delta_{-1}=d_A$. The relations $\mathsf{R}^d_{1,\dots,1,2}$ are equivalent to
$$m^d_{1,\ldots,1}=\hdbd{\Delta_{2d+2n-3}}{n}{\dots}\quad \textrm{for} \quad d\ge 0,\, n\ge 1.$$
Therefore the operations $m^d_{1,\ldots,1}$ are completely characterised by the $\Delta_{2d+2n-3}$. Now the relation
$\mathsf{R}^d_{1,\dots,1,2}$ for $d=0$ says that $m^0_{1,\ldots,1}$ is a derivation with respect to $\bullet$, which is equivalent, by the above equality, to $\Delta_{2n-3}$ having order $\leq n$. We can consider the sum
$\Delta:=\Delta_{-1}+\Delta_1+\Delta_{3}+\cdots$, which is well defined by the definition of a homotopy BV-algebra. The relations $\mathsf{R}^d_1$ are equivalent for $\Delta$ to square to zero. So $A$ is a commutative $\BV_\infty$-algebra.

In the other direction, given a commutative $\BV_\infty$-algebra $A$, we define
$$m^0_2:=\bullet, \quad m^d_1:=\Delta_{2d-1},\quad m^d_{1,\ldots,1}=\hdbd{\Delta_{2d+2n-3}}{n}{\dots}$$
and the other $m^d_{p_1,\ldots, p_t}$ to be $0$. By the previous analysis, all the relations of a homotopy BV-algebra are satisfied, except for the $\mathsf{R}^d_{1,\dots,1}$, which is satisfied by the definition of the $m^d_{1,\ldots,1}$ in terms of higher brackets and because $\Delta^2=0$.
\end{proof}

Therefore in a commutative $\BV_\infty$-algebra, there is a family of $L_n$-algebras, one for each $n$, given up to suspension by
$$\lbrace \Delta_{2n-1}, \hdbd{\Delta_{2n-1}}{2}{\dots}, \ldots, \hdbd{\Delta_{2n-1}}{n}{\dots}\rbrace \quad
\textrm{or equivalently by} \quad
\lbrace m^n_1, m^{n-1}_{1,1}, \ldots, m^0_{1, \ldots, 1}\rbrace.$$

Tamarkin and Tsygan~\cite{TamarkinTsygan00}
 gave a more general definition of $\BV_\infty$-algebra structure.
This is also expressed in terms of higher order operators, not on $A$
itself but on the free Gerstenhaber algebra $\G(s^{-1}A^*)$
on the desuspension of the linear dual of $A$.
Using a dual formulation of Theorem \ref{4 def theo}, a homotopy Gerstenhaber algebra
structure is specified by a differential
$\Delta_{-1}$ of degree $-1$ on the free Gerstenhaber algebra $\G(s^{-1}A^*)$ when $A$ is finite dimensional. They define a homotopy BV-algebra
 structure extending this $\G_\infty$ structure as an odd square zero operator $\Delta$ on $\G(s^{-1}A^*)$ which
can be written as a sum of homogeneous components
$$
\Delta\;\;=\;\;\Delta_{-1}+\Delta_{\textrm{Lie}}+\Delta_{1}+\Delta_{3}+\Delta_{5}+\cdots,
$$
where  $\Delta_{2d-1}$ has degree $2d-1$ and order $\le d+1$ with respect to the commutative product of
$\G(s^{-1}A^*)$, for $d\ge 1$. The operator $\Delta_{\textrm{Lie}}$ is the unique order $\le 2$ operator on $\G(s^{-1}A^*)\cong \Com \circ \Lie_1(s^{-1}A^*)$, with respect to the commutative product, which extends the bracketing of two elements of $\Lie_1(s^{-1}A^*)$, the free $\Lie_1$-algebra on $s^{-1}A^*$. So it has degree $1$. They also require that each operator $\Delta_{1-2d}$ respects the decreasing filtration given by the length of $\Lie$-words in $\G(s^{-1}A^*)$.

The following Proposition was
observed in~\cite[Remark 2.8]{TamarkinTsygan00}, though Koszul duality
for inhomogeneous quadratic operads was not yet known.

\begin{prop}
Specifying a $\BV_\infty$-algebra structure $\Delta$ in the sense
of~\cite{TamarkinTsygan00} for which each $\Delta_{2d-1}$ has order
$\le 1$ and is a derivation with respect to the $\Lie_1$-algebra structure of $\G(s^{-1}A^*)$ is equivalent to specifying a homotopy BV-algebra structure in the sense of this paper.
\end{prop}

\begin{proof}
Let $(A, \Delta)$ be a $\BV_\infty$-algebra in the sense of \cite{TamarkinTsygan00} such that each
 $\Delta_{2d-1}$ has order $\le 1$ and is a derivation with respect to the $\Lie_1$-algebra structure of $\G(s^{-1}A^*)$. This is equivalent to requiring that each  $\Delta_{2d-1}$ is a derivation with respect the Gerstenhaber algebra structure on $\G(s^{-1}A^*)$. Considering the linear dual, each $\Delta_{2d-1}$ induces a degree $2d-1$ coderivation ${^t}\Delta_{2d-1}$ on the cofree Gerstenhaber coalgebra $\G^c(sA)$. Hence the sum $d:={}^t\Delta_{-1}+\delta^{-1}\otimes{}^t\Delta_{1}+\delta^{-2}\otimes{}^t\Delta_{3}+\delta^{-3}\otimes{}^t\Delta_{5}+\cdots$ is a generic degree $-1$ coderivation on $\qBV^{\ac}(A)\cong \KK[\delta]\otimes (s^{-1} \G^c(sA))$. The above definition of $\Delta_{\textrm{Lie}}$ and Lemma~\ref{dvarphi} show that the induced map $\delta^{-1}\otimes {}^t\Delta_{\textrm{Lie}}$ on $\qBV^{\ac}(A)$ is equal to the degree $-1$ coderivation $d_\varphi\circ A$. Finally, the degree $-1$ coderivation $d+d_\varphi\circ A$ squares to zero, which is equivalent to a homotopy BV-algebra structure by Theorem~\ref{4 def theo}.
\end{proof}

To conclude, the definition of homotopy BV-algebra given here lies between the definition of \cite{Kravchenko00}, which is given by an operad that is not a cofibrant replacement of the operad $\BV$, and the definition of \cite{TamarkinTsygan00}, which is actually an algebra over a properad (or a prop).\\

The Koszul resolution $\Po_\infty \qi \Po$ of an operad $\Po$ defines the notion of homotopy $\Po$-algebra. Recall that, in Koszul duality theory, the symmetry of the defining operations of $\Po$ remains strict in the homotopy version; only the relations between these operations being relaxed up to homotopy.

In \cite{SalvatoreWahl03}, Salvatore and Wahl began to extend the Koszul duality theory to operads in the category
of H-modules, where H is a cocommutative Hopf algebra. The Hopf algebra is meant to encode other group actions than the permutations of elements and it often comes from the homology algebra of a topological group. In the case of the circle $S^1$, we work with the algebra of dual numbers $D\cong H_\bullet(S^1)$ with the obvious coproduct. So the category of dg $D$-modules is the category of mixed chain complexes. With a proper action of $D$ on $\G$, a $\G$-algebra in dg  $D$-modules is a BV-algebra.

In this context, they defined the Koszul dual of a $D$-equivariant quadratic operad. One can continue much further and define the cobar construction of $D$-equivariant cooperads. Applying this construction to the Koszul dual of $\G$, one would get a notion of ``homotopy BV-algebra'' where the $\Delta$ operator will still strictly square to zero, this relation not being relaxed up to homotopy. This notion would be equivalent to a homotopy BV-algebra, as defined here, where only $m^1_1$ and the $m^0_{p_1,\ldots, p_t}$ would not be trivial. The only non-trivial relations would be
$\mathsf{R}_{p_1,\dots,p_t}^0$, which are the $\G_\infty$-algebra relations, and $\mathsf{R}_{p_1,\dots,p_t}^1$, which would say that the failure for $m^1_1$ to be a derivation with respect to $m^0_{p_1,\ldots, p_t}$ is
$$\sum_{r=1}^t\sum_{p'_r+p''_r=p_r}\pm
\,m^{0}_{p_1,\dots,p_{r-1},p'_r\!,p''_r,p_{r+1},\dots,p_t}.$$
Notice that this structure differs from that of Proposition \ref{qBV Deligne} exactly by this last term.


\section{Deformation theory of homotopy BV-algebras}\label{Def theo for BV}

We apply the methods of \cite{MerkulovVallette08I, MerkulovVallette08II} to the resolution of Section~$1$ in order to  give a Lie theoretic description of homotopy BV-algebras. This allows us to define the deformation theory, that is, the cohomology theory,  for BV-algebras and homotopy BV-algebras. We apply the general obstruction theory of algebras over a Koszul properad of Section~\ref{Obstruction Theory} to this case. This method allows us to prove an extended version of the Lian--Zuckerman conjecture: there exists a homotopy BV-algebra structure on vertex operator algebras which extends the operations defined by Lian--Zuckerman that induce a BV-algebra structure on homology. Finally, we show that the obstructions to lift a $G_\infty$-algebra structure to a homotopy BV-algebra structure live in negative even cohomology groups of the $G_\infty$-algebra.

\subsection{Convolution Lie algebra}\label{Convolution Lie BV}

Recall from \ref{Twisting morphisms} that for any dg operad $\Qo$, the space of equivariant maps $\g:=\Hom_\Sy({\BV}^{\ac}, \Qo)$ carries a dg Lie algebra structure, called the \emph{convolution Lie algebra}, such that any  morphism of dg operads $\gamma \,:\, \BV_\infty \to \Qo$ is equivalent to a twisting morphism $\BV^{\ac} \to \Qo$, still denoted $\gamma$.
Hence if $\Qo=\End_A$ for  a dg module $A$, twisting morphisms encode homotopy $BV$-algebra structures on $A$. After Section~\ref{ResolutionBV} and Section~\ref{Homotopy BV-algebra Explicit}, this gives a third equivalent definition of homotopy BV-algebras: homotopy BV-algebra structures are in one-to-one correspondence with Maurer--Cartan elements in the convolution Lie algebra $\g$, see Theorem~\ref{4 def theo}. \\

Since the Koszul dual cooperad $\BV^{\ac}$ is weight graded, the convolution algebra  $\g$ is also weight-graded, $\g\cong \prod_{n\ge 0} \g^{(n)}$, where $\g^{(n)}=\Hom_\Sy({{\BV}^{\ac}}^{(n)}, \Qo)$, from Proposition~\ref{graded convolution Lie}. The differential on $\Qo$ induces a weight $0$ differential $\pa_0$ on $\g$ and the internal differential $d_\varphi$ on $\BV^{\ac}$ induces a weight $+1$ differential $\pa_1$ on $\g$.
This result is used in Section~\ref{Vertex Algebras} to endow vertex algebras with a homotopy BV-algebra structure. \\

By Proposition~\ref{graded relative convolution Lie}, we know that the Koszul dual cooperad $\BV^{\ac}$ admits a relative grading ${\BV^{\ac}}^{[n]}$ such that ${\BV^{\ac}}^{[0]}=\G^{\ac}$ and Proposition~\ref{form BV dual} shows that $\BV^{\ac}\cong \KK[\delta]\circ \G^{\ac}$ as graded  $\Sy$-modules. In this case, the relative grading is equal to the power of $\delta$. So $\BV^{\ac}$ is an extension of the Koszul dual cooperad $\G^{\ac}$ of the Gerstenhaber operad. This proves that the convolution Lie algebra governing homotopy BV-algebra structures is a formal extension of the convolution Lie algebra governing homotopy Gerstenhaber algebras, twisted by an extra differential.

\begin{prop}\label{Convolution BV formal extension}
The convolution dg Lie algebra $\g:=\Hom_\Sy({\BV}^{\ac},
\Qo)$ is isomorphic to
$\g^\G[[\hbar]]:=\g^\G\otimes \KK[[\hbar]]$, where $\g^\G$ is the convolution dg Lie algebra $\Hom_\Sy({\G}^{\ac}, \Qo)$ associated to the Gerstenhaber operad $\G$. The formal parameter $\hbar$ has degree $-2$. The differential $\pa$ on $\g^\G[[\hbar]]$ is the sum of two terms $\pa=\pa_0+\pa_1$, where $\pa_0$ is the derivation freely extended from that of $\g^\G$ and where $\pa_1$ raises the power of $\hbar$ by $1$.
\end{prop}

\begin{proof}
By Proposition~\ref{form BV dual}, we have
$$\Hom_\Sy(\BV^{\ac}, \Qo)\cong \Hom_\Sy(\KK[\delta]\circ \G^{\ac}, \Qo)\cong \Hom_\Sy(\G^{\ac}, \Qo)\otimes \KK[[\hbar]],$$ with $\hbar=\delta^*$. By the cooperad structure on $\BV^{\ac}$, this is an isomorphism of Lie algebras. The result about the differentials is then straightforward.
\end{proof}

\subsection{Deformation complex} 

Section~$2$ of \cite{MerkulovVallette08II}, applied here, defines the deformation complex of morphisms from the  operad $\BV_\infty$. In this section, we present only the case $\Qo=\End_A$ for simplicity, but the general case is similar. Consider the convolution Lie algebra $\g=\Hom_\Sy({\BV}^{\ac}, \End_A)$ whose  twisting morphisms correspond to homotopy BV-algebra structures. Once given such a twisting morphism $\gamma$, we consider the
\emph{twisted differential} on $\g$
defined by
$$\partial^\gamma(f):=\partial(f)+[\gamma, f]=d_{\End_A} \circ f - (-1)^{|f|}f \circ d_\varphi +[\gamma, f].$$

\begin{defi}[Deformation complex]
The chain complex $\g^\gamma:=\big(\Hom_\Sy({\BV}^{\ac}, \End_A), \partial^\gamma\big)$ is called the \emph{deformation complex} of the homotopy BV-algebra structure $\gamma$ on $A$.
\end{defi}

This defines the ``cohomology of $A$ with coefficients in itself'', also called the \emph{tangent homology}, for any BV-algebra or homotopy BV-algebra. As usual, these (co)homology groups are the obstructions to formal deformations of the BV-algebra or homotopy BV-algebra structure $\gamma$ as in \cite{Gerstenhaber64}. To study the formal deformations, we use the following Lie algebra structure on the deformation complex.

\begin{prop}
The Lie bracket $[\;,\,]$ induces a \emph{twisted} dg Lie algebra structure on
$\g^\gamma$. This dg Lie algebra is isomorphic to $$\g^\gamma\cong \big(\g^\G[[\hbar]], \pa_0+\pa_1+[\gamma, -], [\;,\,] \big).$$
\end{prop}

\begin{proof}
The differential of any dg Lie algebra twisted by a Maurer--Cartan element always induces a twisted dg Lie algebra structure. The final isomorphism comes from Proposition~\ref{Convolution BV formal extension}.
\end{proof}

Hence the underlying space defining the cohomology of homotopy BV-algebras is the formal extension of that defining the cohomology of homotopy Gerstenhaber algebras. We apply this result in Section~\ref{Obstruction BV}.

\subsection{Obstruction theory applied to topological vertex operator algebras}\label{Vertex Algebras}

In parallel to Theorem~\ref{TCFTBV}, we prove that a large class of vertex operator algebras carry a homotopy BV-algebra structure. To do so, we describe an obstruction theory for homotopy BV-algebras. More precisely, Lian and Zuckerman described in \cite{LianZuckerman93} operations acting on a topological vertex operator algebra, which induce a BV-algebra structure on the BRST (co)homology (the underlying homology). In this paper, they ask whether this structure can be lifted off-shell (on the chain level) with ``\emph{higher homotopies}'' \cite[page 638]{LianZuckerman93}. We fully prove this Lian--Zuckerman conjecture here: there exists an explicit homotopy BV-algebra structure on topological vertex operator algebras which extends the operations defined by Lian and Zuckerman and which induces the BV-algebra structure on the BRST (co)homology. 
Huang has also investigated operadic formulations of vertex operator algebras, in terms of (partial) algebras over the Riemann sphere and framed little {\disc}s operads \cite{Huang97,Huang2003}. \\

The weight filtration on $\BV^{\ac}$ induces a
decomposition of the convolution Lie algebra
$$\g=\prod_{n\ge 0} \g^{(n)}
,\qquad \g^{(n)}=\Hom_\Sy({{\BV}^{\ac}}^{(n)},
\End_A) $$
and homotopy BV-algebra structures on a dg module $(A,d_A)$ are in
one-to-one correspondence with degree $-1$ elements
$$\gamma=\gamma_1+\gamma_2+\cdots\, :\, {\BV^{\ac}} \to \End_{A}
,\qquad \gamma_n:{{\BV}^{\ac}}^{(n)}\to \End_{A},
$$
satisfying the sequence of Maurer--Cartan equations
\begin{align*}\tag{$\textrm{MC}_n$}
\partial_0(\gamma_n) + \partial_1(\gamma_{n-1}) +  \frac{1}{2}\sum_{k+l=n} [\gamma_k, \gamma_l]=0,
\end{align*}
for $n\ge 1$. Notice that here $\partial_0(\gamma_n)$ is equal to $\pa_A (\gamma_n):=d_{\End_A}\circ \gamma_n$.
 See Section \ref{Obstruction Theory} for compete details in the general case of
quadratic-linear operads.


\begin{theo}\label{Obstruction homotopyBV} Let $(A,d_A)$ be a dg module and suppose we are given a commutative
 product of degree $0$, a skew-symmetric bracket of degree $1$ and a
 unary operator of degree $1$ such that the differential $d_A$ is a
 derivation with respect to them, that is, a morphism of $\Sy$-modules
 $$\gamma_1 \, :\, {\BV^{\ac}}^{(1)} \to \End_A$$ of degree $-1$
 satisfying $\pa_A (\gamma_1):=d_{\End_A}\circ \gamma_1=0$. If
 $\mathrm{H}_{-2}(\Hom_{\Sy}({\BV^{\ac}}^{(n)}, \End_A), \pa_A)=0$ for
 all $n\ge 2$, then this structure extends to a homotopy BV-algebra
 structure on $A$.
\end{theo}
\begin{proof}
 By Theorem~\ref{Obstruction operations to full}.\end{proof}

\medskip

We can define a homotopy BV-algebra structure on a topological vertex operator algebra, TVOA for short, closely following the $\G_\infty$-case described in detail in \cite{GalvezGorbounovTonks06}. Consider $(A,\mbox{-}_{(i)}\mbox{-},Q,G,L,F)$ a TVOA, bigraded as usual by the conformal weight $s$ and the fermionic grading $r$
$$A=\bigoplus_{s\in\ZZ} A^s=\bigoplus_{r,s\in \ZZ} A_r^s,$$
see \cite[Definition 2.4]{kivozu} or  \cite[Section 5.9]{kac} for the full definition. We will assume here the conformal weight is always non-negative, that is, $s\in\NN$. The TVOA axioms include the supercommutator relation
$$
G_{(1)}{Q}_{(0)}x+{Q}_{(0)}G_{(1)}x = sx  \text{ if } x\in A^s.
$$
For $s>0$ the complex $(A^s,d_A)$ is therefore contractible, with contracting homotopy $\frac1sG_{(1)}$.

Lian and Zuckerman~\cite{LianZuckerman93}  considered  the algebraic structure on $A$ given by
$$
x\centerdot y:=x_{(-1)}y,\qquad\qquad \{x,y\}:=(-1)^{|x|}(G_{(0)}x)_{(0)}y,\qquad\qquad \Delta x:=G_{(1)}x,
$$
and showed these operations induce a BV-algebra structure on the BRST cohomology $H_\bullet A$ with respect to the differential operator $d_A:={Q}_{(0)}$. They gave explicit first chain homotopies to make up for the failure  of the $\G_\infty$ structure to hold on $A$ itself, and conjectured the existence of a family of higher homotopies \cite[page 638]{LianZuckerman93}. Note that all the operations are of conformal weight zero and it is the fermionic grading that gives us the homological degree.

\begin{theo}\label{vertexBVinfty}
A topological vertex operator algebra $A$, with conformal weights in $\NN$, carries an explicit homotopy BV-algebra structure homogeneous with respect to the conformal weight, which extends the above operations in conformal weight zero and which induces the Lian--Zuckerman BV-algebra structure on the BRST cohomology $H^\bullet A$.
\end{theo}

\begin{proof}
The first key argument, going back to \cite[page 623]{gms-gerbes-2}, is the following. When a TVOA is concentrated in non-negative conformal weight, the (weight zero) operations $\textrm{-} \centerdot \textrm{-}$, 
$\{\textrm{-} ,\textrm{-}\}$ and $\Delta$ endow $A^0$ with a BV-algebra structure. The first homotopies for the BV relations described in \cite[(2.14), (2.16), (2.23), (2.25), (2.28)]{LianZuckerman93} involve terms of the form $v_{(i)}t$, where $v, t \in A^0$ and $i\ge 0$. Since the product ${}_{(i)}$ has weight $-i-1$, all these term vanish. 
In particular, the product and bracket of Lian and Zuckerman are already appropriately (skew) symmetric here.

Then we consider the morphism of dg $\Sy$-modules  $\gamma_1 \, :\, {\BV^{\ac}}^{(1)} \to \End_{A}$ of degree $-1$ by taking the product and bracket defined by the graded (skew-)symmetrizations of $\centerdot$ and $\{\, ,\}$. The image of the arity $1$ summand is equal to the operator defined in conformal weight $s=0$ by $G_{(1)}$ and by the zero operator elsewhere. By the preceding remark, the extension/restriction
$$\bar \gamma \, :\, \BV^{\ac} \to \End_{A^0}\,.$$
of $\gamma_1$ to $\BV^{\ac}$ on the source and  to $\End_{A^0}$ on the target is a twisting morphism.


Next we use the obstruction theory of Theorem~\ref{Obstruction n->n+1:GENERAL} to extend $\bar \gamma$ beyond $A^0$.
Since we are looking for a homotopy BV-algebra structure which respects the conformal weight, we already know the image 
of all the operations on $A^0$. Let us denote $\bar A:=\bigoplus_{s\ge 1} A^s$. We work by induction on the weight grading $n\ge 1$ of the twisting morphism $\gamma_1 +\cdots + \gamma_n+\cdots$ to define its image on $\bar A$. The induction is initiated by the definition of $\gamma_1$. Using the following isomorphisms,
$$\Hom_\Sy(\BV^{\ac}, {\End}_A)\cong \Hom_\KK(\BV^{\ac}(A), A)\cong  \Hom_\KK(\BV^{\ac}(A), A^0) \oplus  \Hom_\KK(\BV^{\ac}(A), \bar A)$$
we see that the component of $\gamma_{n+1}$ on the first direct summand is null, for $n\ge 1$. Since $\frac1sG_{(1)}$ is a contracting homotopy on  the chain complex $(A^s,d_A)$ for each $s>0$, the second direct summand is acyclic with the explicit contracting homotopy  $\frac1s {G_{(1)}}$ on $(\Hom_\KK(\BV^{\ac}(A), A^s), \partial_A)$. We can then conclude the proof by Theorem~\ref{Obstruction homotopyBV}. In fact, from the more general Theorem~\ref{Obstruction n->n+1:GENERAL}, we know that the obstruction to the existence of $\gamma_{n+1}\in \Hom_\KK({\BV^{\ac}}^{(n+1)}(A), \bar A)$ is the homology class of
$$\widetilde\gamma_{n+1}:=\partial_1(\gamma_n)
+\frac{1}{2}\sum_{k+l=n+1\atop k,l \ge 1}[\gamma_k, \gamma_l],$$
which vanishes here thanks to the contracting homotopy $G_{(1)}$.
Let us denote by $\overline \End_A$ the summand of $\End_A$ made up of operations landing in $\bar A$, and by $H$ the contracting homotopy on  $(\Hom_\KK(\BV^{\ac}(A), \bar A), \partial_A)$ induced by the $\frac1s {G_{(1)}}$ on each $A^s$.
Therefore the explicit image of $\gamma_{n+1}$ into $\overline \End_A$ is given by
$$\textrm{proj}_{\overline \End_A} \circ \gamma_{n+1}=H\big( \partial_1(\gamma_n)
+\frac{1}{2}\sum_{k+l=n+1\atop k,l \ge 1}[\gamma_k, \gamma_l]\big). $$
\end{proof}

Following the method of \cite{HirshMilles10}, one can go even further and define the notion of homotopy \emph{unital} BV-algebra. In the same way, a TVOA with non-negative conformal weight should carry such a structure, with strict homotopy unit. Moreover, we conjecture that the converse should be true. Namely, the data of a TVOA with non-negative conformal weight should be equivalent to a strict unital homotopy  BV-algebra. Here is a table, comparing these two algebraic structures, which sustains this conjecture. 

\begin{center}
\begin{tabular}{|c|c|
}\hline
TVOA & strict $uBV_\infty$-algebra\\
\hline
$\textrm{-} {}_{(-1)} \textrm{-}$ & $\textrm{-} \bullet \textrm{-}$\\
$Y : A \to \End_A [[z^{\pm 1}]]$ & higher binary homotopies\\
vaccum $|0\rangle$ & unit $1$\\
operators $L,F$ & bigraded decomposition $V^s_r$\\
operator $Q$ & differential $d$\\
operator $G$ & operator $\Delta$ and higher homotopies\\
\hline
\end{tabular}
\end{center}

This would generalize the well-known fact that the data of a vertex algebra with trivial operations $\textrm{-} {}_{(n)} \textrm{-}\equiv 0$, for $n \ge 0$ is a equivalent to a unital commutative dg algebra \cite[Section~$1.4$]{kac}. \\

We may also associate to any smooth Calabi--Yau $n$-manifold $M$ a sheaf of homotopy BV-algebras,
in a manner parallel to \cite[Section 5]{GalvezGorbounovTonks06} for homotopy Gerstenhaber lgebras.
One first recalls from \cite{gms-gerbes-2,gms-gerbes3} that two particular TVOAs, termed the \emph{chiral de Rham complex} and the \emph{complex of chiral vector fields}, may be constructed from the algebra of functions on $\CC^n$. By the result above, we then have corresponding homotopy BV-algebra structures locally, for each point $x\in M$. Now using \cite[Theorem 4.2]{msv-chiral} one checks, if the first Chern class $c_1(M)$ vanishes, these local structures may be glued to produce a sheaf on $M$.

More generally, this result should be applicable to Chiral algebras, as in Beilinson--Drinfeld \cite{BeilinsonDrinfeld04}, in the Geometric Langlands Programme.

\subsection{Relative obstruction theory}\label{Obstruction BV}

In this section, we show that the obstructions  to  lifting a $G_\infty$-algebra structure to  a homotopy BV-algebra structure already live in particular cohomology groups of the $G_\infty$-algebra.\\

Proposition~\ref{Convolution BV formal extension} proved that the convolution Lie algebra which governs homotopy BV-algebras is a formal extension of that of homotopy Gerstenhaber algebras: $\g\cong\g^\G[[\hbar]]$. Given a homotopy Gerstenhaber algebra structure $\gamma_0$ on $A$, the obstructions to lifting it to a homotopy BV-algebra structure are given by the following proposition.

\begin{prop}\label{Obstruction n->n+1 BV}
Let $\gamma = \gamma_0 + \gamma_1 + \cdots + \gamma_n \in  \prod_{k=0}^n\g^\G \otimes \hbar^k\KK$ be an element which satisfies the $\emph{(MC)}_k$-equations up to $k=n$ in $(\g, \partial)$. We consider  $$\widetilde{\gamma}_{n+1}:=\partial_1(\gamma_n) + \frac{1}{2}\sum_{k+l=n+1\atop k,l \ge 1}[\gamma_k, \gamma_l].$$
\begin{enumerate}
\item In $\g^\G\otimes \hbar^{n+1}\KK$, we have  $\partial^{\gamma_0}\left(  \widetilde{\gamma}_{n+1}  \right)=0$,
that is, $\widetilde{\gamma}_{n+1}$ is a cycle of degree $-2$.
\item There exists an element $\gamma_{n+1}\in \g^\G\otimes \hbar^{n+1}\KK$ such that $\gamma_0 + \gamma_1 + \cdots + \gamma_{n+1}$ satisfies the $\emph{(MC)}_k$-equations up to $k=n+1$ in $(\g, \partial)$ if and only if
the class of $\widetilde{\gamma}_{n+1}$ in $\mathrm{H}_{-2}(\g^\G\otimes \hbar^{n+1}\KK,  \partial^{\gamma_0})$ vanishes.
\end{enumerate}
\end{prop}

\begin{proof}
 It is the application of Theorem~\ref{Obstruction n->n+1:GENERAL} to the convolution Lie algebra of Proposition~\ref{Convolution BV formal extension}.
\end{proof}

Concretely, Proposition~\ref{Obstruction n->n+1 BV} applies as follows.

\begin{theo}
Let $\gamma_0$ be a $G_\infty$-algebra structure on a dg module $A$. If the negative  even cohomology groups of the $G_\infty$-algebra $A$ vanish, that is, $\mathrm{H}_{-2n}(\Hom_{\Sy}(\G^{\ac}, \End_A), \partial^{\gamma_0})=0$ for $n\ge 2$, then this structure can be extended to a homotopy BV-algebra structure.
\end{theo}

\begin{proof}
For $n\ge 1$, the chain complex $(\g^\G\otimes \hbar^{n+1}\KK,  \partial^{\gamma_0})$ is isomorphic to $(\g^\G,  \partial^{\gamma_0})\otimes ( \hbar^{n+1}\KK, 0)$, with $\hbar$ of degree $-2$.
\end{proof}

In other words, the negative even cohomology groups of a $G_\infty$-algebra measure the obstructions to lift it to a homotopy BV-algebra structure. This result comes from the fact that the operator $\Delta$ in the definition of a $BV$-algebra is induced by the topological action of the circle, see \cite{Getzler94, Getzler94bis}.

The same theorem holds for general morphisms of dg operads $\G_\infty\to \Qo$, and provides obstructions to lifting them to morphisms $\BV_\infty \to \Qo$.

\section{Homotopy theory of homotopy BV-algebras}\label{BV Homotopy theory}

The previous resolution of the operad $\BV$ in terms of a small explicit dg cooperad $\BV^{\ac}$ allows us to define the homotopy theory for BV-algebras and homotopy BV-algebras. Following the general methods of \cite{GetzlerJones94}, refined in Appendix~\ref{Homotopy Th General} to apply to inhomogeneous quadratic Koszul operads, we introduce the bar and cobar constructions between BV-algebras and $\BV^{\ac}$-coalgebras. We extend this bar construction to homotopy BV-algebra which allows us to define the notion of $\infty$-morphisms. Finally we show how to transfer homotopy BV-algebra structures across homotopy equivalences. Applied to the homology of a homotopy BV-algebra $A$, this result defines the Massey product  of $H(A)$, which contains the homotopy data of $A$.

\subsection{Bar and cobar constructions}

The morphism of dg operads $\BV_\infty=\Omega \BV^{\ac} \to  \BV$ is equivalent to a twisting morphism $\kappa\, :\, \BV^{\ac}\to \BV$  by Theorem~$2.17$ of \cite{GetzlerJones94}, see also Theorem~\ref{4 def theo}. This twisting morphism is equal to the composite
$$\kappa\, : \, \BV^{\ac}=\qBV^{\ac}
\epi sV \xrightarrow{s^{-1}}V \mono \BV.$$
 Proposition~$2.18$ of \cite{GetzlerJones94} applied to this twisting morphism  proves  the existence of a bar-cobar adjunction  between dg BV-algebras and dg $\BV^{\ac}$-coalgebras:
$$\mathrm{B}_\kappa \ : \ \textrm{dg}\ \BV\textrm{-algebras}\,  \leftrightharpoons\,  \textrm{dg}\ \BV^{\ac}\textrm{-coalgebras}\ : \ \Omega_\kappa.$$

Explicitly, the bar construction of a dg BV-algebra $A$ is equal to
$$\mathrm{B}_\kappa A:=\BV^{\ac}(A)$$ and the differential is given by the isomorphism of chain complexes $$\BV^{\ac}(A)=(\BV^{\ac}\circ_\kappa \BV)\circ_{\BV} A,$$
see Appendix~\ref{Homotopy Th General} for more details. Dually, the cobar construction of a $\BV^{\ac}$-coalgebra $C$ is equal to $\Omega_\kappa C:=\BV(C)$.

Since the operad $\BV$ is Koszul in the quadratic-linear sense given in Appendix~\ref{Appendix Koszul}, this adjunction has the following properties.

\begin{prop}\label{BV bar-cobar}
For any dg BV-algebra $A$, the counit of the bar-cobar adjunction is a quasi-isomorphism of dg BV-algebras:
$\Omega_\kappa \mathrm{B} _\kappa A \qi A$.
\end{prop}

\begin{proof}
It is an application of Proposition~\ref{Bar-cobar resolution}.
\end{proof}

Since we work over a field of characteristic $0$, every $\Sy$-module is exact \cite{GetzlerJones94}, $\Sigma$-split \cite{Hinich97} or $\Sigma$-cofibrant \cite{BergerMoerdijk03}. Therefore the category of dg BV-algebras  admits a model category structure in which the weak equivalences are quasi-isomorphisms. In this language, the previous proposition shows that the bar-cobar construction provides a cofibrant replacement functor for dg BV-algebras.

\subsection{Bar construction for homotopy BV-algebras and $\infty$-morphisms} 

We saw in Theorem~\ref{4 def theo} that a homotopy BV-algebra structure  is equivalently defined by a  square zero coderivation of the cofree $\BV^{\ac}$-coalgebra $\BV^{\ac}(A)$. This last construction is called the \emph{bar construction of the homotopy BV-algebra $A$} and denoted $\mathrm{B}_\iota A:=\BV^{\ac}(A)$.

Any dg BV-algebra $A$ is particular example of a homotopy BV-algebra on which both bar constructions coincide: $\mathrm{B}_\kappa A \cong \mathrm{B}_\iota A$. \\

Let $A$ and $B$ be two homotopy BV-algebras and let $\BV^{\ac}(A)$ and $\BV^{\ac}(B)$ be their bar constructions. An \emph{$\infty$-morphism} between $A$ and $B$ is a morphism of dg $\BV^{\ac}$-coalgebras $\BV^{\ac}(A) \to \BV^{\ac}(B)$. It is equivalent to a map $\BV^{\ac}(A) \to B$ satisfying a certain relation. An $\infty$-morphism is called an \emph{$\infty$-quasi-isomorphism} if the  first component  $I(A)\cong A\to B$ is a quasi-isomorphism. Since $\infty$-morphisms can be composed, the homotopy BV-algebras with $\infty$-morphisms form a category, which is denoted $\BV_\infty\textrm{-alg}$.\\

The aforementioned functors form the following commutative diagram.

$$\xymatrix@R=25pt@C=40pt{\textrm{dg}\ \BV\textrm{-alg}\  \ar@{^{(}->}[d]  \ar@_{->}@<-0.5ex>[r]_{\mathrm{B}_\kappa}& \ar@<-0.5ex>@_{->}[l]_{\Omega_\kappa} \  \textrm{dg}\ \BV^{\ac}\textrm{-coalg} \\
\BV_\infty\textrm{-alg}\ar[ur]_{\mathrm{B}_\iota} &} $$

\begin{prop}[Rectification]
For any homotopy BV-algebra $A$, there is an $\infty$-quasi-isomorphism of homotopy BV-algebras
$$A \xrightarrow{\sim} \Omega_\kappa \mathrm{B}_\iota A,$$
where the right hand side is a dg BV-algebra.
\end{prop}

\begin{proof}
This is direct application of the general Proposition~\ref{rectification GENERAL} to the Koszul operad $\BV$.
\end{proof}

\subsection{Transfer of structure and Massey products for homotopy BV-algebras}\label{Transfer BV}

In this section, we apply the general results of Sections \ref{Transfer} and \ref{Massey products} to the case of homotopy BV-algebras.

Let $(V, d_V)$ and $(W, d_W)$ be two homotopy equivalent chain complexes:
\begin{eqnarray*}
&\xymatrix{     *{ \quad \ \  \quad (V, d_V)\ } \ar@(dl,ul)[]^{h'}\ \ar@<1ex>[r]^{i} & *{\
(W,d_W)\quad \ \  \ \quad } \ar@(dr,ur)[]_h \ar@<1ex>[l]^{p}}&\\
&
\Id_V-p  i =d_V  h'
+ h'  d_V,
\quad
\Id_W-i p =d_W  h
+ h  d_W.
&
\end{eqnarray*}

\begin{theo}[Transfer Theorem for homotopy BV-algebras]\label{TransferThmBV}
Let $i\, :\, (V, d_V) \to (W, d_W)$ be a homotopy equivalence of chain complexes. Any homotopy BV-algebra structure on $W$ induces a homotopy BV-algebra structure on $V$ such that $i$ extends to an $\infty$-quasi-isomorphism.
\end{theo}

\begin{theo}[Massey products for homotopy BV-algebras]\label{MasseyProdBV}
Let $A$ be a homotopy BV-algebra. There is a homotopy BV-algebra structure on the homology $H(A)$ of the underlying chain complex of $A$, which extends its BV-algebra structure and such that the embedding $i\,:\, H(A)\mono A $  extends to an $\infty$-quasi-isomorphism of homotopy BV-algebras.
\end{theo}


\begin{prop}
For any homotopy BV-algebra $A$, one can reconstruct the homotopy type of $A$ from the homotopy BV-algebra structure of $H(A)$.
\end{prop}

Applied to TCFT, these homotopical results show the following theorem.

\begin{theo}
The homology groups $H(X)$ of any TCFT $X$ carry a homotopy BV-algebra structure, which extends the $BV$-structure of \cite{Getzler94} and which allows to reconstruct the homotopy type of $X$.
\end{theo}

After the results of Section~\ref{BVinfty and TCFT}, this theorem gives another way of lifting the BV-algebra structures of \cite{Getzler94}.

\subsection{Homology and cohomology of BV-algebras and homotopy BV-algebras}

After \cite[Section~$4.3$]{GetzlerJones94}, we know that the Andr\'e--Quillen homology of an algebra  over an operad with trivial coefficients is given by the left derived functor of the functor of indecomposable elements. Here we can use the functorial cofibrant resolutions $\Omega_\kappa \mathrm{B}_\kappa A \qi A$ and $\Omega_\iota \mathrm{B}_\iota A\qi A$ to compute Andr\'e--Quillen homology of BV-algebras and homotopy BV-algebras respectively:
\begin{eqnarray*}
\mathrm{H}_\bullet^{\BV}(A)=\mathrm{H}_{\bullet -1}(\overline{\mathrm{B}}_\kappa A)= \mathrm{H}_{\bullet -1} ( {\overline\BV}^{\ac}(A), d_\kappa).\\
\mathrm{H}_\bullet^{\BV_\infty}(A)=\mathrm{H}_{\bullet -1}(\overline{\mathrm{B}}_\iota A)= \mathrm{H}_{\bullet -1} ( {\overline\BV}^{\ac}(A), d_\iota).
\end{eqnarray*}

More generally, one can define the Andr\'e--Quillen homology and cohomology of BV-algebras and homotopy BV-algebras with coefficients, following \cite{GoerssHopkins00,Milles08}. The functorial cofibrant resolutions $\Omega_\kappa \mathrm{B}_\kappa A \qi A$ and $\Omega_\iota \mathrm{B}_\iota A\qi A$ provide explicit chain complexes which compute these homology and cohomology theories, as explained in \cite{Milles08}.




\appendix

\section{Koszul duality theory}\label{Appendix Koszul}

The theory of Koszul duality was originally developed  by Priddy in
\cite{Priddy70} for filtered associative algebras with relations
containing only quadratic and linear terms. The main examples are
the Steenrod algebra and the universal enveloping algebra of a Lie
algebra. This allowed him to provide quasi-free hence projective
resolutions for modules over these algebras. Positselski later
extended this theory to filtered associative algebras whose
relations are written with quadratic, linear and constant terms in
\cite{Positselski93}, see also Polishchuk--Positselski
\cite{PolischukPositselski05}. As a corollary, he gets
quasi-free hence cofibrant resolutions for the algebra itself.
Theorem~$10.2$ of \cite{Brown59} shows that these two kinds of resolutions, projective modules and
cofibrant algebras, are obtained at once. Koszul duality theory for graded quadratic associative algebras
was later extended to graded quadratic operads by Ginzburg--Kapranov
\cite{GinzburgKapranov94} and Getzler--Jones \cite{GetzlerJones94}.\\

Recall that an operad is a monoid in the monoidal category $(\Sy\textrm{-Mod}, \circ, I)$ of $\Sy$-modules with the composition product $\circ$. An operad is an algebraic object which encodes operations with several inputs but one output. To model operations with several outputs, we introduced in \cite{Vallette07} the notion of \emph{properad} which is a monoid in the monoidal category $(\Sy\textrm{-biMod}, \boxtimes, I)$ of $\Sy$-bimodules. The monoidal product, denoted $\boxtimes$ here, is based on graphs which represent the composition scheme of operations with several inputs and several outputs. Since the monoidal category of $\Sy$-modules is a sub-monoidal category of that of
$\Sy$-bimodules, an operad is a particular example of properad, where all the operations have exactly one output.
Koszul duality theory was extended from graded quadratic operads to graded quadratic properads in \cite{Vallette07}. From now on, we assume that the reader is familiar with this theory.  \\

In this section, we further extend Koszul duality theory to
properads, hence to operads, when the relations contain quadratic
and linear terms. Moreover, it applies to properads which are not necessarily
concentrated in homological degree $0$ and which may contain unary operations. As a direct corollary, we
obtain a Poincar\'e--Birkhoff--Witt Theorem for Koszul properads, analogous to that for associative algebras in \cite{PolischukPositselski05, BravermanGaitsgory96}. From this theory, we  explain how
to make the two types of resolutions explicit: at the level of properads and at the level of algebras over a properad.\\

\subsection{Filtered properads with quadratic and linear relations}\label{filtered properad}

Let $\Po$ be a properad which admits a presentation of the form
$\Po=\F(V)/(R)$, where $(R)$ is the ideal generated
by
$R\subset \F(V)^{(1)}
\oplus \F(V)^{(2)}=V \oplus \F(V)^{(2)}
$. The superscript $(n)$ indicates
the weight of the elements. In the free properad $\F(V)$, the weight of an element
is defined by  the number of generating elements from $V$ used to write it. Equivalently, it is
equal to the number of vertices of the underlying graph representing an element.
This means that $\Po$ has a space of generators $V$ and quadratic and linear relations $R$. The space of relations is supposed to be homogeneous with respect to the homological degree. \\

Let $\mathrm{q}\, :\,  \, \F(V) \epi \F(V)^{(2)}$ be the projection
onto the quadratic part of the free properad. We denote by $\qR$ the
image under $\mathrm{q}$ of $R$, so $\qR \subset \F(V)^{(2)}$ is the
quadratic part of the relations of $\Po$. We consider the quadratic
properad
$$\qPo:=\F(V)/(\qR).$$
We further assume that $R\cap V=\lbrace 0\rbrace
$. If it is not the case, one just has  to reduce the space of
generators $V$. Under this assumption, there exists a morphism of
$\Sy$-bimodules $\varphi \, :\,  \qR \to V$ such that $R$ is the
graph of $\varphi$
$$R=\lbrace X - \varphi(X), \ X \in \qR     \rbrace. $$

The natural grading of the free properad $\F(V)$ by the number of
generators induces the following filtration on $\Po$
$$F_n:= \pi\big(\bigoplus_{k\leq n}
\F(V)^{(k)}\big),$$
 where $\pi$ denotes the natural
projection $\F(V) \epi \Po$. We denoted the associated graded operad
by $\gr \Po$. The relations $\qR$ hold in $\gr \Po$. Therefore, there
exists an epimorphism of graded properads
$$p \, : \, \qPo \epi \gr \Po.$$
It is always an isomorphism in weight $0$ and $1$. In weight $2$, $p$ is an isomorphism
if and only if
$\qR=\mathrm{q}(\{V \oplus \F(V)^{(2)}\}\cap(R))$. 
Therefore $p$  is not
an isomorphism in general. At the end of this appendix, we will prove  that $p$ is an isomorphism when
the properad $\Po$ is Koszul.

\subsection{Koszul dual dg coproperad}\label{Koszul dual dg coproperad}

Recall that the Koszul dual coproperad $\qPo^{\ac}$ of $\qPo$ is the
cofree coproperad generated by $sV$ with relators in $s^2\qR$ (see
Section~$2$ of \cite{Vallette08})
$$\qPo^{\ac}:=\Co(sV,s^2\qR),   $$
where $s$ denotes the homological degree shift by $+1$.
 It is a sub-coproperad of the cofree coproperad $\F^c(sV)$ on
 $sV$.

We associate to  $\varphi$ the following composite map
$$\qPo^{\ac} \epi s^2\qR \xrightarrow{s^{-1}}
s\,\qR\xrightarrow{s\varphi}
sV.$$
By Lemma~$15$ of \cite{MerkulovVallette08I}, there exists a unique
coderivation $\widetilde{d}_\varphi\, :\, \qPo^{\ac} \to \F^c(sV)$  which extends this map. The following lemma gives the conditions for $\widetilde{d}_\varphi$ to induce a square zero coderivation on
the Koszul dual coproperad $\qPo^{\ac}$.

We will write $R\otimes V$ for 
the subspace of $\F(V)$ of linear combinations of connected graphs with $2$ vertices, the lower one labelled by an element of $R$ and the upper one labelled by an element of $V$, denoted $R\boxtimes_{(1,1)}V$ in \cite{MerkulovVallette08I}. The subspace $V\otimes R$ is defined similarly.

\pagebreak[2]

\begin{lemm}\label{coderivation varphi}~\linebreak[0]
\begin{itemize}
\item The coderivation $\widetilde{d}_\varphi \, :\, \qPo^{\ac} \to \F^c(sV) $  restricts to a coderivation $d_\varphi$ on
the sub-coproperad $\qPo^{\ac}\subset \F^c(sV)$ if $\{R\otimes V + V \otimes R\}\cap \F(V)^{(2)}\subset \qR$.
\item The coderivation $d_\varphi$ satisfies
${d_\varphi}^2=0$ if $\{R\otimes V + V \otimes R\}\cap \F(V)^{(2)}\subset R\cap \F(V)^{(2)}$.
\end{itemize}
\end{lemm}

\begin{proof}
By the universal property which defines $\Co=\Co(sV,s^2\qR)$ (see
Appendix B of \cite{Vallette08}), it is enough to check that
$\widetilde{d}_\varphi(\Co^{(3)})\subset \Co^{(2)}=s^2\qR$. The space  $s^2\qR \ot sV$ is the  subspace of $\F^c(sV)^{(3)}$
consisting of graphs with $3$ vertices indexed by $sV$ whose lower part
with $2$ vertices lives in $s^2\qR$. Dually, we consider $sV \ot
s^2\qR$. With these notations, we have $\Co^{(3)}=s^2\qR \ot sV \cap
sV \ot s^2\qR$. Hence any element $Y \in \Co^{(3)}$ be can written $Y=\sum
s^2 X \ot s v = \sum s v'\ot s^2 X'$, with $X, X' \in \qR$ and $v, v'\in V$. The
explicit description of $\widetilde{d}_\varphi$ given by  Lemma~$22$
of \cite{MerkulovVallette08I} gives
\begin{eqnarray*}
\widetilde{d}_\varphi(Y) &=&\sum s \varphi(X) \ot s v - \sum s
v'\ot s \varphi(X') \\ &=&\sum  (s\varphi(X)-s^2 X) \ot s v + \sum s
v'\ot (s^2 X'-s \varphi(X')).
\end{eqnarray*}
We identify $\F(V)$ with $\F^c(V)$ and we consider the natural suspension map $\F^c(V)\to \F^c(sV)$. So $\widetilde{d}_\varphi(Y)$ lives in the image of
$\{R\otimes V + V \otimes R\}\cap \F(V)^{(2)}$ in $\F^c(sV)$. Hence the coderivation $\widetilde{d}_\varphi$
restricts to $\Co$ if $\{R\otimes V + V \otimes R\}\cap \F(V)^{(2)}\subset \qR$.\\

Since $d_\varphi$ is a coderivation, to show that ${d_\varphi}^2=0$,
it is enough to prove that the projection of ${d_\varphi}^2$ onto the
generators $sV$ of $\Co$ vanishes. Once again, by the universal property of $\Co$, we only have to check here that
${d_\varphi}^2(\Co^{(3)})=0$. By the previous argument, this condition is equivalent to
$$\{R\otimes V + V \otimes R\}\cap \F(V)^{(2)}\subset \varphi^{-1}(0)=R\cap \F(V)^{(2)}.$$
\end{proof}

Notice that the second condition implies the first one.

\begin{defi}[Koszul dual dg coproperad]
Under the conditions of Lemma~\ref{coderivation varphi}, the \emph{Koszul dual dg coproperad of $\Po$} is the dg coproperad
$$\Po^{\ac}:=(\Co(sV,s^2\qR), d_\varphi).$$
\end{defi}

\subsection{Koszul duality theory}\label{KD Theory Q+L}

We extend Koszul duality theory to the quadratic-linear case as follows.

\begin{defi}
The properad $\Po$ is called a \emph{Koszul properad} if it admits a
quadratic-linear presentation $\Po=\F(V)/(R)$ such that

\begin{enumerate}
\item $R\cap V=\lbrace 0\rbrace$,

\item $\{R\otimes V + V \otimes R\}\cap \F(V)^{(2)}\subset R\cap \F(V)^{(2)}$,

\item its associated quadratic properad $\qPo:=\F(V)/(\qR)$ is Koszul in the
classical sense.
\end{enumerate}
\end{defi}

Condition $(1)$ deals with the minimality of the space of generators. It implies that the space of relations can be written as the graph of a map $\varphi\, : \qR \to V$. Condition $(2)$ corresponds to the maximality of the space of relations. It states that $R$ should contain all the quadratic relations obtained  by composing relations in $R$ with one element above or below. It implies that $\varphi$ extends to a  codifferential $d_\varphi$ 
on the Koszul dual coproperad $\qPo^{\ac}$.\\

In the example of the operad $\BV$ given in Section~$1$, the presentation satisfies these three conditions. But the presentation without the derivation relation between the operator $\Delta$ and the Lie bracket $\langle\, \textrm{-} ,  \textrm{-}\,  \rangle$ does not satisfy condition $(2)$ since this relation is obtained by composing the inhomogeneous relation between $\langle\, \textrm{-} ,  \textrm{-}\,  \rangle$ and $\bullet$-$\Delta$ by $\Delta$ above and below. \\

The bar-cobar constructions of associative algebras and coalgebras was extended to properads in
Section~$4$ of \cite{Vallette07} (see also Section~$3$ of \cite{MerkulovVallette08I}). The cobar construction
is a functor which produces a quasi-free properad from a dg coproperad.
Recall that a quadratic properad $\Po$ is Koszul if and only if the
cobar construction of its Koszul dual coproperad is a resolution of
$\Po$, $\Omega\, \Po^{\ac} \qi \Po$. In the
quadratic-linear case, we consider the cobar construction
 of the Koszul dual \emph{dg} coproperad $\Po^{\ac}=(\qPo^{\ac},d_\varphi)$, which
  we denote by $\Po_\infty:=\Omega \Po^{\ac} $.
More precisely, the underlying $\Sy$-bimodule of the
cobar construction $\Omega \Po^{\ac} $  of the dg coproperad
$(\qPo^{\ac}, d_\varphi)$ is the free properad
$\F(s^{-1}\overline{\qPo}^{\ac})$ on the desuspension of the
augmentation coideal of $\qPo^{\ac}$. Its differential $d$ is the sum
of two terms $d_1$ and $d_2$, where $d_1$ is the unique derivation
which extends $d_\varphi$ and where $d_2$ is the unique derivation
which extends the partial coproduct of the coproperad $\qPo^{\ac}$.
Since $d_\varphi$ is a square zero coderivation of
the coproperad $\qPo^{\ac}$, we get $d^2=0$ by Proposition~$4.4$ of
\cite{Vallette07}.\\

In the Koszul case, the next theorem shows that algebras over
$\Po_\infty$ give the explicit notion of homotopy $\Po$-algebras.

\begin{theo}\label{MainTHM}
Let $\Po$ be a Koszul properad. The cobar construction of its Koszul
dual dg  coproperad $\Po^{\ac}$ is a resolution of $\Po$
$$ \Omega \Po^{\ac} \qi \Po.$$
\end{theo}

\begin{proof}
Let $\Co:=s^{-1}\overline{\qPo}^{\ac}$ be the
desuspension of the augmentation coideal of the coproperad
$\qPo^{\ac}$. So, the underlying $\Sy$-bimodule of
$\Omega \Po^{\ac}$ is $\F(\Co)$.

Let us consider the new ``homological'' degree induced by the weight
of elements of $\qPo^{\ac}$ minus $1$. We call this grading the
\emph{syzygy degree} since it corresponds to the syzygies of the quasi-free resolution of $\qPo$ and therefore of $\Po$. Hence, the syzygy degree of an element of
$\F(\Co)$ is equal to the sum of the weight of the elements which
label its vertices minus the number of vertices. Since the weight of
$\Co$ is greater than $1$, the syzygy degree on $\F(\Co)$ is
non-negative.

The internal differential $d_1$ of $\F(\Co)$ induced by
$d_\varphi$ and the differential $d_2$  induced by the partial
coproduct of the coproperad $\Co$  lower the syzygy
degree by $1$. So we have a well-defined non-negatively graded chain
complex.

We consider the filtration $F_r$ of $\Omega \Po^{\ac} =\F(\Co)$
based on the total weight, which is defined for any graph labelled by elements of $\Co$ by the sum of the
grading of the operations indexing its vertices. We denote this grading by $(r)$. The two components
of the differential map $d=d_1+d_2$ satisfy
$$ d_2 \, : F_r \to F_r \quad \textrm{and} \quad d_1 \, : F_r \to F_{r-1}.$$
The filtration $F_r$ is therefore stable under the differential map
$d$. Since it is bounded below and exhaustive, the associated spectral
sequence $E^\bullet_{rs}$ converges to the homology of
$\Omega \Po^{\ac} $ by the classical convergence theorem of spectral
sequences (see Theorem 5.5.1 of \cite{Weibel}). Hence, $F_r$ induces
a filtration $\mathrm{F}_r$ on the homology of $\Omega \Po^{\ac} $
such that
$$E^\infty_{rs}\cong \mathrm{F}_r(H_{r+s}(\Omega \Po^{\ac} ))/
\mathrm{F}_{r-1}(H_{r+s}(\Omega \Po^{\ac}))=:\gr^{(r)}(H_{r+s}(\Omega \Po^{\ac} )).
$$

The first term of this spectral sequence is equal to
$E^0_{rs}=\F(\Co)^{(r)}_{r+s}$, which is made of the elements of
syzygy degree equal to $r+s$ and grading equal to $(r)$. The
differential map $d^0$ is given by $d_2$. Since the properad $\qPo$
is Koszul, the spectral sequence collapses at page $1$ and it is equal
to $E^1_{rs}=\qPo^{(r)}$. Therefore, $E^1_{rs}$ is concentrated in
the line $r+s=0$: $E^1_{rs}\cong \qPo^{(r)}$, for $r+s=0$ and
$E^1_{rs}=0$, for $r+s\neq 0$.

In conclusion, the convergence theorem gives
$$ E^1_{r, -r}\cong \qPo^{(r)}\cong  E^\infty_{r, -r} \cong
\gr^{(r)}(H_0(\Omega \Po^{\ac} )), $$
$$ E^1_{r,s}\cong 0 \cong  E^\infty_{r,s} \cong
\gr^{(r)}(H_{r+s}(\Omega \Po^{\ac} )),\ \textrm{for}\
r+s\neq 0.$$

With the syzygy degree, we have $H_0(\Omega \Po^{\ac} )=\Po$. And
the quotient $\gr^{(r)}(H_0(\Omega\Po^{\ac}))$ is equal to
$\gr^{(r)}\Po$. Therefore, the morphism $\Omega \Po^{\ac} \qi \Po$ is a
quasi-isomorphism.\end{proof}

This spectral sequence proof gives the following Poincar\'e--Birkhoff--Witt result directly.

\begin{theo}[Poincar\'e--Birkhoff--Witt Theorem for the properad $\Po$] \label{PBW for P}
When $\Po$ is a Koszul properad, the natural epimorphism $\qPo \epi \gr \Po$ is an isomorphism of
bigraded properads, with respect to the weight grading and the
homological degree. Therefore, the following $\Sy$-bimodules, graded
by the homological degree, are isomorphic
$$\Po \cong \gr\Po \cong \qPo.$$
\end{theo}

\begin{proof}
It is a direct corollary of the previous proof.
\end{proof}

This theorem allows us to give the following equivalent definition of the notion of a Koszul properad.

\begin{prop}\label{Max Relations}
The properad $\Po$  is Koszul if and only if it admits a quadratic-linear presentation $\Po=\F(V)/(R)$ satisfying the following conditions

\begin{itemize}
\item[(1)] $R\cap V=\lbrace 0\rbrace$,

\item[$(2')$] $R=\{V \oplus \F(V)^{(2)}\}\cap(R)$,

\item[(3)] its associated quadratic properad $\qPo:=\F(V)/(\qR)$ is Koszul in the
classical sense.
\end{itemize}
\end{prop}

\begin{proof}
When a properad $\Po$ is Koszul, the Poincar\'e--Birkhoff--Witt isomorphism in weight~$2$ implies
$\qR=\mathrm{q}(\{V \oplus \F(V)^{(2)}\}\cap(R))$, which is equivalent to $R=\{V \oplus \F(V)^{(2)}\}\cap(R)$ under
condition~$(1)$, $R\cap V=\lbrace 0\rbrace$. In the other direction, condition~$(2')$ always implies condition~$(2)$. \end{proof}

In this equivalent definition, condition~$(2')$ states that the space of relations $R$ should be maximal among the
generating spaces for the ideal of $\Po$. Such a condition can be hard to check in practice. But this proposition shows that if one finds a quadratic and linear presentation of a properad, which satisfies $(1)$, $(2)$ and $(3)$, then the space of relations $R$ is maximal.

\begin{rema}
This presentation of Koszul duality for properads with quadratic and
linear relations includes and extends the classical case of
quadratic properads. If $\Po$ has a quadratic presentation, then the
map $\varphi$ is equal to $0$. Hence the Koszul dual dg  coproperad
has a trivial differential. It is therefore equal to the classical
Koszul dual coproperad. Since $\qPo=\Po$, the two definitions of a
Koszul properad coincide.
\end{rema}

Recall that a representation of a properad, of the form $\Po\to \End_A$, is called a \emph{$\Po$-gebra structure on $A$} after \cite{Serre93}. (We drop the article ``al" from the Arabic word ``al-jabr" to try and encompass the various notions of algebras, coalgebras, bialgebras, etc.)

\begin{defi}[Homotopy $\Po$-gebra]
When $\Po$ is a Koszul properad, we define a \emph{homotopy $\Po$-gebra} as a gebra over the quasi-free resolution $\Po_\infty:=\Omega \Po^{\ac}$.
\end{defi}

The notion of homotopy $\Po$-gebra can be defined as a gebra over a cofibrant resolution of $\Po$. When $\Po$ is Koszul, there is a canonical choice for such a cofibrant replacement, namely $\Omega \Po^{\ac} $.

\subsection{Resolutions of modules}\label{Resolutions of Modules}

Another application of Koszul duality theory for associative algebras is to provide quasi-free, thus projective, resolutions for modules over a Koszul algebra. There are two main applications:  in algebraic topology to compute the derived functors Tor and Ext, see \cite{Priddy70}, and in algebraic geometry to study syzygies, see \cite{Eisenbud05}. We extend these resolutions to the properadic case.  \\

For three $\Sy$-bimodules 
$M$, $N$ and $O$, let us denote by $(M; N)\boxtimes O$ the part of $(M\oplus N)\boxtimes O$ linear in $N$. It is given by two-levelled graphs with vertices of the first level labelled by elements of $M$, except for one labelled by an element of $N$, and with the vertices of the second level labelled by elements of $O$. Dually, we define $M \boxtimes (N;O)$ as the part of $M\boxtimes (N\oplus O)$ linear in $O$. Let $f\, :\, M \to M'$ and $g\, :\, N \to O$ be two morphisms of $\Sy$-bimodules. We denote by $f\boxtimes ' g \, :\, M\boxtimes N \to M'\boxtimes (N;O)$ the morphism of $\Sy$-bimodules where we apply $f$ to every element of $M$ and $g$ to every element of $N$ but only one each time.
For example, if $(N,d)$ is a dg  $\Sy$-bimodule then the obvious  extension of the differential  to $M\boxtimes N$ may be expressed precisely as
$$
M\boxtimes N\xrightarrow{\text{Id}_M\boxtimes'd}M\boxtimes (N;N)\mono M\boxtimes (N\oplus N)\to M\boxtimes N.
$$ 

Any twisting morphism $\alpha \in \Tw(\Co, \Qo)$ between a (coaugmented) dg coproperad $\Co$ and a dg properad $\Qo$ induces a unique square zero derivation $\bar d_\alpha$ on the free left $\Qo$-module $\Qo \boxtimes \Co$ on $\Co$ which extends the following composite
$$\Co \xrightarrow{{}_{(1)}\Delta} (I; \overline\Co)\boxtimes \Co \xrightarrow{(\Id_I;\alpha)\boxtimes \Id_\Co} (I; \Qo)\boxtimes \Co \to \Qo\boxtimes \Co  , $$
where ${}_{(1)}\Delta$ is the component of the decomposition map $\Delta\, :\, \Co \to \Co\boxtimes \Co$  on
$(I; \overline\Co)\boxtimes \Co$.

\begin{defi}[Twisted composite product]
The chain complex $\Qo\boxtimes_\alpha \Co:=(\Qo \boxtimes \Co, d_\alpha=d_{\Qo \boxtimes \Co}+\bar d_\alpha)$ is called the \emph{left twisted composite product of $\Qo$ and $\Co$}. Similarly, one can define a \emph{right twisted composite product $\Co\boxtimes_\alpha \Qo$} and a \emph{two-sided composite product  $\Qo \boxtimes_\alpha \Co\boxtimes_\alpha \Qo$}.
\end{defi}
We refer the reader to Sections $3.4$ and $7.3$ of \cite{Vallette07} for more details and complete formulas. \\

The morphism $\Po_\infty=\Omega(\Po^{\ac})\to \Po$ is equivalent to a twisting morphism $\kappa\, :\, \Po^{\ac}\to\Po$ according 
\cite[Proposition~$17$]{MerkulovVallette08I}, which is given by
$$
\kappa\, : \, \Po^{\ac}
\epi sV \xrightarrow{s^{-1}}V \mono \Po.
$$
This twisting morphism defines a twisted composite product $\Po\boxtimes_\kappa \Po^{\ac} \boxtimes_\kappa \Po$ called the \emph{Koszul complex}.

\begin{theo}
For any  Koszul properad $\Po$, the twisted composite product is quasi-isomorphic to $\Po$ as $\Po$-bimodule:
$$\Po\boxtimes_\kappa \Po^{\ac} \boxtimes_\kappa \Po
\qi \Po.$$
\end{theo}

\begin{proof}
Once again, we consider another ``homological degree'', the total weight of the elements of $\Po^{\ac}$ this time. The two parts $d_\varphi$ and $\bar d_\kappa$ of the differential map lower this degree by $1$. The natural filtration on $\Po$ and the weight grading on $\Po^{\ac}$ induce an exhaustive and bounded below filtration $F_r$ on $\Po\boxtimes_\kappa \Po^{\ac} \boxtimes_\kappa \Po$. The differential maps satisfy $\bar d_\kappa\, :\, F_r \to F_r$ and $d_\varphi \, :\, F_r \to F_{r-1}$. Therefore, $E^0$ is equal to $\gr \Po\boxtimes_{\bar{\kappa}} \Po^{\ac} \boxtimes_{\bar{\kappa}} \gr \Po$ where $\bar{\kappa}\, :\, \Po^{\ac}\to \gr\Po$ is the associated twisting morphism and where  $d^0=\bar d_{\bar{\kappa}}$. By Corollary~\ref{PBW for P}, $E^0$ is the twisted composite product $\qPo\boxtimes_{\widetilde{\kappa}} \qPo^{\ac} \boxtimes_{\widetilde{\kappa}} \qPo$ of the Koszul quadratic properad $\qPo$, with ${\widetilde{\kappa}}\, : \, \qPo^{\ac} \to \qPo$ the Koszul twisting morphism. Therefore, it is quasi-isomorphic to $\qPo$ as an $\Sy$-bimodule. So, the convergence theorem for spectral sequences implies
$$ E^1_{r, -r}\cong \qPo^{(r)}\cong  E^\infty_{r, -r} \cong
\gr^{(r)}\Po, $$
$$ E^1_{r,s}\cong 0 \cong  E^\infty_{r,s}\  \textrm{for}\
r+s\neq 0.$$
Finally, $\Po\boxtimes_\kappa \Po^{\ac} \boxtimes_\kappa \Po$ is quasi-isomorphic to $\gr \Po$, which is equal to $\Po$ as an $\Sy$-bimodule by the PBW theorem of Corollary~\ref{PBW for P}.
\end{proof}

\begin{defi}[Relative composite product]
We define the \emph{relative composite product $M \boxtimes_\Po N$} of a right dg $\Po$-module $M$, $\rho\, :\, M \boxtimes \Po \to M$, \ and a left dg $\Po$-module $N$, $\lambda\, :\, \Po\boxtimes N \to N$,  by the following coequalizer:
$$\xymatrix@C=40pt{M \boxtimes \Po \boxtimes N  \ar@<0.5ex>[r]^(0.55){\Id_M \boxtimes \lambda} \ar@<-0.5ex>[r]_(0.55){\rho \boxtimes \Id_N}
& M \boxtimes N \ar@{->>}[r]& M \boxtimes_\Po N.  }$$
\end{defi}

\begin{coro}\label{Module Resolution}
Let $\Po$ be a Koszul properad. For any right $\Po$-module $M$ and any left $\Po$-module $N$, we have a quasi-isomorphism
$$M \boxtimes_\kappa \Po^{\ac} \boxtimes_\kappa N:=M\boxtimes_\Po(\Po\boxtimes_\kappa \Po^{\ac} \boxtimes_\kappa \Po) \boxtimes_\Po N \qi M\boxtimes _\Po N. $$
\end{coro}

\begin{proof}
Once again, we use the  total weight  of the elements of $\Po^{\ac}$  for the ``homological degree''.
Let us first prove the result with $N=\Po$, that is
$M \boxtimes_\kappa \Po^{\ac} \boxtimes_\kappa \Po:=M\boxtimes_\Po(\Po\boxtimes_\kappa \Po^{\ac} \boxtimes_\kappa \Po) \qi M$. Since $M\boxtimes_\Po(\Po\boxtimes_\kappa \Po^{\ac} \boxtimes_\kappa \Po)$ is given by the following short exact sequence
$$0\to  M\boxtimes\Po  \boxtimes (\Po\boxtimes_\kappa \Po^{\ac} \boxtimes_\kappa \Po)
\to M\boxtimes (\Po\boxtimes_\kappa \Po^{\ac} \boxtimes_\kappa \Po)  \to
M\boxtimes_\Po(\Po\boxtimes_\kappa \Po^{\ac} \boxtimes_\kappa \Po)\to 0,$$
it induces a long exact sequence in homology. By Proposition~$3.5$ of \cite{Vallette07}, the homology groups of the two first chain complexes vanish in degree greater than $1$. Therefore we have $H_n(M\boxtimes_\Po(\Po\boxtimes_\kappa \Po^{\ac} \boxtimes_\kappa \Po))=0$ for $n\ge 1$. In degree $0$, the long exact sequence is equal to
$ M\boxtimes \Po \boxtimes I \to M \boxtimes I \to M\boxtimes_\Po I$, so $H_0(M\boxtimes_\Po(\Po\boxtimes_\kappa \Po^{\ac} \boxtimes_\kappa \Po))=M\boxtimes_\Po I$. We apply the same method again to $M \boxtimes_\kappa \Po^{\ac} \boxtimes_\kappa N\cong (M\boxtimes_\Po(\Po\boxtimes_\kappa \Po^{\ac} \boxtimes_\kappa \Po))\boxtimes_\Po N$  to conclude the proof.

\end{proof}

For instance, with $M=\Po$ and $N=I$ (resp. with $M=I$ and $N=\Po$), this corollary proves that the Koszul complex $\Po \boxtimes_\kappa \Po^{\ac}$ (resp. $\Po^{\ac} \boxtimes_\kappa \Po$) is quasi-isomorphic to $I$. For $M=\Po$ and any $N$, this construction provides a quasi-free left $\Po$-module which is a resolution of  $N$.

When $\Po=R$ is concentrated in arity $(1,1)$, it is a Koszul algebra $R$, for instance the universal enveloping algebra of a Lie algebra, the Steenrod algebra and the free symmetric algebra \cite{Priddy70}. In this case, the above construction provides us with a functorial projective resolution of $R$-modules $M$: $R\otimes_\kappa R^{\ac}\otimes_\kappa M \qi M$. Notice that in algebraic geometry, people are looking for the \emph{minimal (quasi-)free} resolution of a module over an algebra, see \cite{Eisenbud05}. Condition~$(1)$ is the minimality condition and the resolution provided here is (quasi-)free.

Finally, when $A$ is a $\Po$-gebra, the construction
$$\Po\boxtimes_\kappa \Po^{\ac}(A):=(\Po\boxtimes_\kappa \Po^{\ac}\boxtimes_\kappa \Po)\boxtimes_\Po A \qi A$$
is a resolution of $A$ as a quasi-free left $\Po$-module. So, if $\Po$ is an operad, it gives a quasi-free $\Po$-algebra equivalent to $A$, which we analyse in the next appendix.

\section{Homotopy theory for algebras over a Koszul operad}\label{Homotopy Th General}

In this appendix, we suppose that $\Po$ is an \emph{operad}. (For a general properad $\Po$, as opposed to an  operad, the notion of free $\Po$-gebra does not exist and there is no model category structure on the category of $\Po$-gebras;  this category does not admit coproducts).

In this context, we define the \emph{bar and cobar constructions} for algebras over an inhomogeneous Koszul operad using the general theory of twisting morphism for (co)algebras over a (co)operad of \cite{GetzlerJones94}, see also \cite{Milles08} for a Lie theoretical approach. This allows us to define a \emph{bar construction} for homotopy $\Po$-algebras and a weaker notion of morphism, called \emph{$\infty$-morphism}. We show that any $\Po_\infty$-algebra is $\infty$-quasi-isomorphic to a $\Po$-algebra. We prove that $\Po_\infty$-algebra structures transfer through homotopy equivalences with explicit formulae. Applied to the homology of a $\Po_\infty$-algebra, this result defines the \emph{Massey products} for $\Po_\infty$-algebras.

\subsection{Bar and cobar constructions}\label{bar-cobar constructions GENERAL}

Let  $\alpha \, :\, \Co \to \Qo$  be a twisting morphism between a coaugmented dg cooperad $\Co$ and a dg operad $\Qo$. To any  dg $\Co$-coalgebra $(C, \Delta_C)$ and any dg $\Qo$-algebra $(A, \gamma_A)$, we associate the following  operator on $\Hom(C,A)$:
$$\star_\alpha (\varphi) \ : \ C \xrightarrow{\Delta_C} \Co \circ C
\xrightarrow{\alpha \circ \varphi} \Qo\circ A
\xrightarrow{\gamma_A} A.$$
A \emph{twisting morphism with respect
to $\alpha$} is a map $\varphi \, : \, C \to A$ of degree $0$ which is a
solution of the \emph{Maurer--Cartan}
equation:
$$\partial(\varphi) + \star_\alpha(\varphi)=0.$$ We denote the space of twisting
morphisms with respect to $\alpha$ by $\Tw_\alpha(C,
A)$. This bifunctor $\Tw_\alpha\, :\, \textrm{dg}\ \Co\textrm{-coalgebras}\times \textrm{dg}\ \Qo\textrm{-algebras}\to \textrm{Sets}$ is represented by the following functors.\\


To any dg $\Qo$-algebra $(A,\gamma_A)$, we define its \emph{bar construction}
 $\mathrm{B}_\alpha  A$   on the underlying  module $\Co(A)$. Since it is a cofree
$\Co$-coalgebra, there is a unique coderivation $d_2$ which
extends
$$\Co(A)=\Co\circ A \xrightarrow{\alpha\circ \Id_A} \Qo \circ A \xrightarrow{\gamma_A} A.$$
 The coderivation $d_2$  is equal
to the composite
\begin{eqnarray*}
\Co\circ A \xrightarrow{\Delta_{(1)}\circ \Id_A}
(\Co\circ (I; \overline \Co)) \circ A \xrightarrow{(\Id_{\Co} \circ
(\Id_I;\alpha))\circ \Id_A}
(\Co \circ (I; \Qo))\circ A  \cong
 \Co\circ(A; \Qo\circ A) \xrightarrow{\Id_{\Co} \circ(\Id_A; \gamma_A)} \Co \circ
A,
\end{eqnarray*}
where $\Delta_{(1)}$ is the component of the decomposition map $\Delta\, :\, \Co \to \Co\circ \Co$ on
$\Co \circ (I; \overline\Co)$. Hence $d_2$ squares to $0$. Finally, we endow $\Co(A)$ with the total square zero coderivation  $d:=d_\Co \circ \Id_A + \Id_\Co \circ' d_A + d_2$ to define $\mathrm{B}_\alpha A :=(\Co(A), d)$.
 The bar construction $\mathrm{B}_\alpha$ is a functor from dg $\Qo$-algebras to dg $\Co$-coalgebras such that $$\mathrm{B}_\alpha A =(\Co(A), d)\cong ((\Co\circ_\alpha \Qo)\circ_\Qo A, d_\alpha).$$ (The notations $\circ'$, $\circ_\alpha$ and $d_\alpha$ are the analogues for operads of $\boxtimes'$, $\boxtimes_\alpha$ and $d_\alpha$ of Appendix~\ref{Appendix Koszul}).


The \emph{cobar construction} $\Omega_\alpha C$ of a dg $\Co$-coalgebra $(C, \Delta_C)$ is defined dually. There is a unique derivation $d_2$ on $\Qo(C)$ which extends
$$C \xrightarrow{\Delta_C} \Co \circ C \xrightarrow{\alpha\circ \Id_C} \Qo \circ C.$$
This derivation is equal to the composite
\begin{eqnarray*}
\Qo\circ C \xrightarrow{\Id_\Qo \circ' \Delta_C} \Qo \circ
(C; \Co\circ C)\xrightarrow{\Id_\Qo \circ
(\Id_C;   \alpha\circ \Id_C)}  \Qo \circ
(C; \Qo\circ C) \cong (\Qo \circ (I, \Qo))(C)
\xrightarrow{\gamma_{(1)} \circ \Id_C} \Qo \circ C,
\end{eqnarray*}
so it squares to $0$.
We consider the total square zero derivation $d:=d_\Qo \circ \Id_C + \Id_\Qo \circ' d_C + d_2$, which defines the cobar construction of $C$, $\Omega_\alpha C :=(\Qo(C), d)$.
 The cobar construction $\Omega_\alpha$ is a functor from dg $\Co$-coalgebras to dg $\Qo$-algebras such that
$$\Omega_\alpha C =(\Qo(C), d)\cong ((\Qo\circ_\alpha \Co)\circ_\Co C, d_\alpha).$$

\begin{theo}[{\cite[Theorem~$2.18$]{GetzlerJones94}}]\label{bar-cobar alpha adjunction}
Let $\alpha \, : \, \Co \to \Qo$ be a twisting morphism from a
conilpotent dg cooperad $\Co$ to an operad $\Qo$. The bar and cobar constructions form a pair of adjoint functors
$$\mathrm{B}_\alpha \ : \ \textrm{dg} \ \Qo\textrm{-algebras}\,  \leftrightharpoons\,  \textrm{dg} \ \Co\textrm{-coalgebras}\ : \ \Omega_\alpha,$$
whose natural bijection satisfies
$$\Hom_{\emph{dg} \   \Qo\emph{-alg.}}\left(\Omega_\alpha C,\, A\right) \cong
\Tw_\alpha(C,\, A) \cong \Hom_{
\emph{dg} \ \Co\emph{-coalg.}}\left(C,\, \mathrm{B}_\alpha A\right),$$
for any dg $\Qo$-algebra $A$ and any dg $\Co$-coalgebra $C$.
\end{theo}

Let  $\Qo=\Po$ be an inhomogeneous quadratic operad and let $\Co=\Po^{\ac}$ be its Koszul dual dg cooperad. These results applied to the twisting morphism $\kappa \,:\, \Po^{\ac} \to \Po$ define the following pair of adjoint functors
$$\mathrm{B}_\kappa \ : \ \textrm{dg} \ \Po\textrm{-algebras}\,  \leftrightharpoons\,  \textrm{dg} \ \Po^{\ac}\textrm{-coalgebras}\ : \ \Omega_\kappa.$$

\begin{rema} When $\Po$ is a finitely generated binary quadratic operad, there is Lie theoretical interpretation of the aforementioned results. For any dg $\Po^{\ac}$-coalgebra $C$ and any dg $\Po$-algebra $A$, the space of maps $\Hom(C,A)$  can be endowed with a Lie bracket $[\; ,\, ]$ of degree $-1$, such that $\star_\alpha (\varphi) = \frac{1}{2}[\varphi, \varphi]$ \cite[Theorem~$2.1.1$]{Milles08}. In this case, the equation $\partial(\varphi) + \star_\alpha(\varphi)=0$ is the Maurer--Cartan equation in this dg Lie algebra. These properties still hold in the inhomogeneous quadratic case.
\end{rema}


\begin{lemm}\label{adjunction equivalence}
When the left twisted composite product $\Qo \circ_\alpha \Co$ is acyclic, then the counit of the adjunction is a quasi-isomorphism $\epsilon_A\, :\, \Omega_\alpha \mathrm{B}_\alpha A \xrightarrow{\sim} A$ for every dg $\Qo$-algebra $A$.

Dually, when the right twisted composite product $\Co \circ_\alpha \Qo$ is acyclic then the unit of the adjunction is a quasi-isomorphism $\upsilon_C \, :\, C \xrightarrow{\sim} \mathrm{B}_\alpha \Omega_\alpha C$ for every dg $\Co$-coalgebra $C$.

\end{lemm}

\begin{proof}
First observe that   $\Omega_\alpha \mathrm{B}_\alpha A\cong (\Qo \circ_\alpha \Co)(A)$ and $\mathrm{B}_\alpha \Omega_\alpha C\cong (\Co \circ_\alpha  \Qo)(C)$. We can then use standard filtrations and K\"unneth formula to conclude the proof.
\end{proof}

The following proposition is the extension of  Theorem~$2.25$ of \cite{GetzlerJones94} to the inhomogeneous case.

\begin{prop}\label{Bar-cobar resolution}
When $\Po$ is a Koszul operad, the counit of the bar-cobar adjunction is a quasi-isomorphism of dg $\Po$-algebras,
$$\Omega_\kappa \mathrm{B}_\kappa  A \qi A$$
and the unit of the adjunction is a quasi-isomorphism of dg $\Po^{\ac}$-coalgebras
$$C \xrightarrow{\sim}  \mathrm{B}_\kappa \Omega_\kappa C.$$
\end{prop}

\begin{proof}
It is direct Corollary of Lemma~\ref{adjunction equivalence} and Corollary~\ref{Module Resolution}.
\end{proof}

Since we work over a field of characteristic $0$, every $\Sy$-module is exact \cite{GetzlerJones94}, $\Sigma$-split \cite{Hinich97} or
$\Sigma$-cofibrant \cite{BergerMoerdijk03}. Therefore the category of $\Po$-algebras  admits a model category structure in which \emph{weak equivalences} are quasi-isomorphisms and \emph{fibrations} are epimorphisms. In this homotopical language, the weight of $\Po^{\ac}$ induces a suitable filtration on $\Po^{\ac}(A)$ which makes the bar-cobar construction into a cofibrant replacement functor for $\Po$-algebras.

\subsection{Bar construction for  $\Po_\infty$-algebras and $\infty$-morphisms}

In this section we apply Theorem~\ref{bar-cobar alpha adjunction} to a twisting morphism different from $\kappa$, in order to define a \emph{bar construction for homotopy $\Po$-algebras} and the notion of \emph{$\infty$-morphisms}.\\

Recall from Proposition~$18$ of \cite{MerkulovVallette08I} that any twisting morphism factors through two universal twisting morphisms. Applied to $\kappa$, it gives
$$\xymatrix{& {\Omega \Po^{\ac}} \ar@{-->}[dr]^{g_\kappa}& \\
\Po^{\ac} \ar@{-->}[rd]_{f_\kappa} \ar[ur]^{\iota} \ar[rr]^{\kappa} & & \Po \\
& {\mathrm{B} \Po}, \ar[ur]_{\pi} & }$$
where $g_\kappa\, :\, \Omega \Po^{\ac} \to \Po$ is a morphism of dg operads and $f_\kappa\, :\,  \Po^{\ac} \to \mathrm{B} \Po$ is a morphism of dg cooperads. Theorem~\ref{bar-cobar alpha adjunction} applied to the twisting morphism $\iota\, :\, \Po^{\ac} \to \Omega \Po^{\ac}=\Po_\infty$ gives the following  a pair of adjoint functors
$$\mathrm{B}_\iota \ : \ \Po_\infty\textrm{-algebras}\,  \leftrightharpoons\,  \textrm{dg} \ \Po^{\ac}\textrm{-coalgebras}\ : \ \Omega_\iota.$$
The bar construction $B_\iota$ is called the \emph{bar construction of a homotopy $\Po$-algebra $A$}. It is equal to
$$\mathrm{B}_\iota A=(\Po^{\ac}(A), d)=((\Po^{\ac}\circ_\gamma \End_A )\circ_{\End_A} A, d_\gamma),$$
where $\gamma \, :\, \Po^{\ac} \to \End_A$ denotes the twisting morphism defining the $\Po_\infty$-algebra structure on $A$.  This construction corresponds to the fourth equivalent definition in Theorem~\ref{4 def theo} of  a homotopy $\Po$-algebra structure on a dg module $A$ in terms of a square zero coderivation on the quasi-free dg $\Po^{\ac}$-coalgebra $\Po^{\ac}(A)$.

A homotopy $\Po$-algebra is a ``strict'' $\Po$-algebra if and only if the twisting morphism  $\gamma\, :\, \Po^{\ac} \to \End_A$  factors as $\Po^{\ac} \xrightarrow{\kappa} \Po \to \End_A$. In this case, both bar constructions agree since
$$\mathrm{B}_\iota A= (\Po^{\ac}\circ_\gamma \End_A )\circ_{\End_A} A \cong
((\Po^{\ac}\circ_\kappa \Po) \circ_{\Po}\End_A )\circ_{\End_A} A\cong (\Po^{\ac}\circ_\kappa \Po)\circ_{\Po} A=\mathrm{B}_\kappa A.$$

The various bar and cobar functors form the following commutative diagram.

$$\xymatrix@R=25pt@C=40pt{\textrm{dg}\ \Po\textrm{-alg.}\  \ar@{^{(}->}[d]  \ar@_{->}@<-0.5ex>[r]_{\mathrm{B}_\kappa}& \ar@<-0.5ex>@_{->}[l]_{\Omega_\kappa} \ \textrm{dg}\  \Po^{\ac}\textrm{-coalg.} \\
\Po_\infty\textrm{-alg.}\ar[ur]_{\mathrm{B}_\iota} &} $$

A morphism  between algebras over an operad is a map $f\, : A \to A'$ which strictly commutes with the operad action. When $\Po$ is a Koszul operad, we can define a weaker notion of morphisms as follows.
Let $A$ and $A'$ be two $\Po_\infty$-algebras.  An \emph{$\infty$-morphism} between $A$ and $A'$ is a morphism of dg $\Po^{\ac}$-coalgebras between the associated bar constructions $\mathrm{B}_\iota A=\Po^{\ac}(A) \to \mathrm{B}_\iota A' =\Po^{\ac}(A')$. It is equivalent to a map $\Po^{\ac}(A) \to A'$ satisfying a certain relation, see Chapter~$10$ of \cite{LodayVallette09} for an exhaustive study. An  $\infty$-morphism is called an \emph{$\infty$-quasi-isomorphism} if the  first component  $I(A)\cong A\qi A'$ is a quasi-isomorphism. By definition, $\infty$-morphisms can be composed. They share nice homotopy properties,  see for instance Section~\ref{Transfer}. From now on, we will only consider the category of $\Po_\infty$-algebras with their $\infty$-morphisms.\\

The following theorem was proved in \cite[Proposition~$3$]{DolgushevTamarkinTsygan07} for homogeneous Koszul operads.

\begin{prop}[Rectification]\label{rectification GENERAL}
Let $\Po$ be a Koszul operad. For any homotopy $\Po$-algebra $A$, there is an $\infty$-quasi-isomorphism of homotopy $\Po$-algebras
$$A \qi \Omega_\kappa \mathrm{B}_\iota A,$$
where the right hand side is a dg $\Po$-algebra.
\end{prop}

\begin{proof}
Proposition~\ref{Bar-cobar resolution} applied to the dg $\Po^{\ac}$-coalgebra $\mathrm{B}_\iota A$ provides a quasi-isomorphism of dg $\Po^{\ac}$-coalgebras $\mathrm{B}_\iota A \qi B_\kappa \Omega_\kappa \mathrm{B}_\iota A=\mathrm{B}_\iota \Omega_\kappa \mathrm{B}_\iota A$.
\end{proof}

\subsection{Transfer of homotopy structures}\label{Transfer}

Let $(V, d_V)$ and $(W, d_W)$ be two homotopy equivalent chain complexes:
\begin{eqnarray*}
&\xymatrix{     *{ \quad \ \  \quad (V, d_V)\ } \ar@(dl,ul)[]^{h'}\ \ar@<1ex>[r]^{i} & *{\
(W,d_W)\quad \ \  \ \quad } \ar@(dr,ur)[]_h \ar@<1ex>[l]^{p}}&\\
&
\Id_V-p  i =d_V  h'
+ h'  d_V,
\quad
\Id_W-i p =d_W  h
+ h  d_W.
&
\end{eqnarray*}

\begin{theo}[Transfer Theorem]\label{TransferThm}
Let $\Po$ be a Koszul operad and let $i\, :\, (V, d_V) \to (W, d_W)$ be a homotopy equivalence of chain complexes. Any $\Po_\infty$-algebra structure on $W$ induces a $\Po_\infty$-algebra structure on $V$ such that $i$ extends to an $\infty$-quasi-isomorphism.
\end{theo}

One can start with a $\Po$-algebra structure on $W$, but the transferred structure on $V$ will  be a $\Po_\infty$-algebra structure in general. This theorem provides a homotopy control of the transferred structure: the starting $\Po_\infty$-algebra and the transferred one are related by an explicit $\infty$-quasi-isomorphism. This result follows from the general principle that algebra structures over cofibrant operads transfer through quasi-isomorphisms, see \cite{BoardmanVogt73, BergerMoerdijk03}. We do not use model category arguments to prove this result here. Instead, we give explicit formulas. The proof relies on the following lemma and on the equivalent definitions of $\Po_\infty$-algebra structures given in Theorem~\ref{4 def theo}.

\begin{lemm}[{\cite[Theorem~$5.2$]{VanDerLaan03}}]\label{HomoMorphPsi}
 Let $V$ be a chain complex homotopy equivalent to a chain complex $W$. There is a morphism of coaugmented dg cooperads $\Psi \, :\, \mathrm{B} \End_W \to \mathrm{B} \End_V$, which extends
 $$\mu \in {\End}_W(n) \mapsto  p\, \mu\,  i^{\otimes n} \in {\End}_V(n).$$
\end{lemm}

Since the category of augmented (unital) dg operads is equivalent to the category of non-unital dg operads, we can equivalently apply the bar construction to non-unital dg operads. It is the case in the aforementioned lemma, where the underlying $\Sy$-module of $\mathrm{B} \End_W$ is $\F^c(s\End_W)$.

Since $\Psi$ is a morphism of cooperads to a cofree cooperad, it is completely characterized by its projection onto the cogenerators $\F^c(s\End_W)\to s\End_V$, which is defined as follows. A basis of $\F^c(s\End_W)$ is given by trees labelled by elements of $s\End_W$. Let $\TTT:=\TTT(s\mu_1,\ldots,s\mu_k)$ be such a tree. The image of $\TTT$ under $\Psi$ is defined by the suspension of the following composite: we label every leaf of the tree $\TTT(\mu_1,\ldots,\mu_k)$ with $i \, :\, V \to W$, every internal edge by $h$ and the root by $p$.

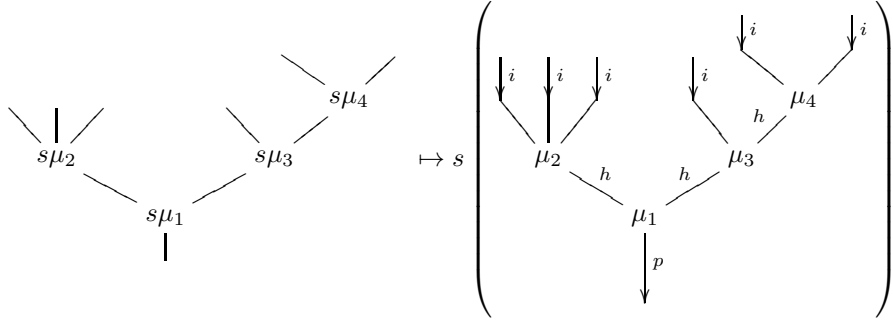
$$ \vcenter{\xymatrix@C=7pt@R=10pt{&&&&&&& \\ &&&&& \ar@{-}[rd]  && \ar@{-}[ld] \\
\ar@{-}[rd] & \ar@{-}[d]& \ar@{-}[ld] && \ar@{-}[rd]&& s\mu_4  \ar@{-}[ld] & \\
& s\mu_2  \ar@{-}[rrd]  &&&&  s\mu_3  \ar@{-}[lld] && \\
&&& s\mu_1  \ar@{-}[d] &&&& \\
&&&&&&&  \\&&&&&&&  }}
\mapsto
s\left(\vcenter{\xymatrix@C=7pt@R=10pt{ &&&&&\ar[d]^{i}&&\ar[d]^{i} \\ \ar[d]^{i} & \ar[d]^{i}& \ar[d]^{i}&& \ar[d]^{i}& *{}\ar@{-}[rd]  && *{}\ar@{-}[ld] \\
*{}\ar@{-}[rd] &*{} \ar@{-}[d]& *{}\ar@{-}[ld] && *{}\ar@{-}[rd]&& \mu_4  \ar@{-}[ld]_{h} & \\
& \mu_2  \ar@{-}[rrd]^{h}  &&&&  \mu_3  \ar@{-}[lld]_{h} && \\
&&& \mu_1  \ar[dd]^{p} &&&& \\
&&& *{}   &&&&  \\
&&&&&&&  }}\right)   $$

This composite scheme defines a map in $s\End_V$. Since the degree of $h$ is $+1$, the degree of $\Psi$ is $0$. Notice that this result is independent of the operad $\Po$.

\begin{proof}
[Proof of Transfer Theorem]
Let $\gamma \in \Tw(\Po^{\ac}, \End_W)$ be a $\Po_\infty$-algebra structure on $W$. By Theorem~\ref{4 def theo}, this twisting morphism is equivalent to a morphism of augmented dg cooperads $f_\gamma \, :\, \Po^{\ac}\to \mathrm{B}\End_W$. We compose it with the morphism of augmented dg cooperads $\Psi \, : \, \mathrm{B} \End_W \to \mathrm{B} \End_V$. The resulting composite $\Psi\,   f_\mu$ is a morphism of augmented dg cooperads which endows $V$ with a  $\Po_\infty$-algebra structure.

The formula for the extension of $i$ to a $\infty$-quasi-isomorphism between the transferred structure on $V$ and the $\Po_\infty$-algebra structure on $W$ is given by the same formula, except that we now label the root of the tree by the homotopy $h$ and not by $p$. The proof that this gives an $\infty$-morphism is straightforward. For more details, we refer the reader to Chapter~$10$ of \cite{LodayVallette09}.
\end{proof}

\subsection{Massey products}\label{Massey products}

Since the homology of the operad $\Po_\infty$ is equal to the operad $\Po$, the homology $H(A)$ of a $\Po_\infty$-algebra $A$ carries a natural $\Po$-algebra structure. By doing this, we lose a large amount of data, namely the homotopy type of $A$. We apply the preceding section to endow $H(A)$ with a $\Po_\infty$-algebra structure which extends this $\Po$-algebra structure. Moreover we show how to recover the homotopy type of $A$ from this data.\\

Let $(A, d)$ be a chain complex. As usual, we denote by $B_n:=\mathop{\rm Im }(d\, :\, A_{n+1} \to A_n)$ the image of the boundary map $d$. Since we work over a field $\KK$, each $A_n$ is isomorphic to $B_n \oplus H_n \oplus B_{n-1}$, after a choice of section. Under this isomorphism, the differential $d$ sends the direct summand $B_n \oplus H_n$ to $0$ and its restriction to  $B_{n-1}$ is an isomorphism. This splitting shows that the homology $H(A)$ is a deformation retract of $A$:
\begin{eqnarray*}
&\xymatrix{  (H(A), 0)\ \ar@{>->}@<1ex>[r]^{i} & *{\
(A,d)\quad \ \  \quad } \ar@(dr,ur)[]_h \ar@{->>}@<1ex>[l]^{p}
},&
\end{eqnarray*}
where the homotopy $h\, : \, A_{n} \to A_{n+1}$ is $0$ on $H_{n}\oplus B_{n-1}$ and is the inverse of $d$ on $B_n$.

\begin{theo}[Massey products]\label{MasseyProd}
Let $\Po$ be a Koszul operad and let $A$ be a $\Po_\infty$-algebra. There is a $\Po_\infty$-algebra structure on the homology $H(A)$ of the underlying chain complex of $A$, which extends its $\Po$-algebra structure and such that the embedding $i\,:\, H(A)\mono A $  extends to an $\infty$-quasi-isomorphism of $\Po_\infty$-algebras.
\end{theo}

\begin{proof}
The first two points follow directly from Theorem~\ref{TransferThm}. Because the morphism of coaugmented dg cooperads $\Psi$ extends $\mu \in {\End}_A(n) \mapsto  p\, \mu\,  i^{\otimes n} \in {\End}_{H(A)}(n)$, the image of ${\Po^{\ac}}^{(1)}$ in $\End_{H(A)}$, given by the transferred structure, corresponds to the $\Po$-algebra structure on $A$.
\end{proof}

One can also prove that the $\Po_\infty$-algebra structure on the homology $H(A)$ is independent of the choice of section for the homology of $A$: any two such transferred structures are related by an $\infty$-isomorphism whose first map is the identity on $H(A)$. We refer to Chapter~$10$ of \cite{LodayVallette09} for a proof of this fact. \\

When $A=C^\bullet_\textrm{sing}(X)$ is the singular cochain complex of a topological space, it is endowed with an associative cup product. This associative algebra structure transfers to an  $A_\infty$-algebra structure on the singular cohomology $H^\bullet_{\textrm{sing}}(X)$. These operations were originally defined by Massey in \cite{Massey58}, so  they are called the \emph{Massey products}. In general, we call the $\Po_\infty$-operations on the homology of a $\Po_\infty$-algebra the Massey products. \\

Whereas the differential on $H(A)$ is equal to $0$, the $\Po_\infty$-algebra structure on $H(A)$ is \emph{not} trivial. In this case, the relations satisfied by the $\Po_\infty$-algebra operations on $H(A)$ do not involve any differential. Hence the operations of weight $1$ satisfy the relations of a $\Po$-algebra. But the higher operations exist and contain the homotopy data of $A$. This result is well known in algebraic topology, where one uses the Massey product with three inputs to detect the non-trivial linking of the Borromean rings \cite{Massey69}.

\begin{prop}
For any Koszul operad $\Po$ and any $\Po_\infty$-algebra $A$, one can reconstruct the homotopy type of $A$ from the $\Po_\infty$-algebra structure of $H(A)$.
\end{prop}

\begin{proof}
Let us denote by $\widetilde \imath\, :\, H(A) \to A$ the extension of $i$ to an $\infty$-quasi-isomorphism. By definition of $\infty$-quasi-isomorphism and since $\Omega_\kappa$ preserves quasi-isomorphisms between quasi-cofree $\Po^{\ac}$-coalgebras (use the filtration based on the weight of $\Po^{\ac}$), it induces the following quasi-isomorphism of dg $\Po$-algebras
$$\Omega_\kappa \mathrm{B}_\iota H(A) \xrightarrow{\Omega_\kappa\mathrm{B}_\iota  \widetilde \imath} \Omega_\kappa \mathrm{B}_\iota A.$$
The latter dg $\Po$-algebra is $\infty$-quasi-isomorphic to $A$ by rectification (Theorem~\ref{rectification GENERAL}).
\end{proof}

\subsection{Homology and cohomology of $\Po$-algebras and $\Po_\infty$-algebras}

After \cite[Section~$4.3$]{GetzlerJones94}, we know that the Andr\'e--Quillen homology of an algebra  over an operad with trivial coefficients is given by the left derived functor of the functor of indecomposable elements. Here we can use the functorial cofibrant resolutions $\Omega_\kappa \mathrm{B}_\kappa A \qi A$ and $\Omega_\iota \mathrm{B}_\iota A\qi A$ to compute Andr\'e--Quillen homology of $\Po$-algebras and $\Po_\infty$-algebras respectively:
\begin{eqnarray*}
\mathrm{H}_\bullet^{\Po}(A)=\mathrm{H}_{\bullet -1}(\overline{\mathrm{B}}_\kappa A)= \mathrm{H}_{\bullet -1} ( {\overline\Po}^{\ac}(A), d_\kappa).\\
\mathrm{H}_\bullet^{\Po_\infty}(A)=\mathrm{H}_{\bullet -1}(\overline{\mathrm{B}}_\iota A)= \mathrm{H}_{\bullet -1} ( {\overline\Po}^{\ac}(A), d_\iota).
\end{eqnarray*}

More generally, one can define the Andr\'e--Quillen homology and cohomology of $\Po$-algebras and $\Po_\infty$-algebras with coefficients, following \cite{GoerssHopkins00,Milles08}. The functorial cofibrant resolutions $\Omega_\kappa \mathrm{B}_\kappa A \qi A$ and $\Omega_\iota \mathrm{B}_\iota A\qi A$ provide explicit chain complexes which compute these homology and cohomology theories, as explained in \cite{Milles08}.

\section{Deformation and obstruction theory}

In this third appendix, we develop the deformation theory and  the obstruction theory for gebras over a Koszul properad.

\subsection{Deformation Theory}\label{Deformation}

In this section, we give a Lie theoretic description of homotopy $\Po$-gebra structures when $\Po$ is a Koszul properad: we make precise the dg Lie algebra governing homotopy $\Po$-gebra structures. \\

Let $\Po$ be a Koszul properad. Recall that for any dg properad $\Qo$, the differential on the space of
$\Sy$-equivariant maps $\Hom_\Sy(\Po^{\ac}, \Qo)$ is defined by
$\partial(f):=d_\Qo\circ f - (-1)^{|f|} f\circ d_\varphi$.
It is shown in \cite{MerkulovVallette08I} that
$(\Hom_\Sy({\Po}^{\ac}, \Qo), \partial)$ is endowed with a  dg Lie algebra structure. We call it the \emph{convolution dg Lie algebra}, or simply the convolution algebra, and we denote it by $\g:=(\Hom_\Sy({\Po}^{\ac}, \Qo), \partial, [\,,\; ])$. Its Maurer--Cartan equation is
$$ \partial(\gamma)+ \frac{1}{2}[\gamma, \gamma]=0.$$
The solutions of the Maurer--Cartan equation in this convolution algebra which vanish on the coaugmentation map of the cooperad $\Po^{\ac}$ are called \emph{twisting morphisms} and the associated set is denoted by $\mathrm{Tw}(\Po^{\ac}, \Qo)$. They are in one-to-one correspondence with morphisms of dg properads from $\Po_\infty$ to $\Qo$ (Proposition~$17$ of \cite{MerkulovVallette08I}):
$$\mathrm{Tw}(\Po^{\ac}, \Qo)\cong \Hom_{\textrm{dg}\ \textrm{properads}}(\Po_\infty, \Qo).$$

Applied to the endomorphism properad, $\Qo=\End_A$, of a dg module $A$, it shows that structures of homotopy $\Po$-gebras on $A$ are in one-to-one correspondence with twisting morphisms of $\mathrm{Tw}(\Po^{\ac}, \End_A)$, cf. Theorem~\ref{4 def theo} in the operad case. \\

Since $\Po$ is Koszul, the underlying coproperad $\qPo^{\ac}$ of its Koszul dual dg coproperad $\Po^{\ac}=(\qPo^{\ac}, d_\varphi)$ is graded by an extra weight induced by that of the cofree coproperad. We denote this weight grading by ${\Po^{\ac}}^{(n)}$. With this convention, $\Po^{\ac}$ satisfies  ${\Po^{\ac}}^{(0)}=I$, ${\Po^{\ac}}^{(1)}=sV$ and ${\Po^{\ac}}^{(2)}=s^2\qR$. Therefore, the convolution Lie algebra $\g$ is also graded by this weight, $\g^{(n)}:=\Hom_\Sy({\Po^{\ac}}^{(n)}, \Qo)$. The differential on  the Koszul dual dg coproperad lowers this grading by one, $d_\varphi \, :\,  {\Po^{\ac}}^{(n)} \to {\Po^{\ac}}^{(n-1)}$. These properties of $\Po^{\ac}$ transfer to the convolution algebra as follows.

\begin{prop}\label{graded convolution Lie}
For any Koszul properad $\Po$ and any dg properad $\Qo$, the
underlying Lie algebra structure of the convolution dg  Lie algebra
$\g:=(\Hom_\Sy({\Po}^{\ac}, \Qo), \partial, [\;,\, ])$ has an
extra grading such that $\g=\prod_{n\ge 0} \g^{(n)}$.
Its differential
is the sum of two anti-commuting square-zero derivations $\partial=
\partial_0+\partial_1$ such that the first one preserves this grading
$\partial_0\, : \g^{(n)}\to \g^{(n)}$ and the second one raises it by
one, $\partial_1\, : \g^{(n-1)}\to \g^{(n)}$.
\end{prop}

\begin{proof}
We define $\g^{(n)}$ by $\g^{(n)}:=\Hom_\Sy({\Po^{\ac}}^{(n)}, \Qo)$ for $n\ge 0$. The differential  $\partial_0$ is given by the differentials $d_\Qo$ of $\Qo$ and the differential $\partial_1$ is given by the differential $d_\varphi$ on $\Po^{\ac}$. Since $d_\Qo$ is a square-zero derivation on the properad $\Qo$, $\partial_0$ is a square-zero derivation on $\g$ which preserves the weight. Since $d_\varphi$ is a square-zero coderivation on the coproperad $\Po^{\ac}$, $\partial_1$ is a square-zero derivation on $\g$ which raises the weight by one, $\partial_1 \, :\, \g^{(n-1)}\to  \g^{(n)}$.
\end{proof}

When $\Po$ is a quadratic Koszul properad the differential $\pa_1$ vanishes.

\subsection{Obstruction Theory}\label{Obstruction Theory}
 
We prove a general theorem about Maurer--Cartan elements in graded dg Lie algebras. We apply it to the convolution Lie algebras associated to Koszul properads. This defines an obstruction theory and a relative obstruction theory for homotopy $\Po$-gebras. \\

Let $(\g,\pa)$ be a dg Lie algebra whose underlying Lie algebra structure is isomorphic to $\g=\prod_{n\ge 0} \g^{(n)}$. We further assume that the differential $\pa$ is the sum of two square-zero derivations $\pa=\pa_0+\pa_1$ such that $\partial_0 \, :\, \g^{(n)}\to  \g^{(n)}$ and $\partial_1 \, :\, \g^{(n-1)}\to  \g^{(n)}$.

Any element of $\g$ can be written as a series  $\gamma = \gamma_0+\gamma_1+\cdots+\gamma_n+\cdots$.
 Under this decomposition, the Maurer--Cartan equation  is equivalent to
\begin{eqnarray*}
(\textrm{MC}_n) \ :\quad  \partial_0(\gamma_n) + \partial_1(\gamma_{n-1}) +  \frac{1}{2}\sum_{k+l=n} [\gamma_k, \gamma_l]=0
\end{eqnarray*}
in $\g^{(n)}$ for any $n\ge 1$ and to the equation $(\textrm{MC}_0)$:  $\partial_0(\gamma_0) +\frac{1}{2} [\gamma_0, \gamma_0]=0$, this last one being the Maurer--Cartan equation in the dg Lie algebra $(\g^{(0)}, \pa_0)$. Recall that a solution $\gamma$ to the Maurer--Cartan equation is required to have homological degree $-1$.

Notice that $\g$ is an extension of the Lie algebra $\g^{(0)}$ by $\prod_{n\ge 1} \g^{(n)}$.  Given a Maurer--Cartan element $\gamma_0$ in $(\g^{(0)}, \pa_0)$, we associate the twisted differential $\pa^{\gamma_0}:=\pa_0 + [\gamma_0, -]$ on $\g$.  It defines a square-zero derivation on $\g$, which preserves the weight $\g^{(n)}$.

\begin{theo}\label{Obstruction n->n+1:GENERAL}
Let $(\g=\prod_{n\ge 0} \g^{(n)}, \partial=\partial_0+\partial_1)$ be a  dg Lie algebra satisfying $\partial_0 \, :\, \g^{(n)}\to  \g^{(n)}$ and $\partial_1 \, :\, \g^{(n-1)}\to  \g^{(n)}$. Let $\gamma =  \gamma_0+\gamma_1 + \cdots + \gamma_n \in  \prod_{k=0}^n\g^{(k)}$ be an element which satisfies the $(\emph{MC})_k$-equations up to $k=n$ in $(\g, \partial)$. We consider
$$\widetilde{\gamma}_{n+1}:=\partial_1(\gamma_n) + \frac{1}{2}\sum_{k+l=n+1\atop k,l \ge 1}[\gamma_k, \gamma_l].$$

\begin{enumerate}
\item In $\g^{(n+1)}$, we have  $\partial^{\gamma_0}\left(  \widetilde{\gamma}_{n+1}  \right)=0$,
that is, $\widetilde{\gamma}_{n+1}$ is a cycle of degree $-2$ in the dg Lie algebra $(\g, \partial^{\gamma_0})$.
\item There exists an element $\gamma_{n+1}\in \g^{(n+1)}$ such that $\gamma_1 + \cdots + \gamma_{n+1}$ satisfies the $(\emph{MC}_k)$-equations up to $k=n+1$ in $(\g, \partial)$ if and only if
the class of $\widetilde{\gamma}_{n+1}$ in $\mathrm{H}_{-2}(\g^{(n+1)},  \partial^{\gamma_0})$ vanishes.
\end{enumerate}
\end{theo}

\begin{proof}
By the definition of $\widetilde{\gamma}_{n+1}$, we have
\begin{eqnarray*}
\partial^{\gamma_0}(\widetilde{\gamma}_{n+1})
&=& \pa^{\gamma_0}(\pa_1(\gamma_n)) +  \frac{1}{2}
\sum_{k+l=n+1\atop  k,l \ge 1 } \pa^{\gamma_0} (  [\gamma_k, \gamma_l]) .
\end{eqnarray*}

Since $\pa_0$ and $\pa_1$ anticommute, we get $\pa_0(\pa_1(\gamma_n))=-\pa_1(\pa_0(\gamma_n))$. Since $\pa^{\gamma_0}$ is a degree $-1$ derivation, we get
$$ \pa^{\gamma_0} (  [\gamma_k, \gamma_l])=    [\pa^{\gamma_0}(\gamma_k), \gamma_l] -   [\gamma_k, \pa^{\gamma_0}(\gamma_l)]=
 [\pa^{\gamma_0}(\gamma_k), \gamma_l] +   [\pa^{\gamma_0}(\gamma_l)   , \gamma_k ]. $$
And since $\pa_1$ is a degree $-1$ derivation, we have
$$ \pa_1 (  [\gamma_0, \gamma_n])=    [\pa_1(\gamma_0), \gamma_n] -   [\gamma_0, \pa_1(\gamma_n)].$$
Finally, we get
\begin{eqnarray*}
\partial^{\gamma_0}(\widetilde{\gamma}_{n+1}) =
-\pa_1(\pa^{\gamma_0}(\gamma_n)) + [\pa_1(\gamma_0), \gamma_n]+
\sum_{k+l=n+1\atop  k,l \ge 1 }   \left( [\pa^{\gamma_0}(\gamma_k), \gamma_l]
 \right).
\end{eqnarray*}
For any $n$, the equation $(\textrm{MC}_{n})$ is
$$\partial^{\gamma_0}(\gamma_{n})= - \pa_1(\gamma_{n-1}) -  \frac{1}{2}\sum_{k+l=n\atop k,l \ge 1 } [\gamma_k, \gamma_l]=- \widetilde{\gamma}_{n}.$$

So the induction hypothesis gives
\begin{eqnarray*}
\partial^{\gamma_0}(\widetilde{\gamma}_{n+1}) &=&
\frac{1}{2}\sum_{k+l=n\atop k,l \ge 1 } \pa_1[\gamma_k, \gamma_l]+[\pa_1(\gamma_0), \gamma_n]
-\sum_{k+l=n+1\atop k,l \ge 1} [\pa_1(\gamma_{k-1}), \gamma_l] -\frac{1}{2}\sum_{a+b+c=n+1 \atop a,b,c \ge 1} [[\gamma_a, \gamma_b], \gamma_c].
\end{eqnarray*}
Once again, since $\pa_1$ is a degree $-1$ derivation, we get
$$\partial_0(\widetilde{\gamma}_{n+1})= -\frac{1}{2}\sum_{a+b+c=n+1 \atop a,b,c \ge 1} [[\gamma_a, \gamma_b], \gamma_c], $$
which vanishes by Jacobi identity. \\

Since the equation $(\textrm{MC}_{n+1})$ is $\partial^{\gamma_0}(\gamma_{n+1})= -\widetilde{\gamma}_{n+1}$, it proves the second assertion.
\end{proof}
Therefore, the homology groups $\mathrm{H}_{-2}(\g^{(n)},  \partial^{\gamma_0})$, for $n\ge 1$, are the obstructions to the existence of Maurer--Cartan elements. Concretely, Theorem~\ref{Obstruction n->n+1:GENERAL} applies as follows.

\begin{theo}\label{Obstruction operations to full}
Let $\Po=\F(V)/(R)$ be a Koszul properad and let $(A, d_A)$ be a dg module. Suppose we are given a set of operations on $A$, $\gamma_1 \, :\, V \to\End_A$ such that $d_A$ is a derivation with respect to them, $\pa_A(\gamma_1)=d_{\End_A}\circ \gamma_1=0$. If
$\mathrm{H}_{-2}(\Hom_\Sy({\Po^{\ac}}^{(n)}, \End_A),  \partial_A)=0$ for $n\ge 2$, then $\gamma_1$ extends to a homotopy $\Po$-gebra structure on $A$.
\end{theo}

\begin{proof}
By Proposition~\ref{graded convolution Lie}, the convolution Lie algebra $\g:=(\Hom_\Sy({\Po}^{\ac}, \End_A)$ is graded and any twisting morphism satisfies $\gamma_0=0$. It remains to apply Theorem~\ref{Obstruction n->n+1:GENERAL} to this particular case where $\gamma_0=0$.
\end{proof}

This general method applies to the classical quadratic case as well. It was used in the particular case of homotopy Lie algebras in \cite{BarnichFulpLadaStasheff98} and in the case of homotopy Frobenius bialgebra structures on differential forms of a closed oriented  manifold in \cite{Wilson07}. We apply this extended version to prove the existence of a homotopy BV-algebras algebra structure of vertex algebras in Section~\ref{Vertex Algebras}.\\

\subsection{Relative Obstruction Theory} 

We study now a relative version of the previous method. Let $\Po=\F(V)/(R)$ be a Koszul properad composed of two properads in the following way. Suppose that the space of generators $V$ splits into two, $V=V_0\oplus V_1$, such that the space of relations splits into three $R=R_0\oplus R_1\oplus R_{01}$ with $R_0\subset \F(V_0)^{(2)}$, $R_1\subset V_1\oplus \F(V_1)^{(2)}$, $R_{01}\subset V_0\oplus \F(V_0\oplus V_1)^{(2)}$. (The subscripts $0$ and $1$ should not be confused with the homological degree). We further assume that $\qR_{01}$ is generated by elements given by a sum of 2-vertex graphs with exactly one vertex indexed by $V_0$ and the other indexed by $V_1$.  We denote by $\Po_0:=\F(V_0)/(R_0)$ the associated quadratic properad and by $\Po_1:=\F(V_1)/(R_1)$ the associated quadratic-linear properad. Therefore the properad $\Po$ is a quotient of the coproduct, or free product, $\Po_0\vee \Po_1$ by the ideal generated by $R_{01}$.
$$\Po\cong \frac{\Po_0\vee \Po_1}{(R_{01})}$$

In this case, the Koszul dual coproperad $\Po_0^{\ac}$ is a sub-coproperad of $\Po^{\ac}$. Moreover, the following diagram of graded coproperads commutes
$$\xymatrix@C=15pt@R=15pt@M=8pt{\Po_0^{\ac} \ar@{^{(}->}[r] \ar@{^{(}->}[d]& \Po^{\ac}\ar@{^{(}->}[d] \\
\F^c(sV_0)   \ar@{^{(}->}[r]& \F^c(sV_0\oplus sV_1)}$$
The assumption on $R_{01}$ proves that $\Po^{\ac}$ is a coproperad graded by the number of $sV_1$.

\begin{lemm}
Under the previous assumptions on the presentation of the properad $\Po$, the underlying coproperad of $\Po^{\ac}$  has an extra grading ${\Po^{\ac}}^{[n]}$, called the \emph{relative grading}, which satisfies ${\Po^{\ac}}^{[0]}=\Po_0^{\ac}$.
Its differential $d_\varphi$ lowers this relative grading by $1$, $d_\varphi \, :\, {\Po^{\ac}}^{[n]} \to {\Po^{\ac}}^{[n-1]}$.
\end{lemm}

\begin{proof}
Under the assumptions on the presentation of $\Po$, any map $\varphi$ associated to it satisfies $s^{-1}\varphi\, :\,
\F^c(sV_1)^{(2)} \to sV_1$ and $s^{-1}\varphi\, :\, \F^c(sV_0\oplus sV_1)^{(2)} \to sV_0$. Hence, it lowers the number of elements of $sV_1$ by $1$. And so does the differential $d_\varphi$.
\end{proof}

We get the following direct corollary.

\begin{prop}\label{graded relative convolution Lie}
For any Koszul properad $\Po$ satisfying the above assumptions and for any dg properad $\Qo$, the underlying Lie algebra structure of the convolution dg Lie algebra $\g:=(\Hom_\Sy({\Po}^{\ac}, \Qo), \partial, [\;,\, ])$ has an extra grading, $g=\prod_{n\ge 0} g^{[n]}$. Its differential is the sum of two anti-commuting square zero derivations $\partial= \partial_0+\partial_1$ such that the first one preserves this grading $\partial_0\, : \g^{[n]}\to \g^{[n]}$ and the second one raises it by one, $\partial_1\, : \g^{[n-1]}\to \g^{[n]}$.

Moreover, the sub dg Lie algebra $(\g^{[0]}, \pa_0)$ is equal to the convolution dg Lie algebra $\Hom_\Sy({\Po_0}^{\ac}, \Qo)$.
\end{prop}

\begin{proof}
The proof is similar to that of Proposition~\ref{graded convolution Lie} with the relative grading this time.
\end{proof}

Once again, any element of $\g$ is a series:  $\gamma = \gamma_0+\gamma_1+\cdots+\gamma_n+\cdots$.
And  the Maurer--Cartan equation  decomposes with respect to the weight as aforementioned.

\begin{theo}\label{Obstruction operations to full }
Let $\Po$ be a Koszul properad with a presentation satisfying the above assumptions and let $(A, d_A)$ be a dg module. Suppose we are given a homotopy $\Po_0$-gebra structure on $A$: $\gamma_0 \, :\, \Po_0^{\ac} \to\End_A$ such that $\pa_0(\gamma_0)+\frac{1}{2}[\gamma_0, \gamma_0]=0$. If $\mathrm{H}_{-2}(\Hom_\Sy({\Po^{\ac}}^{[n]}, \End_A),  \partial^{\gamma_0})=0$ for $n\ge 1$, then $\gamma_0$ extends to a homotopy $\Po$-gebra structure on $A$.
\end{theo}

\begin{proof}
It is a direct corollary of Proposition~\ref{graded relative convolution Lie} and Theorem~\ref{Obstruction n->n+1:GENERAL}.
\end{proof}

Conceptually, the coproperad $\Po^{\ac}$ is an extension of the coproperad $\Po_0^{\ac}={\Po^{\ac}}^{[0]}$ by $\bigoplus_{n\ge 1} {\Po^{\ac}}^{[n]}$. Hence, $\bigoplus_{n\ge 1} {\Po^{\ac}}^{[n]}$ is a $\Po^{\ac}_0$-comodule. Dually,
 the convolution Lie algebra $\g$ is an extension of $\g^{[0]}$  by $\bigoplus_{n\ge 1} \g^{[n]}$. So $\bigoplus_{n\ge 1} \g^{[n]}$ is a twisted dg module over the twisted dg Lie algebra $(\g^{[0]}, \pa^{\gamma_0})$. It naturally carries the obstructions to lift a Maurer--Cartan element of $\g^{[0]}$ to the whole of $\g$.

This theorem applies, for instance, when $\Po_0$ and $\Po_1$ are Koszul properads and when $\qR_{01}$ defines a distributive law, as in Theorem~\ref{qBV Koszul}. In this case, the Koszul dual coproperad has the following form
$\Po^{\ac}\cong \Po^{\ac}_1 \boxtimes \Po^{\ac}_0$. And the convolution Lie algebra is equal to
$\g\cong \Hom_\Sy(\Po_1^{\ac}\boxtimes \Po_0^{\ac}, \End_A)$ (see Section~\ref{Convolution Lie BV} and Section~\ref{Obstruction BV} for an application).


\begin{center}
\textsc{Acknowledgements}
\end{center}
The third author would like to warmly thank the University of
Barcelona and its Institut de Matem\`atica for the invitations and for the excellent working
conditions there. He is grateful to Ezra Getzler for many interesting
discussions, sometimes related to the subject of the article. Finally, we would like to thank the referees for the useful comments. 



\def\cprime{$'$}

\end{document}